\documentclass[final]{siamltex}
\usepackage{amsmath,amsfonts,amssymb}
\usepackage{latexsym}
\usepackage{color}
\usepackage[table, xcdraw, dvipsnames]{xcolor}
\usepackage[hidelinks]{hyperref}
\hypersetup{
    colorlinks=true,
    linkcolor={red!50!black},
    citecolor={blue!50!black},
    urlcolor={blue!80!black}
}
\usepackage{mathtools} 

\let\eps\varepsilon
\newcommand{\R}{\mathbb R}

\newcommand{\cA}{\mathcal A}
\newcommand{\hcA}{{\widehat{\mathcal A}}}

\newcommand{\cB}{\mathcal B}
\newcommand{\bC}{\mathbf C}

\newcommand{\bS}{\mathbf S}

\newcommand{\ba}{\mathbf a}

\newcommand{\bg}{\mathbf g}
\newcommand{\blf}{\mathbf f}
\newcommand{\bn}{\mathbf n}
\newcommand{\be}{\mathbf e}

\newcommand{\bu}{\mathbf u}

\newcommand{\bv}{\mathbf v}

\newcommand{\bx}{\mathbf x}

\newcommand{\bz}{\mathbf z}

\newcommand{\cZ}{\mathcal Z}

\newcommand{\balpha}{\mbox{\boldmath$\alpha$\unboldmath}}
\newcommand{\hbalpha}{\mbox{\boldmath$\widehat{\alpha}$\unboldmath}}
\newcommand{\bsigma}{\mbox{\boldmath$\sigma$\unboldmath}}
\newcommand{\bbeta}{\mbox{\boldmath$\beta$\unboldmath}}

\newcommand{\tsigma}{\widetilde{\sigma}}
\newcommand{\halpha}{\widehat{\alpha}}
\newcommand{\bhalpha}{\widehat{\boldsymbol{\alpha}}}

\newcommand{\bphi}{\boldsymbol{\phi}}
\newcommand{\bseta}{\boldsymbol{\eta}}
\newcommand{\bpsi}{\boldsymbol{\psi}}

\newcommand{\bPhi}{\boldsymbol{\Phi}}
\newcommand{\bPsi}{\boldsymbol{\Psi}}

\newcommand{\rA}{\mathrm A}
\newcommand{\rB}{\mathrm B}
\newcommand{\rC}{\mathrm C}

\newcommand{\rH}{\mathrm H}
\newcommand{\rI}{\mathrm I}
\newcommand{\rK}{\mathrm K}
\newcommand{\rM}{\mathrm M}
\newcommand{\rP}{\mathrm P}
\newcommand{\rQ}{\mathrm Q}
\newcommand{\rR}{\mathrm R}
\newcommand{\rS}{\mathrm S}
\newcommand{\rU}{\mathrm U}
\newcommand{\rV}{\mathrm V}
\newcommand{\rW}{\mathrm W}
\newcommand{\rZ}{\mathrm Z}

\newtheorem{algorithm}{Algorithm}

\def\rev#1{{#1}}

\begin{document}
\title{Interpolatory tensorial reduced order models for parametric dynamical systems}
\author{Alexander V. Mamonov\thanks{Department of Mathematics, University of Houston, 
Houston, Texas 77204 (avmamonov@uh.edu).} \and
Maxim A. Olshanskii\thanks{Department of Mathematics, University of Houston, 
Houston, Texas 77204 (maolshanskiy@uh.edu).}
}
\maketitle

\begin{abstract}
The paper introduces a reduced order model (ROM) for numerical integration of a dynamical system 
which depends on multiple parameters. The ROM is a projection of the dynamical system on a low dimensional
space that is both problem-dependent and \emph{parameter-specific}. 
The ROM exploits compressed tensor formats to find a low rank representation for a sample of high-fidelity 
snapshots of the system state. This tensorial representation provides ROM with an orthogonal basis in a universal 
space of all snapshots \emph{and} encodes information about the state variation in parameter domain. 
During the online phase and for any incoming parameter, this information is used to find  a  reduced basis that 
spans a parameter-specific subspace in the universal space. The computational cost of the online phase then 
depends only on  tensor compression ranks, but not on space or time resolution of high-fidelity computations.  
Moreover, certain compressed tensor formats enable to avoid the adverse effect of parameter space dimension 
on the online costs (known as the curse of dimension). 
The analysis of the approach includes an estimate for the representation power of the acquired ROM basis. 
We illustrate the performance and prediction properties of the ROM with several numerical experiments, 
where tensorial ROM's complexity and accuracy is compared to those of conventional POD-ROM. 
\end{abstract}


\section{Introduction}

In numerical optimal control, inverse modeling or  uncertainty quantification, one commonly needs to integrate 
a parameter-dependent dynamical system for various values of the parameter vector. 
For example, inverse modeling may require repeated solutions of the forward problem represented 
by a dynamical system, along the search path in a high-dimensional parameter space. 
This may lead to extreme-scale computations that, if implemented straightforwardly, often result in overwhelming 
computational costs. Reduced order models (ROMs) offer a possibility to alleviate these costs by replacing a 
high-fidelity model with a low-dimensional surrogate model~\cite{antoulas2000survey, gugercin2004survey}.
Thanks to this practical value and apparent success, ROMs for parametric dynamical systems have already 
attracted considerable attention; see, e.g., 
\cite{benner2015survey,hesthaven2016certified, bui2008model, baur2011interpolatory, benner2014robust,brunton2016discovering}.

In this paper, we are interested in projection based ROMs that build the surrogate model by projecting a 
high-fidelity model onto a low-dimensional problem-dependent vector space~\cite{benner2015survey}. 
Projection-based ROMs for dynamical systems include such well-known model order reduction techniques 
as proper orthogonal decomposition (POD) ROMs~\cite{lumley1967structure,sirovich1987turbulence} 
(and its variants such as POD-DEIM~\cite{chaturantabut2010nonlinear} and balanced POD~\cite{rowley2005model}) and
PGD-ROMs~\cite{chinesta2010recent,chinesta2013proper}.
In these approaches the basis for the projection space is computed by building on the information about 
the dynamical system provided through high-fidelity solutions sampled for certain time instances and/or 
parameter values, the so-called solution snapshots. However, building a general low-dimensional space for all 
times and parameters of interest might be challenging, if possible at all, for wide parameter ranges and long times. 
Several  studies aimed to address this challenge: 
In \cite{eftang2010hp,eftang2011hp,amsallem2012nonlinear} the authors considered partitioning strategies 
which introduce a subdivision of the parameter domain and assign an individual local reduced-order basis 
to each subdomain offline. Another idea~\cite{amsallem2008interpolation,son2013real} is to adapt 
precomputed equal-dimension reduced order spaces by interpolating them for out-of-sample parameters 
along geodesics on the Grassmann manifold. The present paper introduces a different approach that  
builds on recent developments in tensor decompositions and low-rank approximations to quickly compute 
parameter-specific reduced bases for projection based ROMs.

For parametrized systems of time-dependent differential equations, the generated data (the input of ROM) 
naturally takes a form of a multi-dimensional tensor of solution snapshots, with dimensionality $D+2$, 
where $D$ is a the dimension of the parameter space and 2 accounts for the spatial and time-wise distributions. 
Modern ROMs often proceed by unfolding such tensors into a matrix to perform standard POD based on truncated SVD. 
This leads to the loss of information about the dependency of solutions on parameters.
We propose to overcome these \rev{issues} by working directly with tensor data, and by exploiting low-rank tensor 
approximations based on the canonical polyadic (CP), high order SVD (HOSVD), 
and Tensor Train (TT) decompositions. The approach consists of two stages. First, at the offline stage,
the compressed snapshot tensor is computed using one of the three tensor decompositions, thus preserving
the essential information about variation of the solution with respect to parameters. Up to the compression 
accuracy, each of these decompositions provides a (global) basis for the universal space spanned by all 
observed snapshots. The so-called core of the compressed representation is then transmitted to the second stage, 
referred to as the online stage. At the online stage the transmitted part of compressed tensor
allows for a fast computation of a \emph{parameter-specific} reduced  basis for any incoming out-of-sample 
parameter vector through an interpolation and fast linear algebra routines. The reduced order basis is then 
given in terms of its coordinates in the global basis that can be stored offline. For CP and TT formats, 
the cost of these computations is free of exponential growth with respect to the parameter space dimension.
On analysis side of this work, we prove an estimate for prediction power of the parameter-specific reduced order 
basis. The estimate explicitly depends on the approximation accuracy of the original tensor by the compressed 
one, parameter interpolation error, and singular values of a small-size parameter-specific matrix.     

Despite an outstanding recent progress in numerical multi-linear algebra and, in particular, 
in understanding tensor decompositions (see, e.g., review articles
~\cite{ReviewTensor,grasedyck2013literature,sidiropoulos2017tensor}),
the application of tensor methods  in reduced order modeling of dynamical systems is still rather scarce.  
We mention two reports by Nouy \cite{nouy2015low,nouy2017low}, who reviewed tensor compressed formats 
and discussed their possible use for sparse function representation and reduced order modeling, 
as well as a series of publications on the treatment in compressed tensor formats of algebraic systems 
resulting from the stochastic and parametric  Galerkin finite element method, see e.g.
~\cite{benner2015low,benner2016low,benner2017solving,lee2019low, kressner2011low}.
 A POD-ROM was combined with a low-rank tensor representation of a mapping from a  parameter space onto an output domain in~\cite{kastian2020two}.
The authors of survey \cite{benner2015survey} observe that ``The combination of tensor calculus $\ldots$ 
and parametric model reduction techniques for time dependent problems is still in its infancy, but offers a promising 
research direction''. We believe the statement holds true, and the present study contributes to this largely open 
research field.

The remainder of the paper is organized as follows. In Section~\ref{s:setup} we set up a parameter-dependent 
Cauchy problem and recall the basics of POD-ROM approach that is needed for reference purpose later in the text.
Section~\ref{s:tROM} introduces a general idea of the interpolatory tensorial ROM and considers its realization 
using three popular tensor compression formats. Details are worked out for a Cartesian grid-based sampling of 
the parameter domain, and then the approach is extended to a more general parameter sampling scheme. 
A separate subsection discusses online--offline complexity and storage requirements of the method. 
An estimate on the prediction power of the reduced order basis is proved in Section~\ref{sec:analysis}. 
Numerical examples in Section~\ref{sec:num} illustrate the analysis and performance of the method. 
In particular, we compare the delivered accuracy with standard POD-ROM that employs a global low-dimensional basis.  

\section{Parameterized Cauchy problem and the conventional POD-ROM} \label{s:setup}

To fix ideas, consider the following multi-parameter initial value problem. 
For a vector of parameters $\balpha = (\alpha_1,\dots,\alpha_D)$ from the parameter domain 
$\cA \subset \mathbb{R}^D$ find the trajectory $\bu = \bu(t, \balpha) : [0,T) \to \mathbb{R}^M$ solving
\begin{equation}
\label{eqn:GenericPDE}
\bu_t = F(t, \bu, \balpha),  \quad t \in (0,T), \quad \text{and}~ \bu|_{t=0} = \bu_0,
\end{equation}
with a given continuous flow field  $F:(0,T)\times \mathbb{R}^M \times\cA\to \mathbb{R}^M$.
Hereafter we denote all vector quantities by bold lowercase letters.
We assume that the unique solution exists on $(0,T)$ for all $\balpha\in\cA$.
Examples considered in this paper include parameter-dependent parabolic equations,
in which case one can think of \eqref{eqn:GenericPDE} as a system of ODEs for nodal values of the finite 
volume or finite element solution to the PDE problem, where material coefficients, boundary conditions, 
or the computational domain (via a mapping into a reference domain) are parameterized by $\balpha$.

We are interested in projection based ROMs, where for an arbitrary but fixed $\balpha\in\cA$ an 
approximation to $\bu$ is sought as a solution to equations projected onto a reduced space. 
Projection based approaches aim at retaining the structure of the model and thus at preserving 
the physics present in the high-fidelity model~\cite{benner2015survey}.
Among the projection based approaches to model reduction for time-dependent differential equations, 
Proper Orthogonal Decomposition (POD) and its variants are likely the most widely used ROM technique, 
which provides tools to represent trajectories of a dynamical system in a low-dimensional, problem-dependent
basis~\cite{kerschen2005method,rathinam2003new,liang2002properi,liang2002properii}.
We summarize the POD-ROM below for further reference and for the purpose of comparison to our
approach in Section~\ref{sec:num}.

Assume for a moment that $\balpha $ is \emph{fixed}. The POD-ROM computes a representative collection 
of states $\bphi_k(\balpha) = \bu(t_k,\balpha)\in\mathbb{R}^M$ at times $0\le t_1,\dots,t_N<T$, 
referred to as snapshots, through high-fidelity numerical simulations. Next, one finds a parameter-specific 
low-dimensional basis $\{ \bz_i^{\rm pod}(\balpha) \}_{i=1}^{n} \subset \mathbb{R}^M$, $n \ll N$, 
referred to hereafter as the \emph{reduced basis}, such that the projection subspace 
$\mbox{span} \big\{ \bz_1^{\rm pod}(\balpha),\dots, \bz_n^{\rm pod}(\balpha)\big\}$ approximates 
the snapshot space $\mbox{span}\{\bphi_1(\balpha),\dots,\bphi_N(\balpha)\}$ in the best possible way. 

To determine the reduced basis, form a matrix of snapshots
\begin{equation}
\label{eqn:Phi}
\Phi_{\text{pod}}(\balpha) = [\bphi_1(\balpha), \ldots, \bphi_N(\balpha)] \in \R^{M\times N},
\end{equation}
compute its SVD
\begin{equation}
\label{eqn:SVDa}
\Phi_{\text{pod}}(\balpha) = \rU \Sigma \rV^T,
\end{equation}
and define $\bz_i^{\rm pod}(\balpha)$, $i=1,\dots,n$, to be the first $n$ left singular vectors of 
$\Phi_{\text{pod}}(\balpha)$, i.e., the first $n$ columns of $\rU$. Hereafter we denote all matrices 
with upright capital letters. The singular values in $\Sigma$ provide information 
about the approximation power of 
$\mbox{span} \big\{ \bz_1^{\rm pod}(\balpha),\dots, \bz_n^{\rm pod}(\balpha)\big\}$. 
We refer to \cite{liang2002proper,kerschen2005method} and references therein for a  discussion about 
algebraically different ways to define POD and their equivalence.

For parameters $\balpha$ varying in $\cA$, a parametric POD-ROM builds a global reduced basis 
by sampling the parameter domain, generating snapshots for each sampled parameter value and 
proceeding with SVD ~\eqref{eqn:SVDa} for a cumulative matrix of all snapshots. Possible sampling 
strategies include using a Cartesian grid in $\cA$, Monte--Carlo methods, and greedy algorithms based 
on a posteriori error estimates; see, e.g., \cite{benner2015survey,hesthaven2016certified}.
Regardless of the sampling procedure, the resulting basis can accurately reproduce only the data 
from which it originated. Without parameter-specificity, the basis may lack robustness for out-of-sample 
parameters, i.e., away from the reference simulations. 
This is a serious limitation for using POD based ROMs in inverse modeling. We plan to address this 
limitation by introducing tensorial techniques for finding reduced bases that are both problem- and 
parameter-specific.

\section{Tensorial ROMs}\label{s:tROM}

We first consider in Section~\ref{sec:grid} a Cartesian grid-based sampling of the parameter domain $\cA$ 
in the case when $\cA$ is the $D$-dimensional box 
\begin{equation}
\cA = \bigotimes_{i=1}^D [\alpha_i^{\min}, \alpha_i^{\max}],
\label{eqn:box}
\end{equation}  
and the sampling points are placed at the nodes of a Cartesian grid. Next, we describe three tensorial ROMs (\textbf{TROM}s) 
based on three different tensor decompositions, canonical polyadic (CP, Section~\ref{sec:CP-TROM}), high order SVD
(HOSVD, Section~\ref{sec:HOSVD-TROM}) and tensor train (TT, Section~\ref{sec:TT-TROM}).

\subsection{Cartesian grid-based parameter sampling}
\label{sec:grid}

To generate the sampling set $\hcA$, we distribute $n_i$ nodes $\{\halpha_i^j\}_{j=1,\dots,n_i}$ 
within each of the intervals $[\alpha_i^{\min}, \alpha_i^{\max}]$ in \eqref{eqn:box} for $i=1,\dots,D$, 
and define
\begin{equation}
\hcA = \left\{ \bhalpha =(\halpha_1,\dots,\halpha_D)^T\,:\,
\halpha_i \in \{\halpha_i^j\}_{j=1,\dots,n_i}, ~ i = 1,\dots,D \right\}.
\end{equation} 
Hereafter we use hats to denote parameters from the sampling set $\hcA$, and
the cardinality of $\hcA$ is denoted by 
\begin{equation}
K = \prod_{i=1}^{D} n_i.
\label{eqn:Acard}
\end{equation}
The corresponding snapshots
$\bphi_k(\bhalpha) = \bu(t_k, \bhalpha)$, $\bhalpha \in \hcA$, are organized in a multi-dimensional array
\begin{equation}
(\bPhi)_{:,i_1,\dots,i_D,k} = \bphi_k(\halpha_1^{i_1},\dots,\halpha_D^{i_D}),
\end{equation}
which is a tensor of order $D+2$ and size $M\times n_1\times\dots\times n_D\times N$. 
We reserve the first and the last indices of $\bPhi$ for the spatial and temporal distributions, respectively.
All tensors hereafter are denoted with bold uppercase letters.

Unfolding $\bPhi$ along the first index in an $M \times (n_1 \cdots n_D)N$ matrix and applying 
(truncated) SVD to determine the first $n$ left singular vectors is equivalent to the POD with 
grid-based parameter sampling. The disadvantage of this approach for ROM construction is that it 
neglects any information about the dependence of snapshots on parameter variation reflected in the 
tensor structure of $\bPhi$. To preserve this information, we proceed with a compressed approximation 
$\widetilde{\bPhi}$ of $\bPhi$  rather than with the low rank approximation of the unfolded matrix.

\subsection{Tensor compression and universal space}
The notion of a tensor rank and low-rank tensor approximation is ambiguous and later in this section 
we consider three popular compressed tensor formats. For now we only assume that $\widetilde{\bPhi}$ 
satisfies
\begin{equation}
\label{eqn:TensApprox}
\big\| \widetilde{\bPhi} - \bPhi \big\|_F \le \widetilde{\eps}\big\|\bPhi \big\|_F
\end{equation}
for some small $\widetilde{\eps} > 0$, where tensor Frobenius norm is simply
\begin{equation}
\| \bPhi \|_F := \Big( \sum_{j=1}^{M}
\sum_{i_1=1}^{n_1}\dots
\sum_{i_D=1}^{n_D}
\sum_{k=1}^{N} \bPhi_{j, i_1, \dots, i_D, k}^2 \Big)^{1/2}
\end{equation}
The "low-rank" (compressed) tensor $\widetilde{\bPhi}$ is computed during the first, \textit{offline} 
stage of TROM construction and \emph{a part of} $\widetilde{\bPhi}$ is passed on to the second, 
\textit{online} stage which uses this information about variation of snapshots with respect to changes 
in parameters to compute a parameter-specific TROM. 

We call \emph{universal space} the space $\widetilde{V}$ spanned by the first-mode fibers of 
$\widetilde{\bPhi}$, i.e., $\widetilde{V}$ is the column space of the mode-1 unfolding matrix. 
For the exact decomposition (i.e., for $\widetilde{\eps}=0$), $\widetilde{V}$ is the space of all observed 
system states. In general, $\widetilde{V}$ depends on a compression format, dimension of $\widetilde{V}$ 
does not exceed $M$ and depends on $\widetilde{\eps}$ and snapshot variation. We shall see that $\widetilde{V}$ 
does approximate the full space of high-fidelity snapshots, while the CP, HOSVD and TT formats all deliver 
an orthogonal basis for $\widetilde{V}$. In the online stage of TROM we find a local (parameter-specific) 
ROM basis by specifying its coordinates in $\widetilde{V}$.

\subsection{In-sample prediction}

During the online stage we wish to be able to approximately solve \eqref{eqn:GenericPDE} for an arbitrary 
parameter $\balpha \in \cA$ in an \emph{$\balpha$-specific} reduced basis. For this step we need  to introduce
the notion of $k$-mode tensor-vector product $\bPsi \times_k \ba$ of a tensor 
$\bPsi \in \R^{N_1 \times \dots \times N_m}$ of order $m$ and a vector $\ba \in \mathbb{R}^{N_k}$: 
the resulting tensor $\bPsi \times_k \ba $ has order $m-1$ and size 
$N_1 \times \dots \times N_{k-1} \times N_{k+1} \times \dots \times N_m$. 
Specifically, elementwise
\begin{equation}
(\bPsi \times_k \ba)_{j_1,\dots, j_{k-1},j_{k+1},\dots, j_{m}}=
\sum_{j_k = 1}^{N_k} \bPsi_{j_1, \dots, j_{m}} a_{j_k}.
\label{eqn:kmodeprod}
\end{equation}

For a moment, consider some  $\bhalpha = \left( \halpha_1,\dots,\halpha_D \right)^T$ from the 
sampling set $\hcA$ and define $D$ vectors, 
$\be^i (\bhalpha) = \left( e_{ 1}^i (\bhalpha), \dots, e_{ n_i}^i(\bhalpha) \right)^T \in \R^{n_i}$, 
$i=1,\dots,D$, as
\begin{equation}
\label{eqn:eij}
e^i_{ j}(\bhalpha) = 
\begin{cases} 1 & \text{if} ~ \halpha_i = \halpha_i^j \\ 0 & \text{otherwise} \end{cases},
\qquad j=1,\dots,n_i.
\end{equation}
In other words, $\be^i (\bhalpha)$ encodes the position of $\halpha_i$ among the grid nodes on 
$[\alpha_i^{\min}, \alpha_i^{\max}]$, $i=1,\ldots,D$. 

Vectors $\be^i(\bhalpha)$ defined above allow us to extract the snapshots corresponding to a particular 
$\bhalpha \in \hcA$. Specifically, we introduce the following \emph{extraction} operation
\begin{equation}
\Phi_e (\bhalpha) = 
\bPhi \times_2 \be^1(\bhalpha) \times_3 \be^2(\bhalpha) \dots \times_{D+1} \be^D(\bhalpha) 
\in \R^{M \times N},
\label{eqn:extracta}
\end{equation}
which extracts from tensor of all snapshots $\bPhi$ the matrix of snapshots \eqref{eqn:Phi} for the particular 
$\bhalpha \in \hcA$, i.e., $\Phi_e (\bhalpha) = \Phi_{\text{pod}}(\bhalpha)$.

Combining \eqref{eqn:extracta} with compressed approximation \eqref{eqn:TensApprox}, we conclude
that is should be possible to extract from $\widetilde{\bPhi}$ the information about the space spanned by 
the snapshots $\{ \bu(t_i, \bhalpha) \}_{i=1}^N$, for a  particular $\bhalpha \in \hcA$ up to the accuracy of 
approximation in \eqref{eqn:TensApprox}. Indeed, let $\bphi_i(\bhalpha) = \bu(t_i,\bhalpha)$, $i=1,\dots,N$, 
and denote by $\{ \bz_j(\bhalpha) \}_{j=1}^{\widetilde{N}}$, $\widetilde{N}\le N$, 
an orthonormal basis for the column space of
\begin{equation}
\widetilde{\Phi}_e (\bhalpha) = 
\widetilde{\bPhi} \times_2 \be^1(\bhalpha) \times_3 \be^2(\bhalpha) \dots \times_{D+1} \be^D (\bhalpha)
\in \R^{M \times N}, 
\label{eqn:extractat}
\end{equation}
where $\widetilde{N} = \mbox{rank} \big(\widetilde{\Phi}_e(\bhalpha) \big)$. Then, it holds
\begin{equation}
\label{eqn:phibound}
\sum_{i=1}^{N}\left\| \bphi_i - 
\sum_{j=1}^{\widetilde{N}}\langle \bphi_i, \bz_j \rangle \bz_j \right\|^2_{\ell^2} 
\le \widetilde{\eps}^2\big\|\bPhi \big\|_F^2,
\end{equation}
where we use the shortcut $\bphi_i=\bphi_i(\bhalpha)$, $\bz_i=\bz_i(\bhalpha)$.
To establish \eqref{eqn:phibound}, consider (thin) SVD 
$\widetilde{\Phi}_e (\bhalpha) = \widetilde{\rU}\widetilde{\Sigma}\widetilde{\rV}^T$ and compute
\begin{equation*}
\begin{split}
\sum_{i=1}^{N} \left\| \bphi_i - 
\sum_{j=1}^{\widetilde{N}}\langle \bphi_i, \bz_j \rangle \bz_j \right\|^2_{\ell^2} 
& = \left\| (\rI - \widetilde{\rU}\widetilde{\rU}^T){\Phi}_e (\bhalpha) \right\|^2_F \\
& = \left\| (\rI - \widetilde{\rU}\widetilde{\rU}^T) \left( {\Phi}_e (\bhalpha) - \widetilde{\Phi}_e (\bhalpha) \right) \right\|^2_F \\
& \le \left\| \rI - \widetilde{\rU}\widetilde{\rU}^T \right\|^2 
\left\| {\Phi}_e (\bhalpha) - \widetilde{\Phi}_e (\bhalpha) \right\|^2_F \le 
\left\| \Phi_e (\bhalpha) - \widetilde{\Phi}_e (\bhalpha) \right\|^2_F \\ 
& =  \left\| (\bPhi - \widetilde{\bPhi}) 
\times_2 \be^1(\bhalpha) \times_3 \be^2(\bhalpha) \dots \times_{D+1} \be^D(\bhalpha) \right\|^2_F \\
& \le \left\| \bPhi - \widetilde{\bPhi} \right\|^2_F \le \widetilde{\eps}^2\big\|\bPhi \big\|_F^2,
\end{split}
\end{equation*}
where we used linearity of \eqref{eqn:kmodeprod} and the inequality $\| \rA \rB \|_F \le \| \rA \| \| \rB \|_F$ 
for matrices $\rA$, $\rB$, and spectral matrix norm $\| \cdot \|$. We also used $\| \rP \| \le 1$ 
for an orthogonal projection matrix $\rP = \rI - \widetilde{\rU} \widetilde{\rU}^T$. 

Since $\bz_i(\bhalpha)\in \widetilde{V}$ for all in-sample $\bhalpha$, the bound in \eqref{eqn:phibound} 
provides an estimate on how accurate the true snapshots can be approximated in the universal space.
Furthermore, from \eqref{eqn:phibound} we conclude that given $\widetilde{\bPhi}$ and $\bhalpha \in \hcA$ we can 
obtain a parameter-specific (quasi)-optimal reduced basis by taking the first $n$ left singular vectors of 
$\widetilde{\Phi}_e (\bhalpha)$. Representation power of this basis is determined by $\widetilde{\eps}$ 
from \eqref{eqn:TensApprox} and $\sigma_i$, $i>n$, the singular values of $\widetilde{\Phi}_e (\bhalpha)$. 
If $ \widetilde{\eps}$ is sufficiently small, i.e., the snapshot tensor admits an efficient low-rank representation, 
then the computed reduced basis better represents the snapshot space \textit{for a given $\bhalpha$} 
than the first $n$ left singular vectors of the unfolded snapshot matrix (i.e., better \rev{than} the POD basis); 
see numerical results in Section~\ref{sec:num} which show up to several orders of accuracy gain 
for some of the examples. 

The arguments above apply only to parameter values $\bhalpha$ from the sampling set $\hcA \subset \cA$.
Next, we consider ROM basis computation for an arbitrary $\balpha\in \cA$ that may not necessarily come from the
training set $\hcA$, the so-called \emph{out-of-sample} $\balpha$. Below we explore the option \rev{of building} 
ROM basis for an out-of-sample $\balpha$ using interpolation in the parameter space. This approach is 
based on an assumption of smooth dependence of the solution $\bu(t, \balpha)$ of \eqref{eqn:GenericPDE} 
on $\balpha$. We refer to the corresponding tensorial ROMs as \textit{interpolatory} TROMs.

\subsection{Interpolatory TROM}

To construct parameter-specific ROM basis for an arbitrary $\balpha =(\alpha_1,\dots,\alpha_D)^T\in \cA$ 
we introduce the interpolation procedure defined by
\begin{equation}
\label{eqn:bea}
\be^i \,:\, \balpha \to \mathbb{R}^{n_i},\quad i=1,\dots,D.
\end{equation}
Entrywise, we write 
$\be^i(\balpha) = \big(e_{ 1}^i (\balpha), \ldots, e_{ n_i}^i (\balpha) \big)^T \in \R^{n_i}$.
The interpolation procedure should satisfy the following property. For a smooth function 
$\rev{g} : [\alpha_i^{\min}, \alpha_i^{\max}] \to \R$, it holds
\begin{equation}
\rev{g}(\alpha_i) \approx \sum_{j=1}^{n_i} e_{j}^i (\balpha) \rev{g}(\widehat{\alpha}_i^j),
\quad i=1,\ldots,D,
\label{eqn:interpapprox}
\end{equation}
where $\widehat{\alpha}_i^j$, $j=1,\ldots,n_i$, are the grid nodes on $[\alpha_i^{\min}, \alpha_i^{\max}]$.

We further consider Lagrangian interpolation of order $p$: for a given 
$\balpha \in \cA$ let $\widehat{\alpha}_i^{i_1}, \ldots, \widehat{\alpha}_i^{i_p}$ 
be the $p$ closest grid nodes to $\alpha_i$ on $[\alpha_i^{\min}, \alpha_i^{\max}]$, for $i=1,\ldots,D$.
Then 
\begin{equation}
\label{eqn:lagrange}
e_{j}^i (\balpha) = 
\begin{cases} 
\prod\limits_{\substack{m = 1, \\ m \neq k}}^{p}(\widehat{\alpha}_i^{i_m}-\alpha_i) \Big/ 
\prod\limits_{\substack{m = 1, \\ m \neq k}}^{p}(\widehat{\alpha}_i^{i_m}-\widehat{\alpha}_i^j), 
& \text{if } j = i_k \in \{i_1,\ldots,i_p\}, \\
\qquad\qquad\qquad\qquad\qquad\qquad\qquad 0, & \text{otherwise}, \end{cases}
\end{equation}
are the entries of $\be^i(\balpha)$ for $j=1,\dots,n_i$.
For the numerical experiments in 
Section~\ref{sec:num} we use $p=2$ or $3$, i.e., linear or quadratic interpolation. 
\rev{The Lagrangian interpolation is not the only
option and depending on parameter sampling and solution smoothness other fitting procedures can be more suitable.}

Vectors $\be^i$ extend the notion of position vectors $\be^i$ defined in \eqref{eqn:eij}
for out-of-sample vectors. Indeed, from \eqref{eqn:lagrange} it is easy to see that both 
vectors coincide if $\bhalpha \in \hcA$  and so we use the same notation hereafter.
Therefore, we can define \rev{a snapshot matrix $\widetilde{\Phi}_e (\balpha)$ through the 
extraction--interpolation procedure:}
\begin{equation}
\label{eqn:extractbt}
\widetilde{\Phi}_e (\balpha) = \widetilde{\bPhi} 
\times_2 \be^1(\balpha) \times_3 \be^2(\balpha) \dots \times_{D+1} \be^D(\balpha) 
\in \R^{M \times N},
\end{equation}
to generalize \eqref{eqn:extracta} for any $\balpha \in \cA$. 
Note  that \eqref{eqn:extractbt} and \eqref{eqn:extracta} are the same for
$\balpha = \bhalpha \in \hcA \subset \cA$, while \eqref{eqn:extractbt} \rev{defines} 
$\widetilde{\Phi}_e (\balpha)$ also for out-of-sample parameter vectors. 
\rev{If the low-rank representation of the snapshot tensor is exact, i.e., 
$\widetilde{\bPhi}=\bPhi$, then $\widetilde{\Phi}_e (\balpha)$ is the interpolation 
of the snapshot matrices ${\Phi}_{\rm pod} (\widehat\balpha)$.}

Once $\balpha \in \cA$ is fixed and $\widetilde{\Phi}_e (\balpha)$ is given by \eqref{eqn:extractbt}, 
our \emph{parameter-specific reduced basis} $\{\bz_i(\balpha)\}_{i=1}^n$ is defined as 
the first $n$ left singular vectors of $\widetilde{\Phi}_e (\balpha)$. Later we demonstrate that the coordinates 
of this basis in the universal space can be calculated quickly (i.e., using only low-dimensional calculations)
online without actually computing $\widetilde{\Phi}_e (\balpha)$.

In a \textit{non-interpolatory} TROM, a parameter-specific reduced basis can be 
constructed as follows. Choose $p \geq 2$ and fix $\balpha \in \cA$, then let 
$\widehat{\alpha}_i^{i_1}, \ldots, \widehat{\alpha}_i^{i_p}$ be the $p$ closest grid nodes to 
$\alpha_i$ on $[\alpha_i^{\min}, \alpha_i^{\max}]$, for $i=1,\ldots,D$, similarly to the interpolatory
construction above. Define the set 
\begin{equation}
\hcA_p := \left\{ \bhalpha = (\halpha_1,\dots,\halpha_D)^T\,:\,
\halpha_i \in \{\halpha_i^j\}_{j \in \{i_1, \ldots, i_p\}}, ~ i = 1,\dots,D \right\} \subset \hcA.
\end{equation}
Then, assemble a large matrix by concatenating $\widetilde{\Phi}_e (\bhalpha)$
for all $\bhalpha \in \hcA_p$ and take $\{\bz_i(\balpha)\}_{i=1}^{n}$ to be its first $n$ 
left singular vectors. Of course, hybrid strategies (e.g., interpolation only in some parameter directions) 
are also possible. For non-interpolatory or hybrid TROMs it is also possible to compute local basis online 
with only low-dimensional calculations following same steps as considered below.

In the rest of paper we focus on the interpolatory TROM and consider three well-known compressed formats 
for low rank tensor approximation $\widetilde{\bPhi} \approx \bPhi$: canonical polyadic (CP), 
Tucker, a.k.a. \rev{higher} order SVD (HOSVD), and tensor train (TT) decomposition formats. 
Note that the notion of tensor rank(s) differs among these formats. When applied to TROM computation, 
these formats lead to different offline computational costs to build $\widetilde{\bPhi}$, different amounts 
of information transmitted from the offline stage to the online stage (measured by the compression rate,
as explained in Section~\ref{sec:compcomp}), and slightly varying amounts of online computations for finding 
the reduced basis given an incoming $\balpha \in \cA$. 

\subsection{CP-TROM}
\label{sec:CP-TROM}

The first tensor decomposition that we consider is the canonical polyadic decomposition of a tensor 
into the sum of rank one tensors~\cite{hitchcock1927expression,carroll1970analysis,kiers2000towards,ReviewTensor}.
In CP-TROM we approximate  $\bPhi$ by the sum of $R$ (where $R$ is the so-called canonical tensor rank) 
direct products of $D+2$ vectors $\bu^r\in\R^M$, $\bsigma^{r}_i \in \R^{n_i}$, 
$i=1,\dots,D$, and $\bv^r \in \R^N$,
\begin{equation}
\label{eqn:CPv}
\bPhi \approx \widetilde{\bPhi} = 
\sum_{r=1}^{R} \bu^r \circ \bsigma^r_1 \circ \dots \circ \bsigma^r_D \circ \bv^r,
\end{equation}
or entry-wise
\begin{equation*}
(\widetilde{\bPhi})_{j,i_1,\dots,i_D,k}
=\sum_{r=1}^{R} u_j^r \sigma_{1, i_1}^r \dots \sigma_{D, i_D}^r v_{k}^r.
\end{equation*}

CP decomposition often delivers excellent compression. However, there are well-known difficulties 
in determining the accurate canonical rank $R$ and working with the CP format, see, e.g.,
\cite{haastad1990tensor,de2008tensor}. Since we are interested in the \textit{approximation} $\widetilde{\bPhi}$
to $\bPhi$,  the alternating least squares (ALS) algorithm~\cite{harshman1970foundations,ReviewTensor} 
can be used to minimize $\| \bPhi - \widetilde{\bPhi} \|_F$ for a specified target canonical rank $R$ to find the 
approximate factors $\bu^r\in\R^M$, $\bsigma^{r}_i\in\R^{n_i}$, and $\bv^r\in\R^N$, $r=1,\ldots,R$. 
We further assume $R\le M$, where $M$ is the dimension of high-fidelity snapshots.

\rev{Note that the second-mode product of a $D+2$-dimensional rank-one tensor $\bu^r \circ \bsigma^r_1 \circ \dots \circ \bsigma^r_D \circ \bv^r$ and a vector $\be^1(\balpha)\in \mathbb{R}^{n_1}$ is the $D+1$-dimensional rank-one tensor $\left\langle \bsigma^r_1, \be \right\rangle (\bu^r \circ \bsigma^r_2 \circ \dots \circ \bsigma^r_D \circ \bv^r)$. Proceeding with this computation for other modes, in the decomposition \eqref{eqn:CPv} we find that the definition  \eqref{eqn:extractbt}  yields representation of $\widetilde{\Phi}_e (\balpha)$  as the sum of rank one matrices for any  $\balpha \in \cA$:} 
\begin{equation}
\label{eqn:CPa}
\widetilde{\Phi}_e (\balpha) = \sum_{r=1}^{R} s_r\bu^r\circ\bv^r \in \R^{M \times N}, 
~~\text{with}~~
s_r = \prod_{i=1}^{D} \left\langle \bsigma^r_i, \be^i (\balpha) \right\rangle \in \R.
\end{equation}
However, \eqref{eqn:CPa} is \emph{not} the SVD of $\widetilde{\Phi}_e (\balpha)$, since  vectors $
\bu^r$ (and $\bv^r$) are not necessarily orthogonal. To avoid computing $\widetilde{\Phi}_e (\balpha)$ 
and its SVD online, the following preparatory offline step is required. 
Organize the vectors $\bu^r$ and $\bv^r$ from \eqref{eqn:CPv} into matrices
\begin{equation}
\widehat\rU  = [\bu^1, \dots, \bu^{R}] \in \mathbb{R}^{M \times R}, \quad 
\widehat\rV = [\bv^1, \dots, \bv^{R}] \in \mathbb{R}^{N \times R},
\label{eqn:rurv}
\end{equation}
and compute the thin QR factorizations 
\begin{equation}
\widehat\rU = {\rU} \rR_U, \quad \widehat\rV = {\rV} \rR_V, 
\label{eqn:thinqr}
\end{equation}
if $R>N$ let further $\rR_V=\widehat\rV$ and ${\rV}=\rI$. 
\emph{The columns of ${\rU}$ form an orthogonal basis in the universal space $\widetilde{V}$}. 
Matrix ${\rU}$ is stored offline (${\rV}$ is not used and can be dropped), while low-dimensional matrices 
$\rR_U$ and $\rR_V$ together with vectors $\bsigma^r_i$ form the \emph{online} part of $\widetilde{\bPhi}$,
\begin{equation}
\label{eqn:coreCP}
\mbox{online}(\widetilde{\bPhi}) = 
\left\{ \rR_U, \rR_V \in \mathbb{R}^{R\times R},~
\bsigma^r_i\in\mathbb{R}^{n_i},~{\small i=1,\dots,D} \right\},
\end{equation}
which is transmitted to the online stage.

At the online stage and for any incoming $\balpha \in \cA$, we compute the SVD of the 
$R\times R$ \emph{core} matrix 
\begin{equation}
\rC(\balpha) = \rR_U \rS(\balpha) \rR_V^T,
\label{eqn:cpcore}
\end{equation} 
where $\rS(\balpha)=\mbox{diag}(s_1,\dots,s_R)$, with $s_r$ from \eqref{eqn:CPa}:
$ \rC(\balpha) = \rU_c\Sigma_c \rV_c^T$. Since 
\begin{equation}
\widetilde{\Phi}_e (\balpha) = {\rU}\rC(\balpha){\rV}^T = 
\left({\rU}\rU_c\right)\Sigma_c \left({\rV}\rV_c\right)^T
\end{equation} 
is the SVD of $\widetilde{\Phi}_e (\balpha)$,  \emph{the first $n$ columns of $\rU_c$, denoted by 
$\left\{\bbeta_1(\balpha),\dots,\bbeta_n(\balpha)\right\}$, are the coordinates of the local reduced basis 
in the universal space $\widetilde{V}$.} The parameter-specific basis in the physical space is then 
$\{\bz_i(\balpha)\}_{i=1}^{n}$, with $\bz_i(\balpha) = {\rU}\bbeta_i(\balpha)$. 
Note that $\bz_i(\balpha)$ are not actually computed.

Under certain assumptions on $F$, the dynamical system~\eqref{eqn:GenericPDE} is projected offline 
onto $\widetilde{V}$ and passed to the online stage, where for any $\balpha \in \cA$ it is further projected 
onto the local basis $\left\{\bbeta_1(\balpha),\dots,\bbeta_n(\balpha)\right\}$. This avoids any online 
computations with high-dimensional objects used in high-fidelity simulations; 
see further discussion in Section~\ref{rem1}. 

We summarize the above in the following algorithm.

\vskip0.05in
\begin{algorithm}[CP-TROM]
\label{alg:CP-TROM}
\begin{itemize}
\item \textbf{Offline stage}.\\
\emph{Input:} snapshot tensor $\bPhi \in \R^{M \times n_1 \times \ldots \times n_D \times N}$, 
target canonical rank $R$;\\
\emph{Output:} CP decomposition factors, universal basis matrix ${\rU} \in \R^{M \times R}$, 
and upper triangular matrices $\rR_U$, $\rR_V \in \mathbb{R}^{R\times R}$;\\
\emph{Compute:}
\begin{enumerate}
\item Use ALS algorithm to minimize $\| \bPhi - \widetilde{\bPhi} \|_F$ to find
CP decomposition factors $\bu^r$, $\bsigma^{r}_i$, and $\bv^r$ of $\widetilde{\bPhi}$ satisfying \eqref{eqn:CPv};
\item Assemble matrices $\widehat\rV, \widehat\rU$ as in \eqref{eqn:rurv} and compute 
their thin QR factorizations \eqref{eqn:thinqr} to obtain $\rR_U$, $\rR_V$ and ${\rU}$.
\end{enumerate}
\item \textbf{Online stage}. \\
\emph{Input:} $\mbox{\rm online}(\widetilde{\Phi})$ as defined in \eqref{eqn:coreCP}, 
reduced space dimension $n \le R$, and parameter vector $\balpha \in \cA$; \\
\emph{Output:} Coordinates of the reduced basis in $\widetilde{V}$:
$\{ \bbeta_i(\balpha) \}_{i=1}^{n} \subset \mathbb{R}^R$;\\
\emph{Compute:}
\begin{enumerate}
\item Use \eqref{eqn:cpcore} to assemble the core matrix $\rC(\balpha)$;
\item Compute the SVD of the core matrix $\rC(\balpha)=\rU_c\Sigma_c \rV_c^T$, 
with $\rU_c = [\widetilde{\bu}_1, \widetilde{\bu}_2, \ldots, \widetilde{\bu}_R]$;
\item Set $\bbeta_i(\balpha) = \widetilde{\bu}_i$, $i=1,\dots,n$.
\end{enumerate}
\end{itemize}
\end{algorithm}
\vskip0.05in

Note that we do not have direct control over ALS algorithm to enforce a priori CP decomposition accuracy
$\| \bPhi - \widetilde{\bPhi} \|_F < \widetilde{\eps}\| \bPhi  \|_F$. One option is to rerun the offline stage 
\rev{for different} trial values of $R$. Given that the offline stage is computationally expensive, this may become
prohibitive in cases where a desired accuracy $\widetilde{\eps}$ must be strictly enforced. 
The other two variants of TROM presented below are free from this limitation.

\subsection{HOSVD-TROM}
\label{sec:HOSVD-TROM}

As we already mentioned, truncated variant of CP decomposition is not known to satisfy any simple minimization 
property (unlike the SVD decomposition for matrices). A classical tensor decomposition, known to deliver a 
(quasi)-minimization property, is the so-called higher order SVD ({HOSVD})~\cite{de2000multilinear}.  
In HOSVD-TROM variant we approximate the snapshot tensor with a 
Tucker tensor~\cite{tucker1966some,ReviewTensor} $\widetilde{\Phi}$ of the form
\begin{equation}
\label{eqn:TDv}
\bPhi \approx \widetilde{\bPhi} = 
\sum_{j = 1}^{\widetilde{M}}
\sum_{q_1 = 1}^{\widetilde{n}_1}\dots
\sum_{q_D = 1}^{\widetilde{n}_D}
\sum_{k = 1}^{\widetilde{N}}
(\bC)_{j, q_1, \dots, q_D, k} \bu^j \circ \bsigma^{q_1}_1 \circ \dots \circ \bsigma^{q_D}_D \circ \bv^k,
\end{equation}
with $\bu^j \in \R^M$, $\bsigma_i^{q_i} \in \R^{n_i}$, and $\bv^k \in \R^N$.
The numbers $\widetilde{M}$, $\widetilde{n}_1$, $\ldots$, $\widetilde{n}_D$ and $\widetilde{N}$ 
are referred to as Tucker ranks of $\widetilde{\bPhi}$. The HOSVD delivers an efficient compression 
of the snapshot tensor, provided the size of the \emph{core tensor} 
$\bC \in \R^{\widetilde{M} \times \widetilde{n}_1 \times \dots \times \widetilde{n}_D \times \widetilde{N}}$ 
is (much) smaller than the size of $\bPhi$.

In what follows, it is helpful to organize the column vectors from \eqref{eqn:TDv} into matrices
\begin{equation}
\begin{split}
\rU & = [\bu^1, \dots, \bu^{\widetilde{M}}] \in \mathbb{R}^{M \times \widetilde{M}}, \quad 
\rV = [\bv^1, \dots, \bv^{\widetilde{N}}] \in \mathbb{R}^{N \times \widetilde{N}}, \\
\rS_i & = [\bsigma^1_i, \dots, \bsigma^{\widetilde{n}_i}_i]^T \in \mathbb{R}^{\widetilde{n}_i \times {n}_i}, \quad
i = 1,\ldots,D.
\end{split}
\label{eqn:hosvdmat}
\end{equation}

In contrast with CP decomposition, HOSVD computes vectors
$\bu^j$, $j = 1,\ldots,\widetilde{M}$, and $\bv^k$, $k=1,\ldots,\widetilde{N}$, that are orthonormal.
Therefore, \emph{the columns of $\rU$ form an orthogonal basis in the universal reduced space $\widetilde{V}$.} 
The dimension of this space is defined by the first Tucker rank, $\mbox{\rm dim}(\widetilde{V}) = \widetilde{M}$. 
The information about $\widetilde{\bPhi}$ to be transmitted to the online stage includes matrices $\rS_i$ 
and the core tensor $\bC$. Explicitly,
\begin{equation}
\label{eqn:coreHOSVD}
\mbox{online}(\widetilde{\bPhi}) = \left\{
\bC \in \R^{\widetilde{M} \times \widetilde{n}_1 \times \dots \times \widetilde{n}_D \times \widetilde{N}},~
\rS_i \in \mathbb{R}^{n_i\times \widetilde{n}_i},~{\small i=1,\dots,D}\right\}.
\end{equation}

To find coordinates of the local basis for $\balpha \in \cA$, define the $\balpha$-specific 
\emph{core matrix} $\rC_e (\balpha) $ as
\begin{equation}
\label{eqn:C_hosvd}
\rC_e (\balpha) = \bC \times_2 \left(\rS_1 \be^1 (\balpha)\right) \times_3 \left(\rS_2 \be^2 (\balpha)\right) 
\dots\times_{D+1} \left(\rS_D \be^D (\balpha)\right)  \in \R^{\widetilde{M} \times \widetilde{N}}.
\end{equation}
\rev{Using the definition of $k$-mode product, \eqref{eqn:extractbt} and \eqref{eqn:TDv}, one computes
	\[
	\begin{split}
\widetilde{\Phi}_e (\balpha) &= \sum_{j = 1}^{\widetilde{M}}
\sum_{q_1 = 1}^{\widetilde{n}_1}\dots
\sum_{q_D = 1}^{\widetilde{n}_D}
\sum_{k = 1}^{\widetilde{N}}
(\bC)_{j, q_1, \dots, q_D, k}   \langle\bsigma^{q_1}_1,\be^1(\alpha)\rangle \cdot \dotsc \cdot \langle\bsigma^{q_D}_D,\be^D(\alpha)\rangle\,  \bu^j \circ\bv^k\\
&= \sum_{j = 1}^{\widetilde{M}}
\sum_{q_1 = 1}^{\widetilde{n}_1}\dots
\sum_{q_D = 1}^{\widetilde{n}_D}
\sum_{k = 1}^{\widetilde{N}}
(\bC)_{j, q_1, \dots, q_D, k}  \left(\rS_1 \be^1 (\balpha)\right)_{q_1} \cdot \dotsc \cdot \left(\rS_D \be^D (\balpha)\right)_{q_D}  \bu^j \circ\bv^k\\
&= \sum_{j = 1}^{\widetilde{M}}
\sum_{k = 1}^{\widetilde{N}}
\left(\bC \times_2 \left(\rS_1 \be^1 (\balpha)\right) \times_3 \left(\rS_2 \be^2 (\balpha)\right) 
\dots\times_{D+1} \left(\rS_D \be^D (\balpha)\right)\right)_{jk}  \bu^j \circ\bv^k = \rU \rC_e (\balpha) \rV^T
	\end{split}
	\]
Consider the thin SVD of the core matrix
\begin{equation}
\rC_e (\balpha) = \rU_c \Sigma_c \rV_c^T.
\label{eqn:coresvdho}
\end{equation}
Combining this with the representation above we get
\begin{equation}
\widetilde{\Phi}_e (\balpha) = 
\left({\rU}\rU_c\right)\Sigma_c \left({\rV}\rV_c\right)^T,
\label{eqn:Phie_hosvd}
\end{equation}
which is
the thin SVD of $\widetilde{\Phi}_e (\balpha)$ since both matrices $\rU$ and $\rV$ are orthogonal. }
We conclude that \emph{the coordinates $\left\{\bbeta_1(\balpha),\dots,\bbeta_n(\balpha)\right\}$ 
of the local reduced basis in the universal space $\widetilde{V}$ are the first $n$ columns of $\rU_c$ from
\eqref{eqn:coresvdho}.} The parameter-specific basis is then 
$\{ \bz_i(\balpha) \}_{i=1}^{n}$, with $\bz_i(\balpha) = {\rU} \bbeta_i(\balpha)$ 
(not actually computed at the online stage).

To compute the low-rank HOSVD approximation \eqref{eqn:TDv} we employ the 
standard algorithm~\cite{de2000multilinear} based on repeated computations of truncated SVD for unfolded matrices.
In particular, one may compute $\widetilde{\bPhi}$ with either prescribed Tucker ranks or prescribed 
accuracy $\widetilde{\eps}$. Moreover, for fixed Tucker ranks one can show that the recovered 
$\widetilde{\bPhi}$ satisfies a quasi-minimization property~\cite{de2000multilinear}
of the form
\begin{equation}
\label{eqn:HOSVDopt}
\| \bPhi - \widetilde{\bPhi} \| \le \sqrt{D+2} \| \bPhi - \bPhi^{\rm opt} \| \quad 
\text{and} \quad \| \bPhi - \widetilde{\bPhi} \| \le \left(\sum_{i=1}^{D+1} \rev{\Delta_i^2} \right)^{\frac12},
\end{equation}
where $\bPhi^{\rm opt}$ is the best approximation to $\bPhi$ among all \rev{Tucker} tensors of the given rank 
(such approximation always exists), and $\rev{\Delta_i}$ measures truncated SVD error on the $i$-th step of the HOSVD.
We summarize the above in the following algorithm.

\vskip0.05in
\begin{algorithm}[HOSVD-TROM]
\label{alg:HOSVD-TROM}
\begin{itemize}
\item \textbf{Offline stage}.\\
\emph{Input:} snapshot tensor $\bPhi \in \R^{M \times n_1 \times \ldots \times n_D \times N}$ and
target accuracy $\widetilde{\eps}$;\\
\emph{Output:}  Compressed tensor ranks, HOSVD decomposition matrices 
as in \eqref{eqn:hosvdmat}, and core tensor $\bC$; \\
\emph{Compute:} Use algorithm~\cite{de2000multilinear} with prescribed accuracy $\widetilde{\eps}$
to compute HOSVD decomposition matrices and the core tensor.
\item \textbf{Online stage}. \\
\emph{Input:} $\mbox{\rm online}(\widetilde{\bPhi})$ as defined in \eqref{eqn:coreHOSVD}, 
reduced space dimension $n \le \min \{ \widetilde{M}, \widetilde{N} \} $, 
and parameter vector $\balpha \in \cA$; \\
\emph{Output:} Coordinates of the reduced basis in $\widetilde{V}$:
$\{ \bbeta_i(\balpha) \}_{i=1}^{n} \subset \mathbb{R}^{\widetilde{M}}$;\\
\emph{Compute:}
\begin{enumerate}
\item Use the core tensor $\bC$ and matrices $\rS_i$, {\small$i = 1,\ldots,D$},
to assemble the core matrix $\rC_e (\balpha) \in \R^{\widetilde{M} \times \widetilde{N}}$ as in \eqref{eqn:C_hosvd};
\item Compute the SVD of the core matrix $\rC_e (\balpha) = \rU_c \Sigma_c \rV_c^T$ with
$\rU_c = [\widetilde{\bu}_1, \ldots, \widetilde{\bu}_{\widetilde{M}}]$;
\item Set $\bbeta_i(\balpha) =\widetilde{\bu}_i$, {\small$i = 1,\ldots,n$}.
\end{enumerate}
\end{itemize}
\end{algorithm}

\subsection{TT-TROM} 
\label{sec:TT-TROM}

A third low-rank tensor decomposition of interest is the \emph{Tensor Train} ({TT})
decomposition~\cite{TT1}. In TT-TROM we seek to approximate the snapshot tensor 
with $\widetilde{\bPhi}$ in the TT-format, namely
\begin{equation}
\label{eqn:TTv}
\bPhi \approx \widetilde{\bPhi} =
\sum_{j_1=1}^{\widetilde r_1}\dots
\sum_{j_{D+1}=1}^{\widetilde r_{D+1}}
\bu^{j_1} \circ \bsigma^{j_1, j_2}_1 \circ \dots \circ \bsigma^{j_D, j_{D+1}}_D \circ \bv^{j_{D+1}},
\end{equation}
with $\bu^{j_1} \in \R^M$, $\bsigma^{j_i, j_{i+1}}_i \in \R^{n_i}$, and $\bv^{j_{D+1}} \in \R^N$,
where the positive integers $\widetilde r_i$ are referred to as the \emph{compression ranks} (or TT-ranks) 
of the decomposition. For higher order tensors the TT format is in general more efficient compared to HOSVD. 
This may be beneficial for large $D$, the dimension of parameter space.
In~\cite{TT1,TT2} a stable algorithm for finding $\widetilde{\bPhi}$ based on truncated SVD for a sequence 
of unfolding matrices was introduced and the optimality property similar to \eqref{eqn:HOSVDopt} was proved.

Once an optimal TT approximation \eqref{eqn:TTv} is computed,
we organize the vectors and matrices from \eqref{eqn:TTv} into matrices 
\begin{equation}
\rU = [\bu^1, \dots, \bu^{\widetilde{r}_1}] \in \R^{M \times \widetilde{r}_1}, \quad 
\rV = [\bv^1, \dots, \bv^{\widetilde{r}_{D+1}}] \in \R^{N \times \widetilde{r}_{D+1}}, 
\label{eqn:ttuv}
\end{equation}
and third order tensors $\bS_i \in \R^{\widetilde{r}_i \times n_i \times \widetilde{r}_{i+1}}$, 
defined entry-wise as
\begin{equation}
(\bS_i)_{j k q} = (\bsigma^{j q}_i)_k, \quad
j = 1, \ldots, \widetilde{r}_i, \quad
k = 1, \ldots, n_i, \quad
q = 1, \ldots, \widetilde{r}_{i+1},
\label{eqn:ttsi}
\end{equation}
for all $i = 1, \ldots, D$. Note that matrix $\rU$ is orthogonal and so \emph{its columns provide an orthogonal 
basis in the universal space $\widetilde{V}$}. The dimension of $\widetilde{V}$ is defined by the first 
TT-rank, $\mbox{dim}(\widetilde{V})=\widetilde{r}_1$.

While $\rU$ is an orthogonal matrix, the columns of $\rV$ are orthogonal, but not necessarily orthonormal. 
Thus, we introduce a diagonal scaling matrix
\begin{equation}
\rW_c = \mbox{diag} \left( \| \bv^1 \|, \ldots, \| \bv^{\widetilde{r}_{D+1}} \| \right) 
\in \R^{\widetilde{r}_{D+1} \times \widetilde{r}_{D+1}}.
\label{eqn:ttwc}
\end{equation}
The essential information about $\widetilde{\bPhi}$ to be transmitted to the online phase includes $\bS_i$ 
tensors and the scaling factors:
\begin{equation}
\label{eqn:coreTT}
\mbox{online}(\widetilde{\bPhi}) = 
\left\{\bS_i \in \R^{\widetilde{r}_i \times n_i \times \widetilde{r}_{i+1}},~{\small i=1,\dots,D},~
\rW_c \in \R^{\widetilde{r}_{D+1} \times \widetilde{r}_{D+1}} \right\}.
\end{equation}

To find the coordinates of the local basis, we define the parameter-specific 
\emph{core matrix} $\rC_e (\balpha) \in \R^{\widetilde{r}_1 \times \widetilde{r}_{D+1}}$ as the product
\begin{equation}
\label{eqn:C_TT}
\rC_e (\balpha) = \prod_{i=1}^{D} \left( \bS_i \times_2 \be^i (\balpha) \right).
\end{equation}
\rev{Using the definition of $k$-mode product, \eqref{eqn:extractbt} and \eqref{eqn:TTv}, one computes
	\[
	\begin{split}
		\widetilde{\Phi}_e (\balpha) &= 
		\sum_{j_1=1}^{\widetilde r_1}\dots
		\sum_{j_{D+1}=1}^{\widetilde r_{D+1}}
		\langle\bsigma^{j_1, j_2}_1,\be_1(\balpha)\rangle \cdot \dotsc \cdot \langle\bsigma^{j_D, j_{D+1}}_D,\be_D(\balpha)\rangle \, \bu^{j_1}\circ \bv^{j_{D+1}}\\ 
		&= 
		\sum_{j_1=1}^{\widetilde r_1}\dots
		\sum_{j_{D+1}=1}^{\widetilde r_{D+1}}
		\left( \bS_1 \times_2 \be^1 (\balpha) \right)_{j_1, j_2}  \cdot \dotsc \cdot \left( \bS_D \times_{D+1} \be^D (\balpha) \right)_{j_D, j_{D+1}}  \, \bu^{j_1}\circ \bv^{j_{D+1}}\\
		&=\sum_{j_1=1}^{\widetilde r_1}\sum_{j_{D+1}=1}^{\widetilde r_{D+1}}
		\Big(  \prod_{i=1}^{D} \left( \bS_i \times_2 \be^i (\balpha) \right) \Big)_{j_1j_{D+1}}  \, \bu^{j_1}\circ \bv^{j_{D+1}} = \rU \rC_e (\balpha) \rV^T.
\end{split}
\]	
Consider the SVD of the rescaled core matrix: 
\begin{equation}
\rC_e (\balpha) \rW_c = \rU_c \Sigma_c \rV_c^T.
\label{eqn:coresvdtt}
\end{equation}
Using this and the above representation of  $\widetilde{\Phi}_e (\balpha)$ we compute
\begin{equation}
\widetilde{\Phi}_e (\balpha) = \rU \rC_e (\balpha)\rW_c\rW_c^{-1} \rV^T = 
\left({\rU}\rU_c\right)\Sigma_c \left({\rV}\rW_c^{-1}\rV_c\right)^T.
\label{eqn:Phie_TT}
\end{equation}
The right-hand side of \eqref{eqn:Phie_TT} is the thin SVD of $\widetilde{\Phi}_e (\balpha)$, 
since matrices ${\rU}$, $\rU_c$, ${\rV}\rW_c^{-1}$, and $\rV_c$ are all orthogonal.} We conclude that 
\emph{the coordinates $\left\{\bbeta_1(\balpha),\dots,\bbeta_n(\balpha)\right\}$ of the local reduced basis 
in the universal space $\widetilde{V}$ are the first $n$ columns of $\rU_c$.} The parameter-specific basis is then 
$\{ \bz_i(\balpha) \}_{i=1}^{n}$, with $\bz_i(\balpha)={\rU}\bbeta_i(\balpha)$ 
(not actually computed at the online stage).

We summarize the above in the following algorithm.

\vskip0.05in
\begin{algorithm}[TT-TROM]
\label{alg:TT-TROM}
\begin{itemize}
\item \textbf{Offline stage}.\\
\emph{Input:} snapshot tensor $\bPhi \in \R^{M \times n_1 \times \ldots \times n_D \times N}$ and
target accuracy $\widetilde{\eps}$;\\
\emph{Output:} Compression ranks,  TT decomposition matrices 
and third order tensors as in \eqref{eqn:ttuv}--\eqref{eqn:ttsi}; \\
\emph{Compute:} Use algorithm from~\cite{TT1} with prescribed accuracy $\widetilde{\eps}$
to compute TT decomposition \eqref{eqn:TTv}.
\item \textbf{Online stage}. \\
\emph{Input:} $\mbox{\rm online}(\widetilde{\bPhi})$ as defined in \eqref{eqn:coreTT}, 
reduced space dimension $n\le\min\{\widetilde{r}_1,\widetilde{r}_{D+1}\}$, and parameter vector $\balpha \in \cA$; \\
\emph{Output:} Coordinates of the reduced basis in $\widetilde{V}$:
$\{ \bbeta_i(\balpha) \}_{i=1}^{n} \subset \mathbb{R}^{\widetilde{r}_1}$;\\
\emph{Compute:}
\begin{enumerate}
\item Use tensors $\bS_i$ to assemble the core matrix 
$\rC_e (\balpha) \in \R^{\widetilde{r}_1 \times \widetilde{r}_{D+1}}$ as in \eqref{eqn:C_TT};
\item Compute the SVD of the scaled core matrix $\rC_e (\balpha) \rW_c = \rU_c \Sigma_c \rV_c^T$ with
$\rU_c = [\widetilde{\bu}_1, \ldots, \widetilde{\bu}_{\widetilde{r}_1}]$;
\item Set $\bbeta_i(\balpha) =\widetilde{\bu}_i$, {\small$i = 1,\ldots,n$}.
\end{enumerate}
\end{itemize}
\end{algorithm}

\subsection{General parameter sampling} 

Grid-based sampling of parameter space can be computationally expensive or not applicable if the set of 
admissible parameters $\cA$ is not a box (or an image of a box) in Euclidean space. However, interpolatory 
TROMs introduced above can be extended to accommodate a more general sampling set $\hcA$.
If $\cA$ does have the Cartesian structure (a box or an image of a box), then one way to reduce offline 
computational costs is to compute the snapshots for only a few parameter values from a Cartesian grid 
$\hcA \subset \cA$. To recover the missing entries of the full snapshot tensor $\bPhi$, one may use a 
low-rank tensor completion method, e.g., one of those studied in
\cite{liu2012tensor,gandy2011tensor,huang2014provable,yuan2016tensor,bengua2017efficient}.
The low-rank completion can be performed for any of the three compressed tensor formats considered above. 
We shall investigate this option elsewhere. In this paper, we consider another (more general) approach.

With a slight abuse of notation, let ${\hcA} = \{\bhalpha_1, \dots,\bhalpha_K \} \subset \cA$ 
be a set of  sampled parameter values. We assume that $\hcA$ is a frame in $\mathbb{R}^D$
and so $K \ge D$. Note that $K$ does not obey \eqref{eqn:Acard} for a general sampling.
Given an out-of-sample vector of parameters $\balpha \in \cA$, let
\begin{equation}
\be \, : \, \balpha \to \mathbb{R}^K
\end{equation}
be the representation of $\balpha$ in $\hcA$, i.e.,
\begin{equation}
\balpha = \sum_{j=1}^K a_j \bhalpha_j, \quad 
\be (\balpha) = (a_1,\dots,a_K)^T,
\label{eqn:representa}
\end{equation}
with an additional constraint enforcing uniqueness of the representation.

Similarly to \eqref{eqn:interpapprox} for the Cartesian grid case, we assume that for a smooth 
function $\rev{g}: \cA \to \mathbb{R}$ it holds
\begin{equation}
\rev{g}(\balpha) \approx \sum_{j=1}^{K} a_j \rev{g}(\bhalpha_j).
\end{equation}
In Section~\ref{sec:geninterp} we describe one particular choice of $\be(\balpha)$ that is used in 
all numerical experiments reported in Section~\ref{sec:num}.

To assemble the snapshot tensor, for each $\bhalpha_j \in \hcA$, $j = 1, \ldots, K$,
collect the snapshot vectors 
$\bu(t_k, \bhalpha_j) = \left( u_1(t_k, \bhalpha_j), \ldots, u_M(t_k, \bhalpha_j) \right)^T$, $k=1,\dots,N$,
and arrange them in a third order tensor $\bPhi \in \R^{M \times K \times N}$ with entries
\begin{equation}
\label{eqn:Phi1}
(\bPhi)_{ijk} = u_i (t_k, \bhalpha_j),
\quad i = 1,\ldots,M, \quad j = 1,\ldots,K, \quad k=1,\ldots,N.
\end{equation}
Then, for any $\balpha \in \cA$, the parameter-specific reduced basis is defined as the first $n$ left singular vectors of
\begin{equation}
\label{eqn:PsiAlpha1}
\widetilde{\Phi}_e (\balpha) = \widetilde{\bPhi} \times_2 \be (\balpha),
\end{equation}
where $\widetilde{\bPhi}$ is a low rank approximation of the snapshot tensor $\bPhi$ with entries \eqref{eqn:Phi1}. 
The same three compressed tensor formats considered above (CP, HOSVD and TT) can be used for $\widetilde{\bPhi}$, 
with TT format being inferior to HOSVD (for 3D tensors both HOSVD- and TT-decomposition are Tucker tensors).
An orthogonal basis in the universal space of $\widetilde{\bPhi}$ and coordinates of a local parameter-specific 
basis in it are computed similarly to the Cartesian grid sampling cases considered previously in 
Sections~\ref{sec:CP-TROM}--\ref{sec:TT-TROM}. The only difference is that all the calculations therein 
are performed setting $D=1$ and replacing $\be^1(\balpha)$ with $\be(\balpha)$. 
This includes the optimality result \eqref{eqn:HOSVDopt}, where the factor $\sqrt{D+2}$ becomes $\sqrt{3}$.

\subsection{Complexity and compression analysis} 
\label{sec:compcomp}

Projection-based parametric ROM framework consists in general of the following steps.
\begin{itemize}
\item[(i)] High-fidelity simulations of \eqref{eqn:GenericPDE} to generate the snapshot tensor $\bPhi$;
\item[(ii)] \textbf{Offline stage}: computing the compressed approximation $\widetilde{\bPhi}$ to $\bPhi$ 
in one of low-rank tensor formats; 
\item[(iii)] Passing the $\mbox{online}(\widetilde{\bPhi})$ part of the compressed tensor to the online stage;
\item[(iv)] \textbf{Online stage}: using $\mbox{online}(\widetilde{\bPhi})$ to compute the coordinates of 
the parameter-specific reduced basis for an input $\balpha$;
\item[(v)] Solving \eqref{eqn:GenericPDE} projected onto the reduced space.
\end{itemize}
Since steps (i) and (v) are common for all projection-based ROM approaches, we focus below on the 
computational and storage/transmission costs invoked in steps (ii)--(iv). The necessary details on step (v) 
are included in Section~\ref{rem1}.

First, we discuss briefly the computational costs at the more expensive offline stage. 
For CP-TROM, the standard algorithm for finding $\widetilde{\bPhi}$ in CP format \eqref{eqn:CPv}
is the ALS method~\cite{harshman1970foundations} which for a \emph{given} CP rank $R$ 
iteratively fits a rank $R$ tensor $\widetilde{\Phi}$ by solving on each iteration $D+2$ 
least squares problems for the factors $\bu^r$, $\bsigma^{r}_i$, $i=1,\dots,D$, and $\bv^r$,
$r = 1,\ldots,R$. While straightforward to implement, the method is sensitive to the choice of 
initial guess and may converge slowly. We refer the reader to \cite{ReviewTensor} for a guidance on 
the literature on improving the efficiency of ALS and possible alternatives. On the other hand,
computing $\widetilde{\bPhi}$ in either HOSVD or TT formats relies on finding truncated SVDs 
for matrix unfoldings of $\bPhi$ \cite{de2000multilinear, TT1}. Therefore, the computational 
complexity and cost of step (ii) for HOSVD- and TT-TROM is essentially the same as that of 
standard POD-ROM.

Second, to measure the amount of information transmitted to the online stage at step (iii), 
we introduce the \emph{compression factor} CF, defined as
\begin{equation}
\label{eqn:CF}
\text{CF} = \frac{\#(\bPhi)}{\#(\mbox{online}(\widetilde{\bPhi}))}\,,
\end{equation}
where we denote by $\#(\bPsi)$ the number of floating point numbers needed to store 
a tensor $\bPsi$. Specifically, $\#(\bPhi) = MKN$ is simply the total number of entries in $\bPhi$,
while $\#(\mbox{online}(\widetilde{\bPhi}))$ is the number of entries needed to store all the factors 
passed to the online stage, as defined in \eqref{eqn:coreCP}, \eqref{eqn:coreHOSVD} and \eqref{eqn:coreTT}
for CP-, HOSVD- and TT- TROMs,  which we summarize in Table~\ref{tab:numPhi}.  
 
\begin{table}[h]
\centering
\caption{Number of entries needed to store $\mbox{online}(\widetilde{\bPhi})$ for CP, HOSVD and TT formats.}
\begin{tabular}{l|cc}
 & \multicolumn{2}{c}{$\#(\mbox{online}(\widetilde{\bPhi}))$} \\[0.5ex]
Format & {\small Cartesian grid-based} & \small {General} \\[0.5ex]
\hline\\[-2.2ex]
CP    & $R\big(\sum\limits_{i=1}^{D}n_i+R+1\big)$ & $R\big(K+R+1\big)$\\
HOSVD & $\widetilde{N}\widetilde{M}\prod\limits_{i=1}^{D}\widetilde{n}_i  + \sum\limits_{i=1}^{D}\widetilde{n}_i n_i$ &
$\widetilde{N}\widetilde{M} \widetilde{n}_1  + \widetilde{n}_1 K$ \\ 
TT  & $ \widetilde{r}_{D+1} + \sum\limits_{i=1}^{D} \widetilde{r}_{i} n_i \widetilde{r}_{i+1}$ & 
$ \widetilde{r}_{2} + \widetilde{r}_{1} K \widetilde{r}_{2}$
\\[0.5ex] 
\end{tabular}
\label{tab:numPhi}
\end{table}

Table~\ref{tab:numPhi} shows that the compression factor
is largely determined by the compression ranks. In turn, the ranks depend on $\widetilde{\eps}$ and variability of observed states.

Third, the computational complexity of finding $\balpha$-specific reduced basis in step (iv) is determined 
by the interpolation procedure and the computation of first $n$ left singular \rev{vectors of}  the core matrix.
 Since vectors $\be^i (\balpha)$ contain very few non-zero entries,
e.g., $p = 2$ or $3$ of them for the Cartesian sampling,  the number of operations for 
computing core matrices $\rC(\balpha)$ for CP-, HOSVD- and TT-TROM is
\begin{equation}
O \left( R^2 \right), ~~
O \Big( \widetilde{M} \widetilde{N} \prod\limits_{i=1}^{D} \widetilde{n}_i \Big), ~~\text{and}~~ 
O \Big( \sum\limits_{i=2}^{D} \widetilde{r}_{i-1} \widetilde{r}_{i} \widetilde{r}_{i+1} \Big), 
\end{equation}
respectively. CP-, HOSVD- and TT-TROM algorithms proceed to
compute the SVD of small core matrices of sizes $R\times R$, $\widetilde{M} \times \widetilde{N}$ and 
$\widetilde{r}_{1} \times \widetilde{r}_{D+1}$, respectively. 
If a reduced basis in the physical space is desired, then one finds its vectors as linear combinations 
of columns of $\rU$, which requires $O(M R n)$, $O(M \widetilde{M} n)$ or $O(M \widetilde{r}_1 n)$ 
operations for CP-, HOSVD- or TT-TROM, respectively. Section~\ref{rem1} below discusses how these 
costs can be avoided at the online phase.
 We note that for a fixed compression accuracy $\widetilde{\epsilon}$, it is often observed
in practice that the corresponding ranks of HOSVD and TT formats satisfy $\widetilde{M} \simeq \widetilde{r}_1$, $\widetilde{N} \simeq \widetilde{r}_{D+1}$.

In summary, the computational costs of the offline stage for TROMs are comparable to those of POD-ROM 
for a multi-parameter  problem. At the online stage complexity of all preparatory steps depends only on compressed 
tensor ranks rather than the size of the snapshot tensor $\bPhi$. 
The amount of  information transmitted from offline to online 
stages is  determined  by the compressed tensor ranks, as should be clear from Table~\ref{tab:numPhi}.

\subsection{TROM evaluation}
\label{rem1}

Besides finding a suitable reduced basis,  a fast evaluation of the reduced model for any incoming 
$\balpha\in\cA$ is required for a reduced modeling scheme to be effective. 
Efficient implementation of a projected parametric model is a well-known challenge that have been addressed 
in the literature with various approaches; see, e.g.,
~\cite{BrennerScott,grepl2007efficient,chaturantabut2010nonlinear,drohmann2012reduced,
carlberg2013gnat,benner2015survey,hesthaven2016certified,kramer2019nonlinear}.
The tensorial approach presented here does not directly contribute to resolving this issue, but it does not make 
it harder either and so techniques known from the literature can be adapted in the TROM framework. 

For example, assume that $F(t,\bu,\balpha)$ from \eqref{eqn:GenericPDE} has an affine dependence on 
parameters and linear dependence on $\bu$:
\[
F(t,\bu,\balpha)=\sum_{i=1}^{P}f_i(\balpha)\rA_i\bu,
\]
with some $f_i:\cA\to\mathbb{R}$ and parameter-independent  $\rA_i\in \mathbb{R}^{M\times M}$. 
We assume that $P$ is not too large, at least independent of other dimensions. Then the \emph{offline} 
stage of model reduction consists of projecting matrices onto the universal space by computing 
$\widehat{\rA}_i=\rU^T\rA_i\rU$, where $\rU$ is an orthogonal basis matrix for $\widetilde{V}$ 
provided by the tensor decompositions. 
The new matrices $\widehat{\rA}_i$ have the reduced size $T_r \times T_r$, 
with $T_r \in \{R,\widetilde{M},\widetilde{r}_1\}$ for CP-, HOSVD- and TT-TROMs, respectively. 

\rev{For each of TROMs, denote by $\rU_c(n)$ the matrix of the first $n$ columns of $\rU_c$, left singular vectors of $\alpha$-specific core matrices.} During the \emph{online} stage one solves the system projected further on the parameter-specific local basis:    
\[
\bv_t=\sum_{i=1}^{P}f_i(\balpha)\rev{\rU^T_c(n)}\widehat{\rA}_i\rev{\rU_c(n)}\bv,
\]
where $\bv(t)$ is the trajectory in a space spanned by the columns of $\rev{\rU_c(n)}\in \mathbb{R}^{T_r \times n}$, 
i.e., the corresponding physical states are given by $\bu(t)=\rU\rev{\rU_c(n)}\bv(t)$. We see that online computations 
depend only on reduced dimensions (tensor ranks) and the small dimension $n$ of parameter-specific basis. 
This observation can be extended to the case when $F$ has a low order polynomial non-linearity with respect to $\bu$. 
For example, quadratic nonlinear terms, as in Burgers or Navier-Stokes equations, can be evaluated in 
$O(T_r^2)$ operations on each time step given a vector $\bv$ in the local reduced basis.   

To evaluate more general nonlinear terms, one can use a hyper-reduction technique such as the discrete empirical 
interpolation method (DEIM)~\cite{chaturantabut2010nonlinear}. In this approach, the nonlinear term is approximated 
in a basis of its snapshots. As an example, consider $F(t,\bu,\balpha) = \rA \bu + f(t, \bu(t),\balpha)$, 
with $\rA \bu$ representing linear and $f(t,\bu(t),\balpha)$ representing the non-linear part of $F$. 
Define the snapshots $\blf_i = f(t_i, \bu(t_i), \balpha_i)$, $i=1,\dots,N_{\rm DEIM}$ for some ``greedy'' 
choice of parameters and time instances and high-fidelity solution $\bu$. Denote by $\rQ$ an orthogonal basis 
matrix for $\mbox{span}\{\blf_1,\dots,\blf_{N_{\rm DEIM}}\}$, then DEIM approximates
\[
f(t,\bu,\balpha)\approx \rQ (\rP\rQ)^{-1}\rP f(t,\bu,\balpha),
\]  
where $\rP^T$ is an $N_{\rm DEIM}\times M$ ``selection'' matrix such that for any $f\in\mathbb{R}^M$,  
$\rP f$ contains $N_{\rm DEIM}$ selected entries of $f$. A particular $\rP$ corresponds to the choice of 
spatial interpolation nodes, cf.~\cite{chaturantabut2010nonlinear}. In TROM one may pre-compute 
$\widehat{\rQ} = \rU^T \rQ$ and $\widehat{\rA} = \rU^T \rA \rU$ during the offline stage, 
then solve at the online stage
\[
\bv_t = \rev{\rU_c^T(n)}\widehat{\rA}\rU_c(n)\bv + \rev{\rU^T_c(n)}\widehat{\rQ} (\rP\rQ)^{-1}\rP f(t,\rU\rev{\rU_c(n)} \bv,\balpha),
\] 
with costs depending on compressed tensor ranks, $n$, and the dimension of  DEIM space, but not on the dimensions of high-fidelity simulations. It is an interesting question, whether the tensor technique can be applied to make the  DEIM space  parameter-specific for more efficient reduce online computations. We plan to address this question elsewhere.   

\section{Prediction analysis} 
\label{sec:analysis}

In this section we assess the prediction power of the reduced basis $\cZ_n(\balpha) = \{ \bz_1, \ldots, \bz_n \}$ 
consisting of the first $n$ left singular vectors of $\widetilde\Phi_e (\balpha)$ from \eqref{eqn:extractbt}, 
for a parameter $\balpha = (\alpha_1,\ldots,\alpha_D)^T \in \cA$, not necessarily from a sampling set; 
i.e., $[ \bz_1, \ldots, \bz_n ]=\rU\rev{\rU_c(n)}$.

For the discussion below we also need the following notation. Given an $\balpha \in \cA$, 
we denote by $\bpsi_i = \bu(t_i, \balpha) \in \mathbb{R}^M$, $i = 1,\dots,N$, the snapshots of a high-fidelity solution 
to  \eqref{eqn:GenericPDE} and let $\Psi(\balpha) = [\bpsi_1,\dots,\bpsi_N] \in \R^{M \times N}$ 
be the corresponding snapshot matrix. Note that in practice the snapshots for out-of-sample parameters
are not available, so the matrix $\Psi(\balpha)$ should be treated as unknown.

We estimate the prediction power of $\cZ_n(\balpha)$ in terms of the quantity
\begin{equation}
	\label{eqn:quant}
	E_n(\balpha) = \frac{1}{NM} \sum_{i=1}^{N}
	\left\| \bpsi_i - \sum_{j=1}^{n} \langle \bpsi_i, \bz_j \rangle \bz_j \right\|^2_{\ell^2},
\end{equation}
which measures how accurate the solution $\bu(t, \balpha)$ at time instances $t_i$ can be represented in 
the reduced basis for the arbitrary but fixed $\balpha \in \cA $. The scaling $1/(NM)$ accounts for 
the variation of dimensions $N$ and $M$, which may correspond to the number of temporal and spatial 
degrees of freedom, respectively, if \eqref{eqn:GenericPDE} comes from a discretization of a parabolic PDE 
defined in a spatial domain $\Omega$. In this case and for uniform grids, the quantity in \eqref{eqn:quant} 
is consistent with the $L^2(0, T, L^2(\Omega))$ norm.

Below we prove an estimate for $E_n(\balpha)$ in terms of $\widetilde{\eps}$ from \eqref{eqn:TensApprox}, 
the singular values of $\widetilde{\Phi}_e (\balpha)$ and interpolation properties of $\be^i (\balpha)$, $i=1,\ldots,D$.
To make use of the latter, we introduce the following quantities related to the interpolation procedure. For Cartesian 
grid-based sampling we define the maximum grid step 
\begin{equation}
	\delta_i = \max\limits_{1 \leq j  \leq n_i-1} \left| \widehat{\alpha}_j^i - \widehat{\alpha}_{j+1}^i \right|, 
	\quad i = 1,\ldots,D.
\end{equation}
Relation \eqref{eqn:lagrange} implies that the interpolation procedure \eqref{eqn:bea}--\eqref{eqn:lagrange}
is of order $p$, i.e., for any sufficiently smooth $f:  [\alpha_i^{\min}, \alpha_i^{\max}] \to \R$ it holds
\begin{equation}
	\label{eqn:int_aprox}
	\sup_{a \in [\alpha_i^{\min}, \alpha_i^{\max}]} 
	\Big| f(a) - \sum_{j=1}^{n_i} e^i_j \left( a \be_i \right) f(\halpha_i^j) \Big|\le
	C_a \| f^{(p)} \|_{C( [\alpha_i^{\min}, \alpha_i^{\max}])} \delta_i^p, 
\end{equation}
for $i = 1, \dots, D$, where $\be_i \in \R^{n_i}$ is the $i^{\mbox{\scriptsize th}}$ column of an
$n_i \times n_i$ identity matrix. The constant $C_a$ does not depend on $f$. We let $\delta^p = \sum_{i=1}^{D} \delta_i^p$ 
and also assume that the interpolation procedure is stable in the sense that
\begin{equation}
	\label{eqn:int_stab}
	\left(\sum_{j=1}^{n_i} \left| e^i_j(a \be_i) \right|^2 \right)^{\frac12} \le C_e,
\end{equation}
with some $C_e$ independent of $a \in [\alpha_i^{\min}, \alpha_i^{\max}]$ and $i = 1, \ldots, D$.
For the example of linear interpolation with $p=2$ and $\alpha_i^{\min}$, $\alpha_i^{\max}$ included 
among the grid nodes, bounds \eqref{eqn:int_aprox} and \eqref{eqn:int_stab} hold with $C_a = \frac18$
and $C_e = 1$.

To estimate $E_n(\balpha)$, consider the SVD of $\widetilde{\Phi}_e (\balpha) \in \R^{M \times N}$ given by
\begin{equation}
	\widetilde{\Phi}_e(\balpha) = \widetilde{\rU} \widetilde{\Sigma} \widetilde{\rV}^T,
	~~\text{with}~~ \widetilde{\Sigma} = \text{diag} (\tsigma_1,\dots,\tsigma_N).
\end{equation}
Then $\rZ = [\bz_1, \ldots, \bz_n] \in \R^{M \times n}$ is build \rev{as} the first $n$ columns of 
$\widetilde{\rU}$, the reduced basis vectors, i.e., $\rZ=\rU\rU_c$. Then,
\begin{align}
	E_n(\balpha) &
	= \frac{1}{NM} \left\| (\rI  - \rZ \rZ^T) \Psi(\balpha) \right\|^2_F \nonumber \\ 
	& \le \frac{1}{NM} \left( \left\| (\rI - \rZ \rZ^T) (\Psi(\balpha) - \widetilde{\Phi}_e(\balpha)) \right\|_F
	+ \left\| (\rI - \rZ \rZ^T) \widetilde{\Phi}_e(\balpha) \right\|_F \right)^2 \nonumber \\
	& \le\frac{1}{NM} \left( \left\| {\Psi}(\balpha) - \widetilde{\Phi}_e(\balpha) \right \|_F 
	+ \left\| (\rI - \rZ \rZ^T) \widetilde{\Phi}_e(\balpha) \right\|_F \right)^2, \label{eqn:aux413}
\end{align}
where we used triangle inequality and $\| \rI - \rZ \rZ^T\| \le 1$ for the spectral norm of the projector.
For the last term in \eqref{eqn:aux413}, we observe
\begin{equation}
	\left\| (\rI - \rZ \rZ^T) \widetilde{\Phi}_e(\balpha) \right\|_F =
	\left\| \widetilde{\rU} \; \text{diag}(0,\dots,0,\tsigma_{n+1},\dots,\tsigma_N) \; \widetilde{\rV}^T \right\|_F =
	\left( \sum_{j = n+1}^{N} \tsigma_j^2 \right)^{\frac12}.
\end{equation}
To handle the first term of \eqref{eqn:aux413}, consider the extraction 
\begin{equation}
	\Phi_e (\balpha) = \bPhi \times_2 \be^1 (\balpha) \times_3 \be^2(\balpha) \dots \times_{D+1} \be^D(\balpha)
\end{equation}
and proceed using the triangle inequality
\begin{equation}
	\label{eqn:aux438}
	\big\| \Psi (\balpha) - \widetilde{\Phi}_e(\balpha) \big\|_F \le 
	\big\| \Psi (\balpha) - \Phi_e (\balpha) \big\|_F + 
	\big\| \Phi_e(\balpha) - \widetilde{\Phi}_e (\balpha) \big\|_F.
\end{equation}
We use the stability of interpolation \eqref{eqn:int_stab} and \eqref{eqn:TensApprox} to bound the second term
of \eqref{eqn:aux438}. Specifically,
\begin{equation}
	\label{eqn:aux442}
	\begin{split}
		\left\| \Phi_e (\balpha) - \widetilde{\Phi}_e(\balpha) \right\|_F & =
		\left\| (\bPhi - \widetilde{\bPhi}) 
		\times_2 \be^1(\balpha) \times_3 \be^2(\balpha) \dots \times_{D+1} \be^D (\balpha) \right\|_F \\
		& \le \left\| \bPhi - \widetilde{\bPhi} \right\|_F 
		\| \be^1 (\balpha) \|_{\ell^2} \| \be^2(\balpha) \|_{\ell^2} \dots \| \be^D (\balpha) \|_{\ell^2} \\
		& \le (C_e)^D \left\| \bPhi - \widetilde{\bPhi} \right\|_F 
		\le (C_e)^D \widetilde{\eps}\left\| \bPhi \right\|_F.
	\end{split}
\end{equation}
It remains to handle the first term in \eqref{eqn:aux438}. At this point we need more precise assumptions 
on the smoothness of $\bu(t,\balpha)$, the solution of \eqref{eqn:GenericPDE}. In particular,
\begin{equation}
	\label{eqn:AssU}
	\bu \in C([0,T] \times \overline{\cA})^M, \quad 
	\frac{\partial^{\mathbf{j}} \bu}{\partial \alpha^{j_1}_1\dots\alpha^{j_D}_D} 
	\in C([0,T] \times \overline{\cA})^M, \quad |\mathbf{j}| \le p.
\end{equation}
We note that \eqref{eqn:AssU} is guaranteed to hold if the unique solution to \eqref{eqn:GenericPDE} 
exists on $(0,T_1)$ for all $\balpha \in\cA_1$, with  $T < T_1$ and $\overline{\cA} \subset \cA_1$
and $F$ is continuous with continuous partial derivatives in components of $\bu$ and $\balpha$ 
of order up to $p$~\cite{Hartman}.

Using interpolation property \eqref{eqn:int_aprox}, we compute
\begin{equation}
	\begin{split}
		\left( \bPhi \times_2 \be^1(\balpha) \right)_{:,i_2,\dots,i_D,k} & 
		= \sum_{j=1}^{n_1} e^1_j(\balpha) 
		\bu(t_k, \halpha_1^{j}, \halpha_2^{i_2}, \dots, \halpha_D^{i_D}) \\
		& = \bu(t_k, \alpha_1, \halpha_2^{i_2} \dots, \halpha_D^{i_D}) + \Delta^1_{:,i_2,\dots,i_D,k},
	\end{split}
\end{equation}
with the remainder term obeying a component-wise bound
\begin{equation}
	|\Delta^1_{:,i_2,\dots,i_D,k}| \le C_a \sup_{a \in [\alpha_1^{\min}, \alpha_1^{\max}]}
	\left| \frac{\partial^p \bu}{\partial \alpha^p_1}(t_k, a, \halpha_2^{i_2}, \dots, \halpha_D^{i_D}) \right| \delta_1^p,
\end{equation}
where the absolute value of vectors is understood entry-wise.
Analogously, we compute
\[
\begin{split}
	\big(\bPhi \times_2 \be^1(\balpha) & \times_3 \be^2(\balpha) \big)_{:,i_3,\dots,i_D,k} = 
	\left( (\bPhi \times_2 \be^1(\balpha)) \times_2 \be^2(\balpha) \right)_{:,i_3,\dots,i_D,k}\\
	& = \sum_{j=1}^{n_2} e^2_j(\balpha) 
	\left( \bu (t_k, \alpha_1, \halpha_2^{j}, \halpha_3^{i_3}, \dots, \halpha_D^{i_D}) 
	+ \Delta^1_{:, j, i_3,\dots,i_D,k} \right) \\
	& = \bu (t_k, \alpha_1, \alpha_2, \halpha_3^{i_3}, \dots, \halpha_D^{i_D}) 
	+ \Delta^2_{:,i_3,\dots,i_D,k} + \sum_{j=1}^{n_2}e^2_j(\balpha) \Delta^1_{:,j,i_3,\dots,i_D,k},
\end{split}
\]
with a component-wise bound for the remainder
\begin{equation}
	\begin{split}
		\Big| \Delta^2_{:, i_3, \dots, i_D, k} + 
		\sum_{j=1}^{n_2} e^2_j(\balpha) & \Delta^1_{:, j, i_3, \dots, i_D, k} \Big| \\
		\le C_a & \sup_{a \in [\alpha_2^{\min}, \alpha_2^{\max}]} 
		\left| \frac{\partial^p \bu}{\partial \alpha^p_2}
		(t_k, \alpha_1, a, \halpha_3^{i_3}, \dots, \halpha_D^{i_D}) \right| \delta_2^p \\
		+ \; C_e \; C_a & \sup_{a \in [\alpha_1^{\min}, \alpha_1^{\max}]}
		\left| \frac{\partial^p \bu}{\partial \alpha^p_1}
		(t_k, a, \halpha_2^{i_2}, \dots, \halpha_D^{i_D}) \right| \delta_1^p.
	\end{split}
\end{equation}
Applying the same argument repeatedly, we obtain
\begin{equation}
	\begin{split}
		\left( \Phi_e (\balpha) \right)_{:, k} & = 
		\left( \bPhi \times_2 \be^1(\balpha) \times_3 \be^2(\balpha) \dots \times_{D+1} \be^D(\balpha) \right)_{:, k} \\
		& = \bu (t_k, \alpha_1, \alpha_2, \dots, \alpha_D) 
		+ \Delta_{:, k} = \left( \Psi(\balpha) \right)_{:, k} + \Delta_{:,k},
	\end{split}
\end{equation}
with a component-wise bound for the remainder
\begin{align*}
	\left| \Delta_{:,k} \right| \le C_a \Big( & \sup_{a \in [\alpha_D^{\min}, \alpha_D^{\max}]}
	\left|\frac{\partial^p \bu}{\partial \alpha^p_D}
	(t_k, \alpha_1, \dots, \alpha_{D-1}, a) \right| \delta_D^p + \dots \\
	+\, (C_e)^{D-2} & \sup_{a \in [\alpha_2^{\min}, \alpha_2^{\max}]}
	\left| \frac{\partial^p \bu}{\partial \alpha^p_2}
	(t_k, \alpha_1, a, \halpha_3^{i_3}, \dots, \halpha_D^{i_D}) \right| \delta_2^p \\
	+\, (C_e)^{(D-1)} & \sup_{a \in [\alpha_1^{\min}, \alpha_1^{\max}]}
	\left| \frac{\partial^p \bu}{\partial \alpha^p_1}
	(t_k, a, \halpha_2^{i_2}, \dots, \halpha_D^{i_D}) \right| \delta_1^p \Big).
\end{align*}
Using the definition of the Frobenius norm,  we arrive at
\begin{equation}
	\label{eqn:aux450}
	\left\| \Psi(\balpha) - \Phi_e (\balpha) \right\|_F \le 
	\sqrt{N M} \; C_a C_\bu \max \left\{ (C_e)^{(D-1)}, 1 \right\} \; \delta^p,
\end{equation}
where $C_\bu$ depends only on the smoothness of $\bu$ with respect to the variations of 
parameters $\balpha$. More precisely, we can take $C_\bu = \| \bu \|_{C(0,T; C^p(\cA))}$, 
which is bounded due to assumption \eqref{eqn:AssU}.

Summarizing \eqref{eqn:aux413}--\eqref{eqn:aux450}, we proved the following result.

\begin{theorem}
	\label{Th1}
	Assume the solution $\bu$ to \eqref{eqn:GenericPDE} satisfies \eqref{eqn:AssU}, $\hcA$ is a Cartesian 
	grid in parameter domain $\cA$. Then for any $\balpha \in \cA$ the interpolatory TROM reduced basis 
	$\cZ(\balpha) = \{ \bz_1, \ldots, \bz_n \}$ delivers the following  representation estimate
	\begin{multline}
		\label{eqn:est1}
		\frac1{3NM} \sum_{i=1}^{N} 
		\left\| \bu(t_i, \balpha) - \sum_{j=1}^{n} 
		\left\langle \bu(t_i, \balpha), \bz_j \right\rangle \bz_j \right\|^2_{\ell^2} \\
		\le \frac1{NM} \left( (C_e)^{2D} \widetilde{\eps}^2\left\| \bPhi \right\|_F + 
		\sum_{i = n+1}^{N} \tsigma_i^2 \right) + C_a C_\bu \max \left\{ (C_e)^{2 (D-1)}, 1 \right\} \delta^{2p} ,
	\end{multline}
	with $C_\bu = \| \bu \|_{C(0,T; C^p(\cA))}$ independent of $\balpha$, sampling grid and $n$.
\end{theorem}

We summarize here the definitions of quantities that appear in Theorem~\ref{Th1}: $N$ is a number of 
time steps for snapshot collection, while $M$ is the spatial dimension of snapshots, i.e., 
$\bu(t_i, \balpha) \in \R^M$, $i=1,\ldots,N$; $\widetilde{\eps}$ is the relative accuracy of the snapshot tensor 
compression from~\eqref{eqn:TensApprox}; $\tsigma_i$ are the singular values of $\widetilde{\Phi}_e (\balpha)$ 
from \eqref{eqn:extractbt} (note that in TROMs we have an access to $\tsigma_i$ as the singular values 
of core matrices \eqref{eqn:cpcore}, \eqref{eqn:coresvdho} and \eqref{eqn:coresvdtt} for CP-, HOSVD-, 
and TT-TROM, respectively); 
$\delta = \big( \sum_{i=1}^{D} \delta_i^p \big)^\frac1p$ is the grid step parameter of the Cartesian grid in $\cA$; 
$p$ is both the number of nearest grid points and the order of interpolation of the interpolation procedure 
\eqref{eqn:bea}--\eqref{eqn:lagrange}; $C_e$ is the interpolation stability constant from \eqref{eqn:int_stab}; 
and $D$ is the dimension of parameter space.

For the general parameter sampling, prediction power analysis follows the same lines as above, simply setting $D=1$.
However, the order of interpolation $p$ is slightly more difficult to formalize, so instead of 
\eqref{eqn:int_aprox}--\eqref{eqn:int_stab} we rather assume
\begin{align}
	\label{eqn:int_stab2}	
	\sup_{\balpha\in \cA} \left| f(\balpha) - \sum_{j=1}^{K} (\be(\balpha))_j f(\bhalpha_j) \right| \le
	\|f\|_{C^p(\cA)}\delta,\quad
	\Big(\sum_{j=1}^{K}| (\be(\balpha))_j |^2\Big)^\frac12 \le C_e
\end{align}
with some $\delta$ depending on $\hcA$. The prediction estimate then becomes
\begin{multline}
	\label{eqn:est2}
	\frac1{3NM} \sum_{i=1}^{N} 
	\left\| \bu(t_i, \balpha) - \sum_{j=1}^{n} \left\langle \bu(t_i,\balpha), \bz_j \right\rangle \bz_j \right\|^2_{\ell^2} \\
	\le \frac1{NM} \left( (C_e)^{2} \widetilde{\eps}^2\left\| \bPhi\right\|_F + \sum_{i=n+1}^{N} \tsigma_i^2 \right) + 
	C_\bu \max \{ (C_e)^2, 1\} \delta^2,
\end{multline}
with $C_\bu = \| \bu \|_{C(0,T; C^p(\cA))}$. 

We finally, note that the feasibility of a sufficiently accurate lower rank representation of $\bPhi$ depends on the smoothness of $\bu$ as a function of $\bx,$ $t$ and $\balpha$. This question can be addressed by considering tensor decompositions of multivariate functions, e.g.~\cite{hackbusch2012tensor,nouy2015low}.  For these functional CP, HOSVD and hierarchical \rev{Tucker} (including TT) formats, the dependence of compression ranks  on $\widetilde{\eps}$ from \eqref{eqn:TensApprox} and the regularity (smoothness) of $\bu$ was studied in \cite{griebel2021analysis,trefethen2017multivariate,schneider2014approximation,hackbusch2007tensor,temlyakov1988estimates}. 
This compression property for multivariate functions was exploited to effectively represent solutions of parametric elliptic PDEs using tensor formats in~\cite{khoromskij2011tensor,cohen2011analytic,ballani2017multilevel,eigel2017adaptive,dolgov2018direct} among other publications. 

\section{Numerical experiments}
\label{sec:num}

We perform several numerical experiments to assess the performance of the three TROM approches
and compare them to the conventional POD-ROM. The testing in Section~\ref{sec:linpde} is performed
for a dynamical system originating from a discretization of linear parameter-dependent heat equation. 
In Section~\ref{sec:advdiff} a similar set of tests is carried out for a time-dependent parameterized advection-diffusion
system.

\subsection{General parameter sampling interpolation scheme}
\label{sec:geninterp}

For the numerical examples in the general parameter sampling setting we employ the following interpolation scheme.
Fix an integer $q \geq D+1$ and let $\balpha = (\alpha_1, \ldots, \alpha_D) \in \cA$ be an out-of-sample parameter vector. 
The interpolation scheme is based on the weighted minimum norm fit over $q$ nearest neighbors of 
$\balpha$ in the sampling set. Thus, we denote by $\bhalpha_{i_1},\ldots,\bhalpha_{i_q}$ the $q$ closest 
parameter samples in $\hcA$ to $\balpha$ and set $d_k = \| \bhalpha_{i_k} - \balpha \| > 0$, $k = 1, \ldots, q$. 
Next, define the weighting matrix
\begin{equation*}
	{\mathrm D} = \mbox{diag}(d_1^{-1}, \ldots, d_q^{-1}) \in \R^{q \times q}.
\end{equation*}
Also, assemble the matrix
\begin{equation*}
	{\mathrm X} = \begin{bmatrix}
		(\bhalpha_{i_1})_1 & (\bhalpha_{i_2})_1 & \cdots & (\bhalpha_{i_q})_1 \\
		(\bhalpha_{i_1})_2 & (\bhalpha_{i_2})_2 & \cdots & (\bhalpha_{i_q})_2 \\
		\vdots & \vdots & \vdots & \vdots \\
		(\bhalpha_{i_1})_D & (\bhalpha_{i_2})_D & \cdots & (\bhalpha_{i_q})_D \\
		1 & 1 & \cdots & 1
	\end{bmatrix}  \in \R^{D+1 \times q}.
\end{equation*}
Solve the weighted minimum norm fitting problem to obtain
\begin{equation}
	\widehat{\ba} = {\mathrm D} ({\mathrm X} {\mathrm D})^{\dagger} 
	\begin{bmatrix} \balpha \\ 1 \end{bmatrix} \in \R^q.
	\label{eqn:wminnorm}
\end{equation}
Note that the last row of ${\mathrm X}$ and the last entry of $(\balpha, 1)^T$ enforces the condition
that the entries of $\widehat{\ba}$ sum to one. Meanwhile, the presence of the weighting matrix 
puts more emphasis on the neighbors of $\balpha$ that are closest to it. 

Once $\widehat{\ba}$ is obtained from \eqref{eqn:wminnorm}, we define $a_j$, $j=1,\ldots,K$,
the entries of $\be(\balpha)$ as
\begin{equation*}
	a_j = \begin{cases} 
		\widehat{a}_{k}, & \text{if } j = i_k \in \{i_1, \ldots, i_q \} \\
		0, & \text{otherwise} \end{cases}
\end{equation*}
Clearly, such construction enforces representation \eqref{eqn:representa}.

\subsection{Parameterized heat equation}
\label{sec:linpde}

We first assess performance of the three TROM approaches on a dynamical system resulting from the 
discretization of a  heat equation
\begin{equation}
w_t = \Delta w,
\label{eqn:heat}
\end{equation}
in a rectangular domain with three holes 
$\Omega = \Omega_r \setminus (\Omega_1 \cup \Omega_2 \cup \Omega_3) \subset \mathbb{R}^2$,
where $\Omega_r = [0, 10] \times [0, 4]$, and the holes are $\Omega_1 = [1, 3] \times [1, 3]$,
$\Omega_2 = [4, 6] \times [1, 3]$, $\Omega_3 = [7, 9] \times [1, 3]$. The PDE and geometry
of $\Omega$ follow that of \cite{grepl2005posteriori}, while the boundary conditions are modified 
from those used in \cite{grepl2005posteriori}, as described below.

We parametrize the system with $D = 4$ parameters that enter the boundary conditions.
Convection boundary conditions are enforced on the left side of the rectangle 
$\Gamma_o = 0 \times [0, 4]$ and on the boundaries of each hole $\partial \Omega_j$, $j=1,2,3$.
Explicitly,
\begin{equation}
\left. (\bn \cdot \nabla w + \alpha_1(w-1)\,) \right|_{\Gamma_o} = 0,
\label{eqn:bco}
\end{equation}
and
\begin{equation}
\left. \left( \bn \cdot \nabla w + \frac{1}{2} w \right) \right|_{\partial \Omega_j} = \frac{1}{2} \alpha_{j+1},
\quad j=1,2,3,
\label{eqn:bcj}
\end{equation}
i.e., the first parameter in $\balpha \in \mathbb{R}^4$ is Biot number at $\Gamma_o$ with a fixed outside 
temperature $t_o = 1$, while the other three parameters are the temperatures at $\partial \Omega_j$,
$j=1,2,3$, respectively, with Biot numbers equal to $\frac{1}{2}$ on all three hole boundaries. 
The rest of the boundary of $\Omega$ is assumed to be insulated 
\begin{equation}
\left. (\bn \cdot \nabla w) \right|_{\partial \Omega_r \setminus \Gamma_o} = 0.
\label{eqn:bcins}
\end{equation}
In \eqref{eqn:bco}--\eqref{eqn:bcins}, $\bn$ is the outer unit normal.
Observe that the boundary conditions \eqref{eqn:bco}--\eqref{eqn:bcins} can be combined into
\begin{equation*}
\left. (\bn \cdot \nabla w + q(\bx, \balpha) w) \right|_{\partial \Omega} = g(\bx, \balpha),
\end{equation*}
for the appropriate choices of $q(\bx, \balpha)$ and $g(\bx, \balpha)$ defined on $\partial \Omega$ with
$\balpha \in \cA$, a parameter domain that we take to be the 4D box $\cA = [0.01, 0.5] \times [0, 0.9]^3$.
The initial temperature is taken to be zero throughout $\Omega$.

The system \eqref{eqn:heat}--\eqref{eqn:bcins} is discretized with $P_2$ finite elements on a quasi-uniform 
triangulation of $\Omega$ resulting in $M = 3,562$ spatial degrees of freedom. The choice of standard nodal 
basis functions $\{ \theta_j(\bx) \}_{j=1}^M$  defines the mass
$\rM \in \mathbb{R}^{M \times M}$ and stiffness $\rK \in \mathbb{R}^{M \times M}$ matrices,
 as well as boundary terms $\rQ(\balpha) \in \mathbb{R}^{M \times M}$,
$\bg(\balpha) \in \mathbb{R}^{M}$ with entries given by
\begin{equation*}
(\rQ)_{ij}(\balpha) = \int_{\partial \Omega} q(\bx, \balpha) \theta_j(\bx) \theta_i(\bx) ds_\bx, \quad 
(\bg)_{j}(\balpha) = \int_{\partial \Omega} g(\bx, \balpha) \theta_j(\bx) ds_\bx,
\end{equation*}
for $i,j = 1,\ldots,M$.

The vector-valued function of nodal values $\bu(t, \balpha) : [0,T) \times \cA \to \mathbb{R}^M$ solves
\begin{equation}
\rM \bu_t + \left( \rK + \rQ(\balpha) \right) \bu = \bg(\balpha),
\label{eqn:systemu}
\end{equation}
i.e., it satisfies the dynamical system of the form \eqref{eqn:GenericPDE} with 
\begin{equation*}
F(t, \bu, \balpha) = - \rM^{-1} \left( \rK + \rQ(\balpha) \right) \bu + \rM^{-1} \bg(\balpha)
\end{equation*}
and the initial condition $\bu (0, \balpha) = \bu_0 = \boldsymbol{0} \in \mathbb{R}^M$
corresponding to zero initial temperature condition for $w$.

\begin{figure}[ht]
\begin{center}
\includegraphics[width=0.6\textwidth]{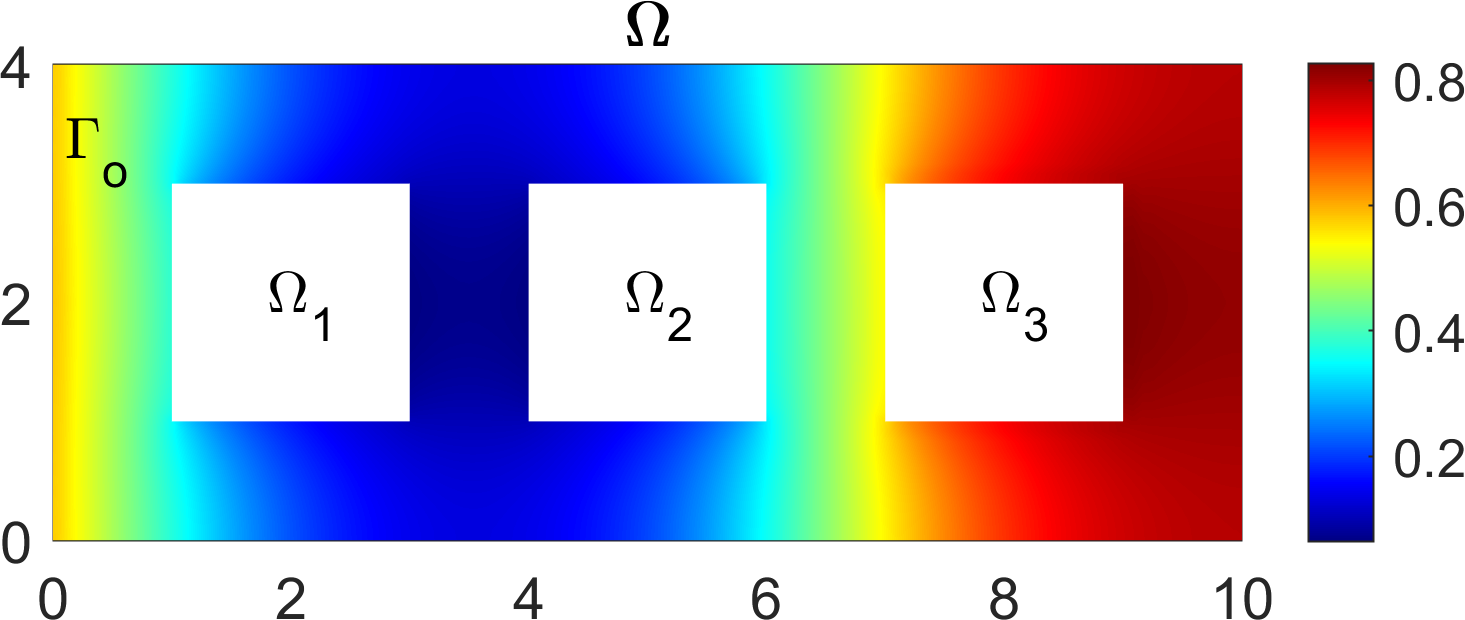} 
\end{center}
\caption{Domain $\Omega$ and the solution $w(T, \bx, \balpha)$ of the heat equation 
\eqref{eqn:heat}--\eqref{eqn:bcins} corresponding to $\balpha = (0.5, 0, 0, 0.9)^T$.}
\label{fig:domain}
\end{figure}

We compute the snapshots $\bphi_k = \bu(t_k, \balpha)$ by time-stepping \eqref{eqn:systemu}
at $t_k = 0.2 k$, $k = 1,2,\ldots,N$, with $N = 100$ time steps and $T = 20$ using
Crank-Nicolson scheme. 
Setting $\Theta(\bx) = [\theta_1(\bx), \ldots, \theta_M(\bx)]$ allows to express
the solution $w(t, \bx, \balpha)$ of \eqref{eqn:heat}--\eqref{eqn:bcins} as 
$w(t, \bx, \balpha) = \Theta(\bx) \bu(t, \balpha)$, hence the solution snapshots are
\begin{equation}
w(t_k, \bx, \balpha) = \Theta(\bx) \bu(t_k, \balpha).
\label{eqn:wsnap}
\end{equation}
The setting is illustrated in Figure~\ref{fig:domain}, where we display the domain $\Omega$ 
along with solution $w(T, \bx, \balpha)$ corresponding to parameter values $\balpha = (0.5, 0, 0, 0.9)^T$.

For an arbitrary but fixed $\balpha\in\cA$ let $\rZ = [\bz_1, \ldots, \bz_n] \in \mathbb{R}^{M \times n}$ be a matrix with columns being vectors constituting the reduced basis, i.e. $\rZ=\rU\rU_c$ for TROM. Then, the projection ROM of \eqref{eqn:systemu} is
\begin{equation}
\widetilde{\rM} \widetilde{\bu}_t + \left( \widetilde{\rK} + \widetilde{\rQ}(\balpha) \right) 
\widetilde{\bu} = \widetilde{\bg}(\balpha),
\label{eqn:systemurom}
\end{equation}
where 
\begin{align*}
\widetilde{\rM} & = \rZ^T \rM \rZ \in \R^{n \times n}, \quad
\widetilde{\rK} = \rZ^T \rK \rZ \in \R^{n \times n}, \\ 
\widetilde{\rQ}(\balpha) & = \rZ^T \rQ (\balpha) \rZ \in \R^{n \times n}, \quad
\widetilde{\bg}(\balpha) = \rZ^T \bg(\balpha) \in \R^{n}, 
\end{align*}
and the initial condition is $\widetilde{\bu}(0, \balpha) = \rZ^T \bu_0 = \boldsymbol{0} \in \R^n$. 
As discussed in Section~\ref{rem1}, the evaluation of \eqref{eqn:systemurom} can be effectively split 
between the offline and online stages.
Solving \eqref{eqn:systemurom} for $\widetilde{\bu}(t, \balpha)$ allows to recover the approximate 
solution at times $t_k$ as
\begin{equation}
\widetilde{w}(t_k, \bx, \balpha) = \Theta(\bx) \; \rZ \; \widetilde{\bu}(t_k, \balpha) \approx w(t_k, \bx, \balpha).
\label{eqn:wromsnap}
\end{equation}

\subsubsection{In-sample prediction and compression study}

We begin TROM assessment with in-sample prediction and compression study for the linear parabolic
system described Section~\ref{sec:linpde}. To measure TROM predictive power and to compare it to 
that of POD-ROM, we sample $\cA$ uniformly in each direction with 
$n_1 \times n_2 \times n_3 \times n_4 = 9 \times 5 \times 5 \times 5$ samples, for a total of 
$K = 1,125$ samples in the set $\hcA = \{ \widehat{\balpha}_1, \dots, \widehat{\balpha}_K \}$.
For each of the three TROMs and for POD-ROM we compute the following in-sample prediction error
\begin{equation}
E_{L^2(\hcA)} =  \left( \frac{1}{M N K} \sum_{j=1}^{K} 
\left\| (\rI - \rZ \rZ^T) \Phi_e (\widehat{\balpha}_j) \right\|_F^2 \right)^{1/2},
\label{eqn:epred}
\end{equation}
to quantify the ability of the  CP-TROM, HOSVD-TROM, TT-TROM local bases and POD-ROM basis to 
represent original snapshots for in-sample parameter values.
Note that $E_{L^2(\hcA)}^2$ is the quantity from \eqref{eqn:phibound} averaged over $\hcA$ and scaled by $(MN)^{-1}$.    

\addtolength{\tabcolsep}{-3pt}  

\begin{table}[h!]
\caption{Cartesian grid-based sampling in-sample prediction and compression study results
reporting prediction errors $E_{L^2(\hcA)}$ for ${\rm HOSVD},~{\rm TT},~{\rm POD}$, 
the number of compressed tensors elements transmitted to the online stage, 
and the corresponding compression factors CF.
$\widetilde{\bPhi}$ ranks comprise: Tucker ranks for HOSVD-TROM in format 
$[\widetilde{M}, \widetilde{n}_1, \widetilde{n}_2, \widetilde{n}_3, \widetilde{n}_4, \widetilde{N} ]$, and
compression ranks for TT-TROM in format
$[\widetilde{r}_1, \widetilde{r}_2, \widetilde{r}_3, \widetilde{r}_4, \widetilde{r}_5]$.}
\label{tab:D4C}
\small
\begin{center}
\begin{tabular}{|l|c|c|c|c|c|c|}
\hline
\multicolumn{2}{|c|}{$\widetilde{\eps}$} & 1e-4 & 1e-5 & 1e-6 & 1e-7 & 1e-9 \\
\hline
\multicolumn{2}{|c|}{$n$} & 12 & 16 & 19 & 23 & 30 \\
\hline
                   & HOSVD & 3.45e-05 & 3.85e-06 & 2.78e-07 & 5.61e-08 & -- \\
$E_{L^2(\hcA)}$ & TT   & 6.15e-05 & 5.58e-06 & 5.21e-07 & 5.03e-08 & 5.41e-10 \\
                   & POD & 1.67e-03 & 9.06e-04 & 5.84e-04 & 3.05e-04 & 7.76e-05 \\
\hline
         & HOSVD & $[34, 4, 2,$  & $[46, 5, 2,$  & $[57, 6, 2,$  & $[66, 7, 2,$  & -- \\
$\widetilde{\bPhi}$ ranks  &             & $2, 2, 12]$  & $2, 2, 16]$ & $2, 2, 20]$ & $2, 2, 23]$ & \\
  & TT        & $[34, 35, 30,$  & $[46, 48, 41,$  & $[57, 61, 51,$  & $[69, 74, 60,$  & $[99, 97, 80,$  \\
         &              & $21, 12]$ & $29, 16]$ & $36, 19]$ & $43, 23]$ & $57, 30]$ \\
\hline
  & HOSVD &  13122 & 29515 & 54804 & 85101 & -- \\
$\#\mbox{online}(\widetilde{\bPhi})$  & TT    &  20382 & 37993 & 53877 & 86022 & 156607 \\
\hline
    & HOSVD &  3.05e+4 & 1.36e+4 & 7.31e+3 & 4.71e+3 & -- \\
CF  & TT    &  1.97e+4 & 1.05e+4 & 7.42e+3 & 4.66e+3 & 2.56e+3 \\
\hline
\end{tabular}
\end{center}
\end{table}

\begin{table}[h!]
\caption{General parameter sampling in-sample prediction and compression study results
reporting prediction errors $E_{L^2(\hcA)}$ for ${\rm HOSVD},~{\rm TT},~{\rm POD}$, 
the number of compressed tensors elements transmitted to the online stage, 
and the corresponding compression factors CF.
$\widetilde{\bPhi}$ ranks comprise: Tucker ranks for HOSVD-TROM in format 
$[\widetilde{M}, \widetilde{N} ]$, and
compression ranks for TT-TROM in format
$[\widetilde{r}_1, \widetilde{r}_2]$.}
\label{tab:D4G}
\small
\begin{center}
\begin{tabular}{|l|c|c|c|c|c|c|}
\hline
\multicolumn{2}{|c|}{$\widetilde{\eps}$} & 1e-4 & 1e-5 & 1e-6 & 1e-7 & 1e-9 \\
\hline
\multicolumn{2}{|c|}{$n$} & 12 & 16 & 19 & 23 & 30 \\
\hline
                   & HOSVD & 5.69e-05 & 5.29e-06 & 4.66e-07 & 7.93e-08 & -- \\
$E_{L^2(\hcA)}$ & TT       & 5.63e-05 & 5.95e-06 & 5.72e-07 & 5.56e-08 & 5.43e-09 \\
                   & POD     & 2.05e-03 & 1.12e-03 & 5.84e-04 & 3.05e-04 & 8.46e-05 \\
\hline
$\widetilde{\bPhi}$ ranks  & HOSVD & $[33, 11]$ & $[45, 16]$ & $[56, 19]$ & $[65, 23]$ & -- \\
 & TT      & $[32, 11]$ & $[43, 15]$ & $[55, 19]$ & $[66, 23]$ & $[95, 29]$ \\
\hline
  & HOSVD &  11904 & 18450 & 26268 & 36680 & -- \\
$\#\mbox{online}(\widetilde{\bPhi})$	 & TT    &  396011 & 725640 & 1175644 & 1633522 & 3099404 \\
\hline
                         & HOSVD & 3.37e+4 & 2.17e+4 & 1.53e+4 & 1.09e+4 & -- \\
CF 						 & TT    & 1.01e+3 & 5.52e+2 & 3.41e+2 & 2.45e+2 & 1.29e+2 \\
\hline
\end{tabular}
\end{center}
\end{table}

\begin{figure}[ht]
\begin{center}
\includegraphics[width=0.48\textwidth]{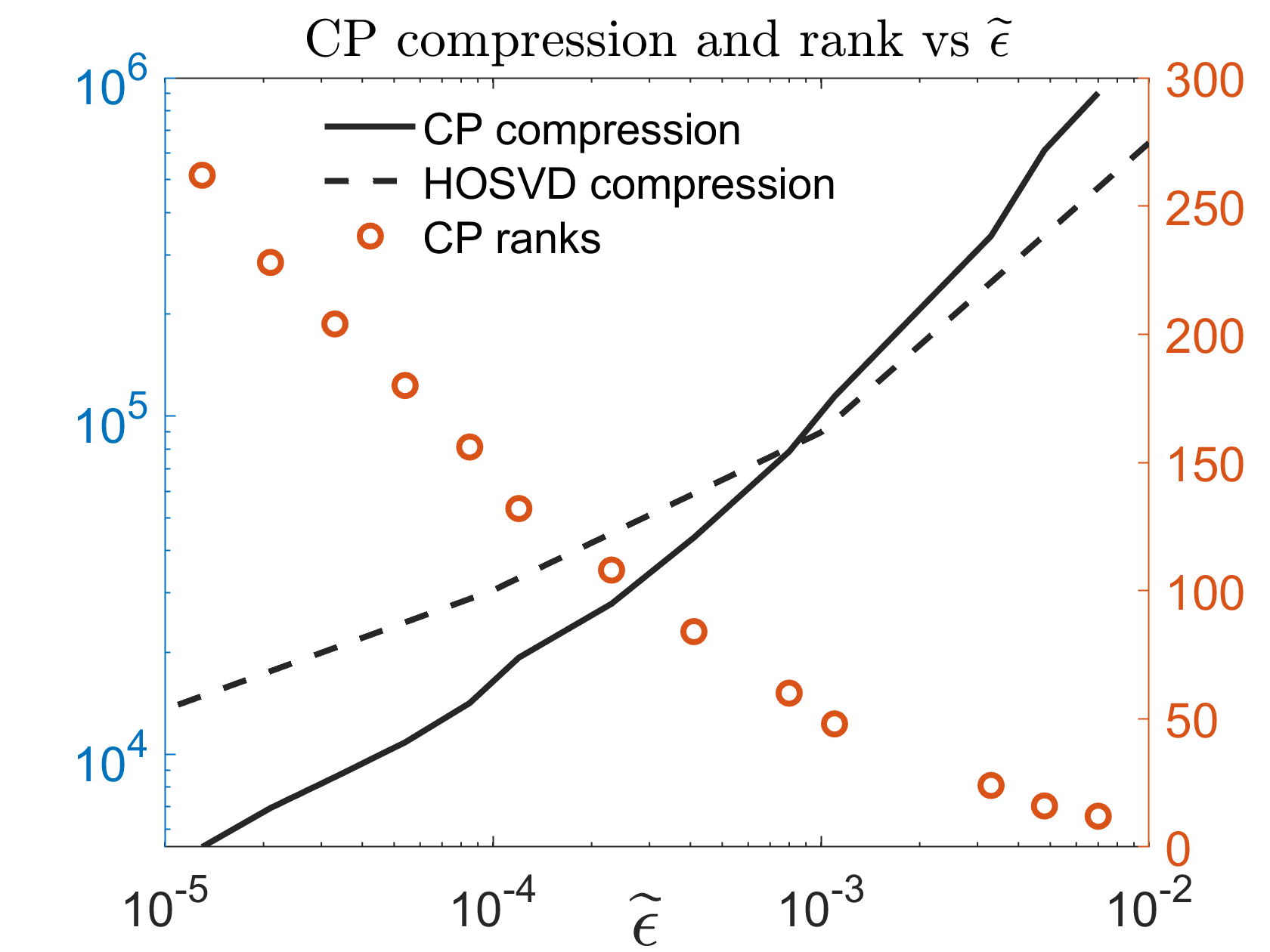} 
\includegraphics[width=0.48\textwidth]{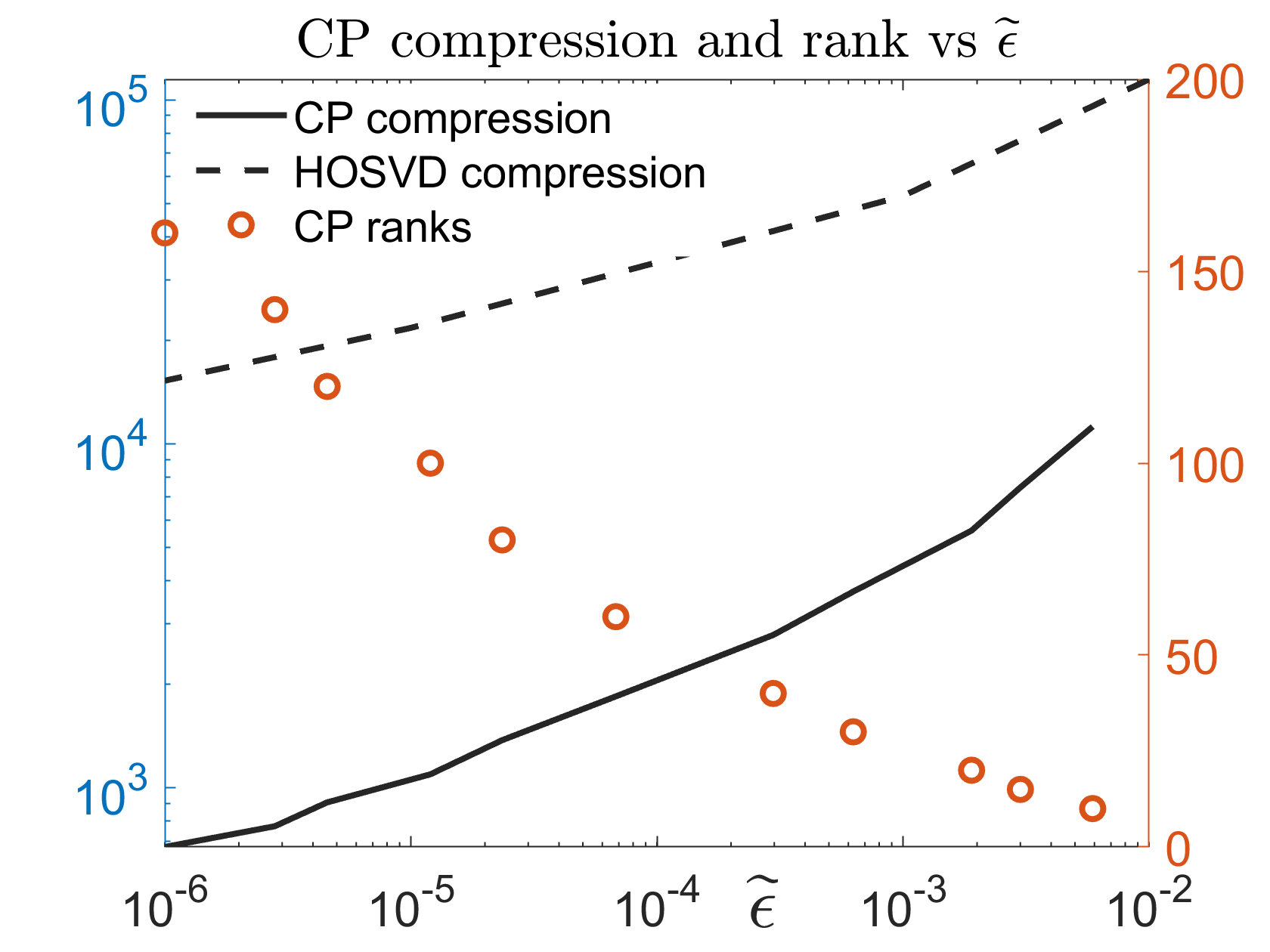} 
\end{center}
\caption{CP-TROM compression factor CF and canonical rank $R$ as functions of $\widetilde{\eps}$. 
Left: Cartesian grid-based sampling; Right: general sampling.}
\label{fig:CPranks}
\end{figure}

We report in Tables~\ref{tab:D4C} and \ref{tab:D4G} the in-sample prediction errors and compression 
factors CF defined in \eqref{eqn:CF} for Cartesian grid-based and general samplings, respectively.
For general sampling we organize the sampling parameters from the Cartesian grid in 1D array leading to 
snapshot tensors or order 3.  
Experimenting with the same number of randomly sampled parameters showed very similar compression 
rates and in-sample prediction errors and so those are not reported here.   
The results are reported in Tables~\ref{tab:D4C} and \ref{tab:D4G}  for a number of  decreasing values 
of $\widetilde{\eps}$ and correspondingly increasing $n$, such that $n \leq \min(\widetilde{N}, \widetilde{r}_{D+1})$, 
where $\widetilde{N}$ and $\widetilde{r}_{D+1}$ are the last Tucker and compression ranks, respectively,
for HOSVD- and TT-TROM. Available compression algorithm failed to deliver the accuracy of 1e-9 for HOSVD, 
so we report only TT statistics for this extreme value of $\widetilde{\eps}$.
Note that we leave $E_{L^2(\hcA)}$ for CP-ROM out since there is no direct way to control its relative error 
$\widetilde{\eps}$, as discussed at the end of Section~\ref{sec:CP-TROM}. Instead, 
compression factors and canonical ranks for CP-TROM are illustrated in Figure~\ref{fig:CPranks}.

We observe in Tables~\ref{tab:D4C} and \ref{tab:D4G} that both HOSVD- and TT-TROM outperform 
POD-ROM in terms of prediction error up to four orders of magnitude for larger $n$ and correspondingly 
small $\widetilde{\eps}$. This result is consistent across both Cartesian grid-based and general 
parameter samplings. The $\#\mbox{online}(\widetilde{\bPhi})$ values tell us that for Cartesian based 
sampling,  HOSVD and TT are comparable in terms of memory and data transmission requirements with HOSVD 
doing somewhat better for lower representation accuracy for the snapshot tensor of order 6. 
Compression achieved varies in $\widetilde{\eps}$ (as should expected) and gives more than 3 orders of saving 
even for the finest available representation accuracy. If the Cartesian structure of $\widehat{\cA}$ is abandoned 
and snapshots are organized in tensors of order 3, then HOSVD format has a clear advantage over TT in terms 
of compression achieved; see $\#\mbox{online}(\widetilde{\bPhi})$ and CF statistics  in Table~\ref{tab:D4G}. 

For CP-TROM we display in Figure~\ref{fig:CPranks} compression factors CF and canonical ranks $R$
for a number of values of $\widetilde{\eps}$. For comparison we also show on the same plots CF for HOSVD, 
which shows that HOSVD-TROM is on par with CP-TROM if Cartesian sampling allows to organize snapshots in a 
higher order tensor. However, we \rev{were} able to compute $\widetilde{\bPhi}$ in CP compressed format only up to 
moderate values of $\widetilde{\eps}$, since the corresponding CP rank was growing fast as $\widetilde{\eps}$ 
decreases (see the left plot in Figure~\ref{fig:CPranks}). Smaller $\widetilde{\eps}$ become feasible with CP 
if the snapshots are organized in tensors of order 3, but in this case HOSVD-TROM achieves much better compression 
than CP-TROM  (see the right plot in Figure~\ref{fig:CPranks}). We conclude that for this example with a relatively 
small number of parameters ($D=4$) HOSVD-TROM appears to be the best performing TROM. We finally note that 
for the same compression accuracy $\widetilde{\eps}$ \rev{(if it was achieved) CP-TROM demonstrated} very similar in-sample
prediction error \rev{as HOSVD- and TT-ROMs}. This observation largely carries over to out-of-sample representation studied next.

\subsubsection{Out-of-sample prediction study}

To quantify the ability of the CP-, HOSVD- and TT-TROM local bases to represent the solution of \eqref{eqn:heat} 
for arbitrary out-of-sample parameter values, we use $E_{L^2(\cA)}$ which is defined as in \eqref{eqn:epred} 
but with $\balpha$ (in place of $\hbalpha$) running through a large number of random points from $\cA$. 
We also use 
\begin{equation*}
	E_{L^\infty(\cA)} = \sup_{\balpha\in\cA} \left( \frac{1}{M N}  
	\left\| (\rI - \rZ \rZ^T) \Phi_e ({\balpha}) \right\|_F^2 \right)^{1/2},
\end{equation*}
for the maximum of representation error over the parameter domain. An estimate of $E_{L^\infty(\cA)}$ is given by
Theorem~\ref{Th1}. \rev{Regarding constants appearing in \eqref{eqn:est1} we note that for our choice of the uniform grid in $\cA$ and $p=2,3$ one computes $C_e=1$, $C_a=\frac18$ ($p=2$) and $C_a=\frac1{48}$ ($p=2$), while constant $C_{\bu}$ is hard to evaluate. Experimentally we found that the grid in the forth parameter should be finer than for the first three parameters to balance the observed error suggesting that the solution is less smooth as a function of $\alpha_4$.} To study the error dependence on the parameter mesh size $\delta$, we reduce the number 
of parameters to two, letting $\alpha_1=\alpha_2=\alpha_3$ in \eqref{eqn:bcj}.
\begin{figure}[ht]
\begin{center}
\begin{tabular}{ccc}
\tiny TROM and POD prediction error vs $n$ & 
\tiny TROM prediction error vs $\delta$ & 
\tiny TROM prediction error vs $\widetilde{\eps}$ \\
\includegraphics[width=0.317\textwidth]{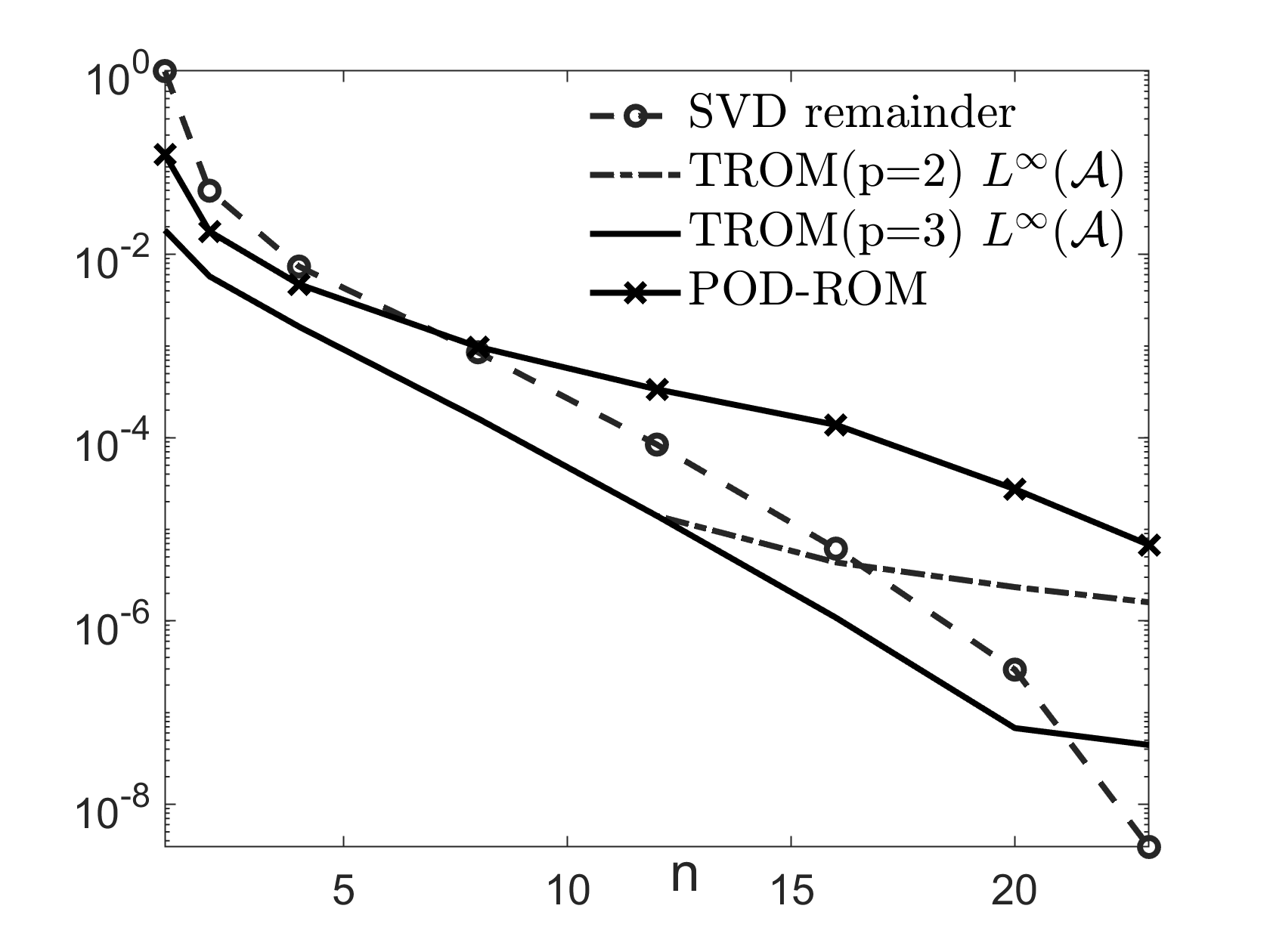} &
\includegraphics[width=0.317\textwidth]{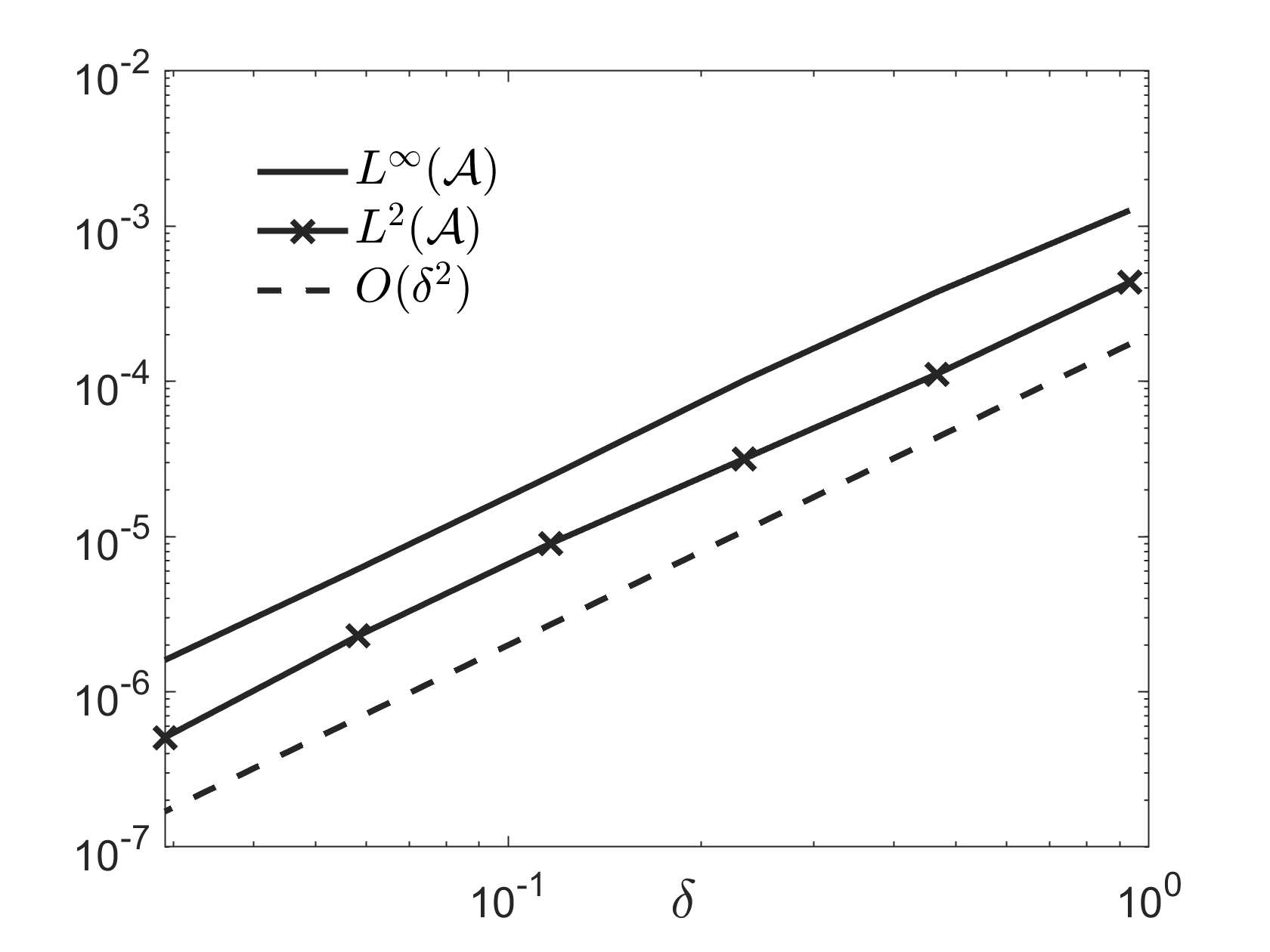} &
\includegraphics[width=0.317\textwidth]{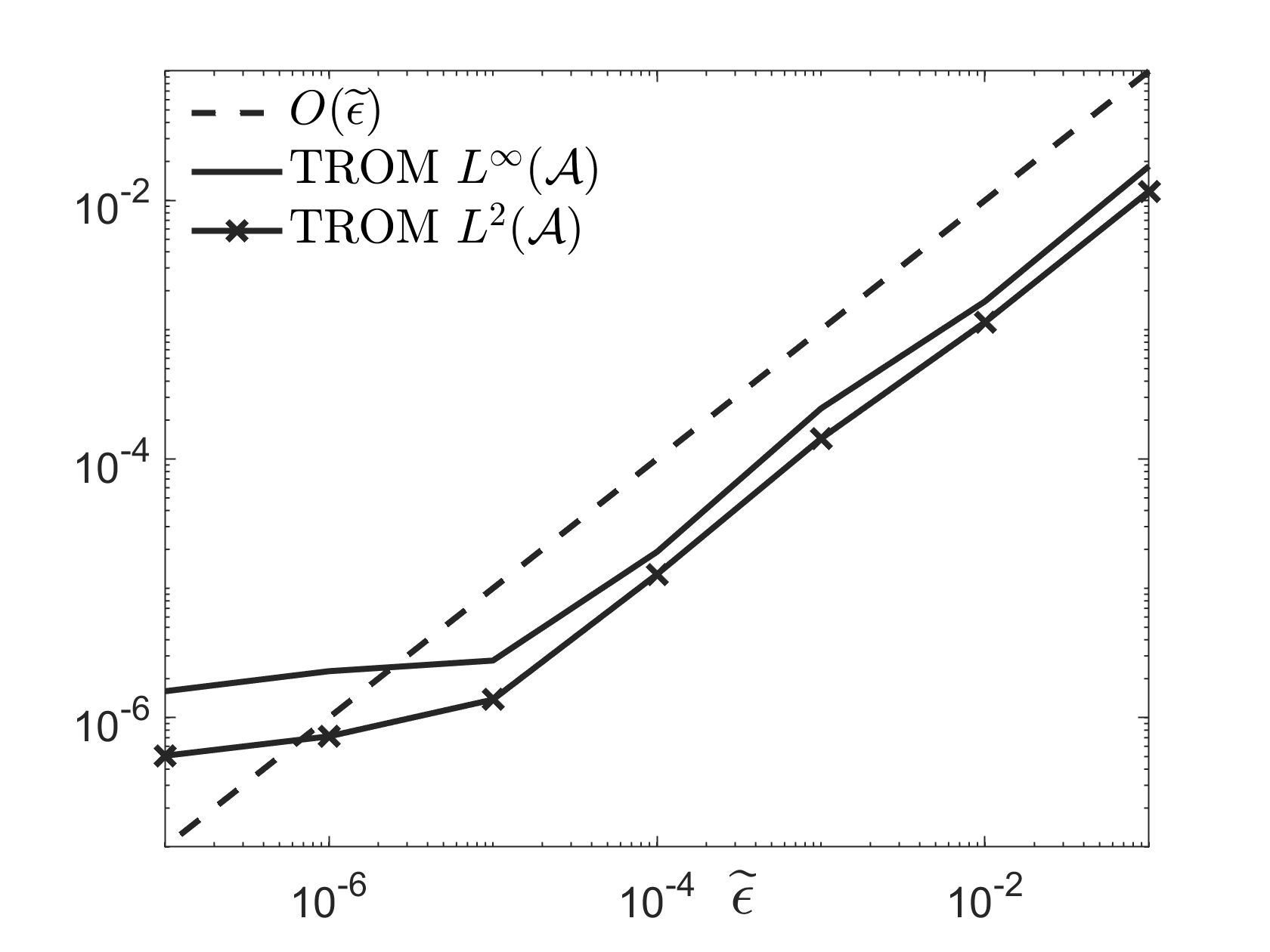} 
\end{tabular}
\end{center}
\caption{Out-of-sample prediction error of TROM versus ROM basis dimension $n$ (left plot), parameter mesh size $\delta$ (middle plot) and tensor compression accuracy $\widetilde{\eps}$ (right plot). }
\label{fig:predict}
\end{figure}

In Figure~\ref{fig:predict} we plot $E_{L^\infty(\cA)}$ and $E_{L^2(\cA)}$ versus ROM basis dimension $n$, 
parameter mesh size $\delta$ and tensor compression accuracy $\widetilde{\eps}$. The results were computed 
using 100 randomly distributed parameters from $\cA$ to evaluate the error quantities for HOSVD-TROM. 
For TT-TROM and CP-TROM the error dependence on  $n$, $\delta$ and $\widetilde{\eps}$ was virtually the 
same and are omitted. Variants of TROM of course may differ by complexity. \rev{While the cost of offline phase of computing $\widetilde{\bPhi}$ varies significantly depending on the format, the online distribution of costs was persistent for all three formats showing between 30\% and 40\% of online time spent on finding $U_c$ (coordinates of local basis), less then 5\% on projection to $\balpha$-specific coordinates, and about 60\% of  time on the integration of the projected system (this last step is common with standard POD approach).} 

The left plot in Figure~\ref{fig:predict} shows the $E_{L^\infty(\cA)}$ error for different values of $n$ and 
compares it to the error of POD-ROM. We use $\widetilde{\eps}$=1e-7 and $65\times33$ parameter grid. 
Such fine compression accuracy and grid allows us to isolate the effect of the second term on the right hand 
side of \eqref{eqn:est1}. Indeed, we see that the error curve follows closely the graph of ``SVD remainder'' 
(the maximum over all $\balpha$ of the second term on the right-hand side of \eqref{eqn:est1}) until the interpolation error starts 
dominating for larger $n$. As can be expected, the interpolation error for $p=3$ starts affecting 
$E_{L^\infty(\cA)}$ for larger $n$ than the interpolation error for $p=2$. 

The middle plot in Figure~\ref{fig:predict} demonstrates that both $E_{L^\infty(\cA)}$ and $E_{L^2(\cA)}$ 
for TROM decrease as $O(\delta^2)$ for $p=2$  (computed with $\widetilde{\eps}$=1e-7 and $n=\widetilde{N}$) 
just as predicted by \eqref{eqn:est1}. Likewise, the right plot in Figure~\ref{fig:predict} gives evidence for 
$O(\widetilde{\eps})$ decrease of $E_{L^\infty(\cA)}$ and $E_{L^2(\cA)}$  (computed with $65\times33$ 
parameter grid and $n=\widetilde{N}$) in accordance to \eqref{eqn:est1} as long the first term on the right 
hand side of  \eqref{eqn:est1} dominates. 

\subsubsection{Out-of-sample TROMs  vs POD-ROM  performance: heat equation}
\label{sec:oosheat}

Performance of TROM and POD-ROM may vary for different out-of-sample parameter values. 
Therefore, assessment of TROM vs POD-ROM performance in this section is conducted in a statistical setting. 
The quantities of interest are computed for $N_r \gg 1$ out-of -sample realizations of 
$\balpha^{(r)} \in \cA$, $r=1,2,\ldots,N_r$, where we use $N_r = 200$ 
for the numerical studies below. Realizations 
$\balpha^{(r)} = \left( \alpha_1^{(r)}, \alpha_2^{(r)}, \alpha_3^{(r)}, \alpha_4^{(r)} \right)^T$ 
are drawn at random from $\cA = [0.01, 0.5] \times [0, 0.9]^3$ with each $\alpha_i$ 
distributed uniformly on $[\alpha_i^{\min}, \alpha_i^{\max}]$, $i=1,2,3,4$.

For statistical tests we use the following quantities to measure performance of TROM.
First, we introduce the relative $L^{\infty}(0, T, L^2(\Omega))$ ROM solution error
\begin{equation}
	\begin{split}
R_X (\balpha) & = 
\frac{\max\limits_{k = 1,\ldots,N} 
\left\| \widetilde{w}(t_k, \bx, \balpha) - w(t_k, \bx, \balpha) \right\|_{L^2(\Omega)}}
{\max\limits_{k=1,\ldots,N} 
\left\| w(t_k, \bx, \balpha) \right\|_{L^2(\Omega)}}  \\
& \approx
\frac{\sup\limits_{t \in [0,T]} 
\left\| \widetilde{w}(t, \bx, \balpha) - w(t, \bx, \balpha) \right\|_{L^2(\Omega)}}
{\sup\limits_{t \in [0,T]} \left\| w(t, \bx, \balpha) \right\|_{L^2(\Omega)}},
\end{split}
\label{eqn:Rx}
\end{equation}
which we compute for each realization $\balpha^{(r)}$, $r=1,2,\ldots,N_r$, 
for both POD-ROM and each of the three TROMs with X$\in\{$POD, CP, HOSVD, TT$\}$.
The true and reduced order snapshots for \eqref{eqn:Rx} are computed as in
\eqref{eqn:wsnap} and \eqref{eqn:wromsnap}, respectively. We report the mean, 
minimum,  and standard deviation of the three relative \emph{gain} distributions
\begin{equation}
G_{\text{X}}^{(r)} = 
\frac{R_{\text{POD}} \big( \balpha^{(r)} \big)}
{R_{\text{X}} \big( \balpha^{(r)} \big)},
\quad r=1,2,\ldots,N_r,
\label{eqn:relgain}
\end{equation}
for X$\in\{$CP, HOSVD, TT$\}$, which quantify the error decrease of CP-TROM, HOSVD-TROM and TT-TROM, 
respectively, relative to POD-ROM. We study the dependency of \eqref{eqn:relgain} with respect to 
$K$, the number of sampled parameter values in $\cA$, and $n$, the dimension of the reduced space. 
The results are reported for both Cartesian grid-based sampling and general parameter sampling.

\begin{table}[h!]
\caption{Statistics of relative gain \eqref{eqn:relgain} for various values of $K$, 
the number of sampled parameter values in $\cA$. 
The study is performed with $n = 10$.}
\label{tab:relgaink}
\small
\begin{center}
\begin{tabular}{|c|c|c|c|c|c|c|c|}
\hline
\multicolumn{2}{|c|}{} & 
\multicolumn{3}{c|}{General parameter sampling} & \multicolumn{3}{c|}{Cartesian grid-based sampling} \\
\hline
$K$ & $G_X$ & {CP} & {HOSVD} & {TT}  & {CP} & {HOSVD}  & {TT} \\
\hline
                           & mean & \bf 24.17 & \bf 24.17 & \bf 24.17 & \bf 24.76 & \bf 25.08 & \bf 25.08 \\
$135 = $              & min   & 0.49  & 0.49 & 0.49 & 0.56 & 0.56 & 0.56 \\
$5 \times 3^3 $    
                            & std   & 17.33 & 17.33 & 17.33 & 16.88 & 17.31 & 17.32 \\
\hline
                           & mean & \bf 34.60 & \bf 34.60 & \bf 34.61 & \bf 35.21 & \bf 35.52 & \bf 35.51 \\
$1000 = $            & min    & 2.80 & 2.80 & 2.80 & 1.72 & 1.72 & 1.72 \\
$8 \times 5^3 $   
                           & std     & 15.36 & 15.36 & 15.37 & 15.03 & 15.11 & 15.11 \\
\hline
                            & mean & \bf 38.14 & \bf 38.15 & \bf 38.15 & \bf 37.80 & \bf 38.80 & \bf 38.80 \\
$3430 = $                   & min    & 3.81 & 3.81 & 3.81 & 4.20 & 4.45 & 4.43 \\
$10 \times 7^3 $           
                            & std     & 14.20 & 14.21 & 14.22 & 12.96 & 13.61 & 13.62 \\
\hline
\end{tabular}
\end{center}
\end{table} 

\begin{table}[h!]
\caption{Statistics of relative gain \eqref{eqn:relgain} for $n = 10$ and $20$. 
The study is performed for $K = 10 \times 7 \times 7 \times 7 = 3430$.}
\label{tab:relgainn}
\small
\begin{center}
\begin{tabular}{|c|c|c|c|c|c|c|c|}
\hline
\multicolumn{2}{|c|}{} & 
\multicolumn{3}{c|}{General parameter sampling} & \multicolumn{3}{c|}{Cartesian grid-based sampling} \\
\hline
$n$ & $G_X$ & {CP} & {HOSVD} & {TT}  & {CP} & {HOSVD}  & {TT} \\
\hline
                            & mean & \bf 38.14 & \bf 38.15 & \bf 38.15 & \bf 37.80 & \bf 38.80 & \bf 38.80 \\
$10$                       & min    & 3.81 & 3.81 & 3.81 & 4.20 & 4.45 & 4.43 \\
                            & std     & 14.20 & 14.21 & 14.22 & 12.96 & 13.61 & 13.62 \\
\hline
        & mean & \bf 155.00 & \bf 158.97 & \bf 161.75 & \bf 49.80 & \bf 155.65 & \bf 154.03 \\
$20$ & min    & 4.59 & 4.54 & 4.55 & 1.51 & 5.26 & 5.23 \\
        & std     & 513.33 & 530.50 & 557.86 & 39.48 & 551.92 & 541.88 \\
\hline
\end{tabular}
\end{center}
\end{table} 

We present in Table~\ref{tab:relgaink} the dependence of relative gain statistics on the values of $K$
 (statistics in the table were computed setting $\widetilde{\eps} = 10^{-5}$ as the targeted accuracy 
of HOSVD-TROM, TT-TROM, and $R=250$ as the targeted rank for CP-TROM). 
We observe that as $K$ increases, TROMs become both more accurate on average and more robust.
The robustness is observed in both the increase of the minimum relative gain and the decrease of 
its standard deviation. On average, for $K = 3,430$, all three TROMs are almost $40$ times more 
accurate compared to POD-ROM. The performance difference between the TROMs themselves is 
basically negligible for this particular study. 

The performance of TROMs in the example above is limited by the relatively small value of $n = 10$. 
The effect of increasing $n$ to $20$ while keeping $K = 3,430$ is shown in Table~\ref{tab:relgainn} 
(statistics in the table were computed setting $\widetilde{\eps} = 10^{-7}$ as the targeted accuracy 
of HOSVD-TROM, TT-TROM, and $R=250$ as the targeted rank for CP-TROM).
While the worst case scenario stays relatively unchanged, the average accuracy gain by TROMs is 
over two orders of magnitude. 
An outlier here is CP-TROM that underperforms HOSVD- and TT-TROM in case of Cartesian grid-based
parameter sampling. Aside from that, performance difference of TROMs for general and Cartesian
samplings is negligible. \rev{It is also interesting to see if the interpolation of approximate snapshots from the 
lower-rank tensor $\widetilde{\bPhi}$ alone, i.e., without finding a reduced local basis and solving the projected 
problem, gives reasonable approximation to high-fidelity solutions for out-of-sample parameters. 
Such interpolation-only predicted solution for incoming $\balpha\in\cA$ is given by columns of $\widetilde{\Phi}(\balpha)$.
Repeating the experiment with Cartesian grid-based sampling and other parameters the same as used 
for results in Tables~\ref{tab:relgaink} and~\ref{eqn:relgain}, we find that the mean relative gain of 
HOSVD-tROM compared to interpolation-only approach is $\{2.10,\, 1.51,\, 1.13\}$ for $n=10$, 
$K=\{135, 1000, 3430\}$ and $\{7.64,\, 7.60,\, 7.54\}$ for $n=20$ and same values of $K$. 
The numbers were very close for other two TROMs. We see that TROM based on solving projected problem in 
general gives more accurate results then pure interpolation, especially if more vectors are included in the reduced basis.
If $n$ is fixed, then for sufficiently fine sampling the interpolation-only approach delivers the same (or even better) accuracy.}

\begin{table}[t!]
\caption{\rev{Statistics of relative gain \eqref{eqn:relgain} for POD-ROM with POD-Greedy reduced basis
for $n = 10$ and $20$. 
The study is performed for $K = 10 \times 7 \times 7 \times 7 = 3430$.}}
\label{tab:relgaing}
\small
\rev{
\begin{center}
\begin{tabular}{|c|c|c|c|c|c|c|c|}
\hline
\multicolumn{2}{|c|}{} & 
\multicolumn{3}{c|}{General parameter sampling} & \multicolumn{3}{c|}{Cartesian grid-based sampling} \\
\hline
$n$ & $G_X$ & {CP} & {HOSVD} & {TT}  & {CP} & {HOSVD}  & {TT} \\
\hline
                            & mean & \bf 53.14 & \bf 53.14 & \bf 53.14 & \bf 54.02 & \bf 54.12 & \bf 54.12 \\
$10$                       & min    & 7.01 & 7.01 & 7.01 & 7.43 & 7.42 & 7.42 \\
                            & std     & 18.63 & 18.63 & 18.63 & 18.03 & 18.03 & 18.03 \\
\hline
        & mean & \bf 188.57 & \bf 195.57 & \bf 198.77 & \bf 62.18 & \bf 193.66 & \bf 191.86 \\
$20$ & min    & 7.39 & 7.35 & 7.36 & 1.87 & 8.51 & 8.47 \\
        & std     & 571.50 & 607.44 & 634.94 & 49.07 & 649.81 & 639.61 \\
\hline
\end{tabular}
\end{center}
}
\end{table} 

We conclude the numerical study for the parameterized heat equation with a comparison between the three
TROMs and another variant of POD-ROM often used in practice, the so-called greedy POD or POD-Greedy approach
to computing the reduced basis \cite{patera2007reduced}. Replacing the conventional POD-ROM computation with 
POD-Greedy algorithm, we perform the same out-of-sample performance study as presented in
Table~\ref{tab:relgainn}. The resulting relative gain statistics are reported in Table~\ref{tab:relgaing}.
Qualitatively, the results are very similar to those in Table~\ref{tab:relgainn}. However, quantitatively we observe
$20\%$ to $40\%$ increased relative gain for all TROMs. This is consistent with the fact that POD-Greedy reduced basis
is sub-optimal compared to the conventional POD-ROM reduced basis computed from the snapshots corresponding to
all parameter values in $\hcA$.

In terms of computational performance for this specific setting, POD-Greedy algorithm was found to be significantly slower 
in the offline stage than all three TROM approaches and the conventional POD-ROM even if one includes into 
the offline cost of TROM and POD-ROM the computation of $KN$ snapshots $\bphi_k(\widehat{\balpha}_j)$, 
$k=1,\ldots,N$, $j=1,\ldots,K$. Indeed, the bulk of computational cost of POD-Greedy approach is in the evaluation 
of error estimator that has to be performed at each of its $n$ iterations for all $\widehat{\balpha}_j$,
$j=1,\ldots,K,$ to determine the sample with the largest error. In turn, each evaluation of the estimator 
requires the computation of the residual of the chosen time-stepping scheme (e.g., Crank-Nicolson) that needs
$O(N)$ matrix-vector products of a dense $M \times n$ matrix and a vector in $\mathbb{R}^n$.
Thus, error estimator evaluations alone account for $O(K N M n^2)$ operations of POD-Greedy offline
stage. {We note that in other situations, when computing the snapshopts is much more  expensive 
compared to the evaluation of the error estimator (e.g., when the matrices of systems to be solved on each time step are dense), 
both POD-ROM and TROM may benefit from a greedy approach to parameter sampling. }


\subsection{Advection-diffusion PDE}
\label{sec:advdiff}

In the second numerical example, for assessing performance of the three TROM approaches we are 
interested in a case with a higher order of parameter space compared to $D=4$ in Section~\ref{sec:linpde}.
To that end we set up a dynamical system resulting from the discretization of a linear advection-diffusion equation
\begin{equation}
w_t = \nu \Delta w - \bseta (\bx, \balpha) \cdot \nabla w + f(\bx),
\label{eqn:advdiff}
\end{equation}
in a unit square domain $\Omega = [0,1] \times [0,1] \subset \R^2$, $\bx = (x_1, x_2)^T \in \Omega$. 
Here $\nu$ is a constant diffusion coefficient, $\bseta: \Omega \times \cA \to \R^2$ is the advection field 
and $f(\bx)$ is a Gaussian source
\begin{equation}
f(\bx) = \frac{1}{2 \pi \sigma_s^2} \exp \left( - \cfrac{ (x_1 - x_1^s)^2 + (x_2 - x_2^s)^2}{2 \sigma_s^2} \right),
\label{eqn:advdiffsrc}
\end{equation}
where we take $\sigma_s = 0.05$, $x_1^s = \rev{x_2^s} = 0.25$. 
We enforce homogeneous Neumann boundary conditions and zero initial condition
\begin{equation}
\left. \left( \bn \cdot \nabla w \right) \right|_{\partial \Omega} = 0, \quad w(0, \bx, \balpha) = 0.
\label{eqn:advdiffbcic}
\end{equation}
The model is parametrized with $D = 9$ parameters with only the advection field $\bseta$ depending on 
$\balpha \in \R^9$. The advection field is given as follows
\begin{equation}
\bseta(\bx, \balpha) = \begin{pmatrix} \eta_1(\bx, \balpha) \\ \eta_2(\bx, \balpha) \end{pmatrix}  = 
\begin{pmatrix} \cos \alpha_9 \\ \sin \alpha_9 \end{pmatrix} 
+ \frac{1}{\pi} \begin{pmatrix} \partial_{x_2} h(\bx, \balpha) \\ - \partial_{x_1} h(\bx, \balpha) \end{pmatrix},
\label{eqn:etacos}
\end{equation}
where $h(\bx)$ is the cosine trigonometric polynomial
\begin{equation}
\begin{split}
h(\bx, \balpha) = & \;\;\;\; \alpha_1 \cos(\pi x_1) + \alpha_2 \cos(\pi x_2) + \alpha_3 \cos(\pi x_1) \cos(\pi x_2) \\
& + \alpha_4 \cos(2\pi x_1) + \alpha_5 \cos(2\pi x_2) + \alpha_6 \cos(2\pi x_1) \cos(\pi x_2) \\
& + \alpha_7 \cos(\pi x_1) \cos(2\pi x_2)  + \alpha_8 \cos(2\pi x_1) \cos(2\pi x_2).
\end{split}
\end{equation}
Here $\alpha_9$ determines the angle of the dominant advection direction, while parameters 
$\alpha_i$, $i=1,\ldots,8$, introduce perturbations into the advection field.
The parameter domain is a 9D box, see Section~\ref{sec:oosadvdiff} for details.

The system \eqref{eqn:advdiff}--\eqref{eqn:advdiffbcic} is discretized using $P_2$ finite elements on a grid 
with either $M = 1,893$ or $M = 4,797$ nodes (depending on a particular experiment) using the standard nodal 
basis functions $\{ \theta_j(\bx) \}_{j=1}^M$ that define the mass $\rM \in \mathbb{R}^{M \times M}$, 
stiffness $\rK \in \mathbb{R}^{M \times M}$ and advection $\rH(\balpha)\in \mathbb{R}^{M \times M}$ 
matrices, and the source vector $\blf \in \R^M$.
The vector-valued function of nodal values $\bu(t, \balpha) : [0,T) \times \cA \to \mathbb{R}^M$ solves
\begin{equation}
\rM \bu_t + \left( \rK + \rH(\balpha) \right) \bu = \blf,
\label{eqn:systemuad}
\end{equation}
i.e., it satisfies the dynamical system of the form \eqref{eqn:GenericPDE} with 
\begin{equation}
F(t, \bu, \balpha) = - \rM^{-1} \left(\rK + \rH(\balpha) \right) \bu 
+ \rM^{-1} \blf,
\end{equation}
and the initial condition $\bu (0, \balpha) = \boldsymbol{0} \in \mathbb{R}^M$.

\begin{figure}[ht]
\begin{center}
\begin{tabular}{cc}
Advection field $\bseta(\bx, \balpha)$ & Solution $w(T, \bx, \balpha)$ \\
\includegraphics[height=0.35\textwidth]{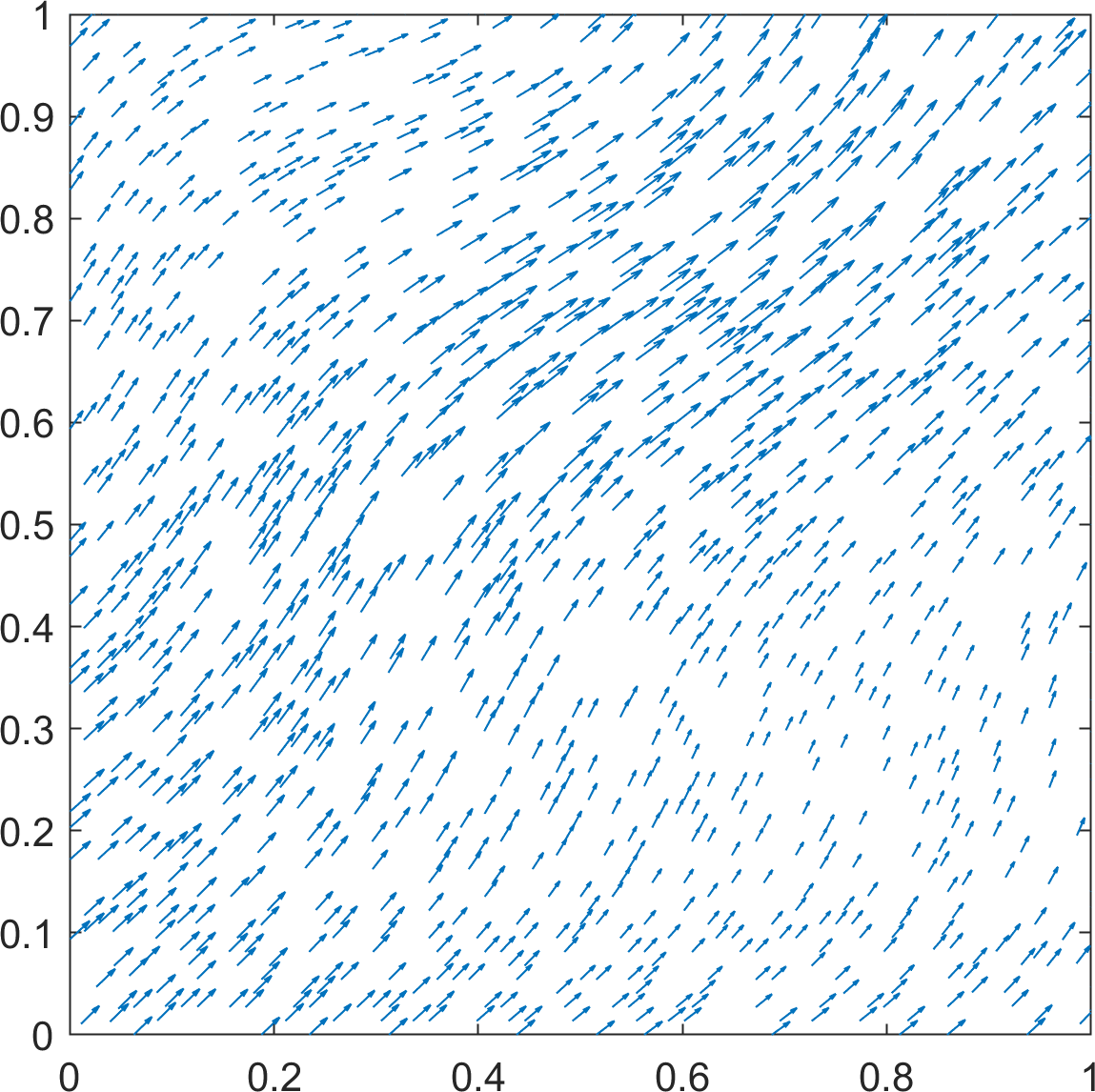} &
\qquad \includegraphics[height=0.35\textwidth]{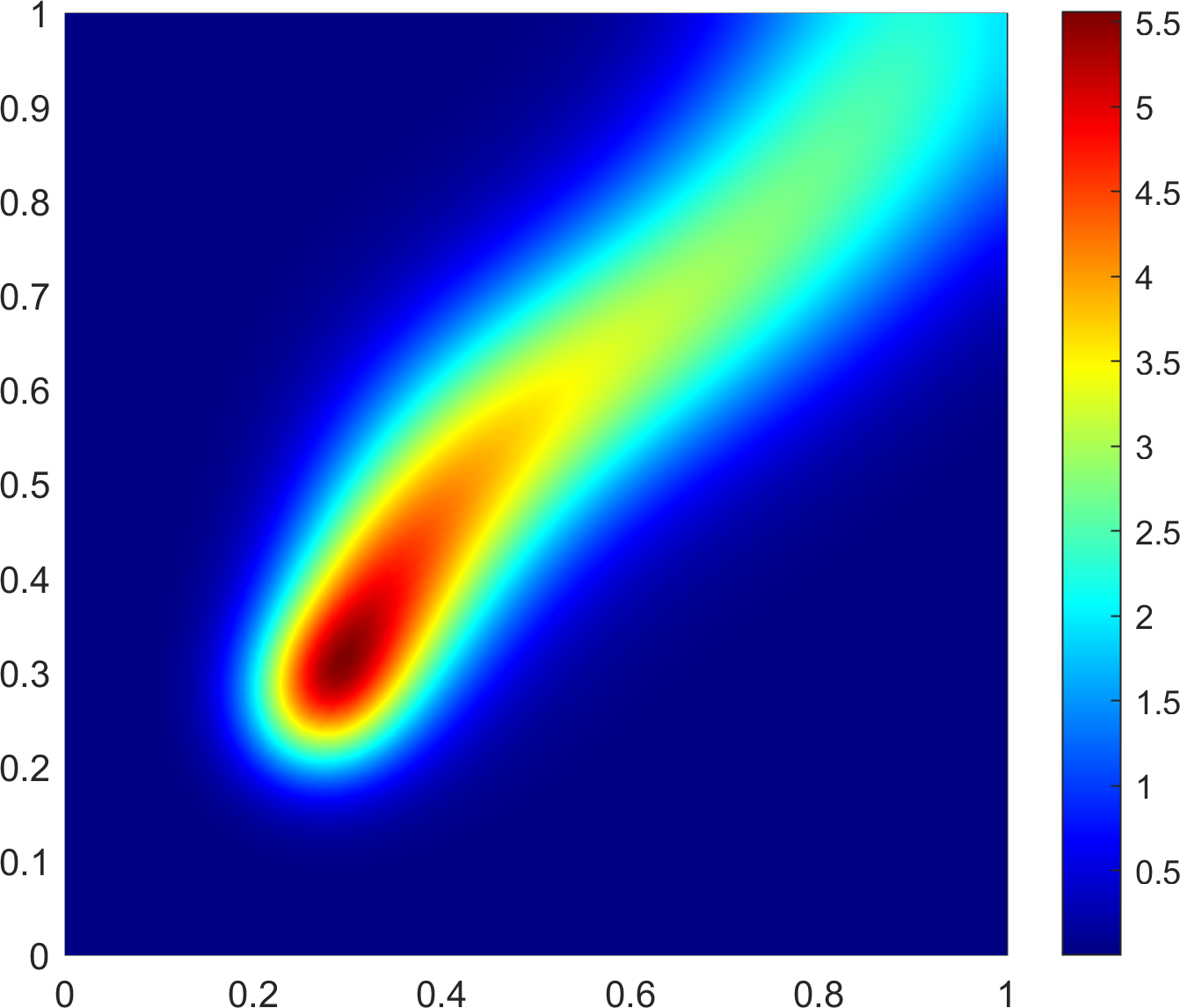} 
\end{tabular}
\end{center}
\caption{Advection field \eqref{eqn:etacos} (left) and the solution $w(T, \bx, \balpha)$
of \eqref{eqn:advdiff}--\eqref{eqn:advdiffbcic} (right) corresponding to a random realization of $\balpha$ from $\cB$
and $\nu = 0.01$.}
\label{fig:advdiff}
\end{figure}

Similarly to the experiments in Section~\ref{sec:linpde}, we compute the snapshots 
$\bphi_k = \bu(t_k, \balpha)$ by time-stepping \eqref{eqn:systemuad}
at $t_k = (1/30) k$, $k = 1,2,\ldots,N$, with $N = 30$ time steps and $T = 1$ using
Crank-Nicolson scheme. Then, the physical solution snapshots have the form \eqref{eqn:wsnap}.
The setting is illustrated in Figure~\ref{fig:advdiff} where we display the advection
field for a random realization of $\balpha \in \cA$ and the corresponding solution
$w(T, \bx, \balpha)$.

Projection ROM of \eqref{eqn:systemuad} is obtained similarly to that of \eqref{eqn:systemu}
using the matrix of reduced basis vectors $\rZ = [\bz_1, \ldots, \bz_n] \in \mathbb{R}^{M \times n}$.
Specifically,
\begin{equation}
\widetilde{\rM} \widetilde{\bu}_t + \left( \widetilde{\rK} + \widetilde{\rH}(\balpha) \right) 
\widetilde{\bu} = \widetilde{\blf},
\label{eqn:systemuadrom}
\end{equation}
where $\widetilde{\rM}$ and $\widetilde{\rK}$ are defined as in section~\ref{sec:linpde}, whereas
$\widetilde{\rH}(\balpha) = \rZ^T \rH(\balpha) \rZ \in \R^{n \times n}$, 
$\widetilde{\blf} = \rZ^T \blf \in \R^n$,
and the initial condition is $\widetilde{\bu}(0, \balpha) = \boldsymbol{0} \in \R^n$. 
Efficient evaluation of \eqref{eqn:systemuadrom} was discussed in Section~\ref{rem1}. 
Solving \eqref{eqn:systemuadrom} for $\widetilde{\bu}(t, \balpha)$ allows to compute the approximate 
solution snapshots $\widetilde{w}(t_k, \bx, \balpha)$ exactly as in \eqref{eqn:wromsnap}.

\subsubsection{Out-of-sample TROM performance: advection-diffusion equation}
\label{sec:oosadvdiff}

The testing of TROMs for out-of-sample parameters for the advection-diffusion system 
\eqref{eqn:advdiff}--\eqref{eqn:advdiffbcic} is performed similarly to that for the heat equation in
Section~\ref{sec:oosheat}. In particular, we compute the average over $\cA$ of the relative 
$L^\infty(0,T;L^2(\Omega))$ and $L^2(0,T;H^1(\Omega))$ errors for $\widetilde{w}$ and
use the statistical behavior of the TROM vs POD gain  
\eqref{eqn:relgain} to compare the three TROM variants to conventional POD-ROM. 

First, we test TROM in the following setting. The flow behavior is balanced between advection 
and diffusion, with a diffusion coefficient $\nu = 0.1$ and $M=1893$.  The parameter domain is 
$\cA =  [-0.05, 0.05]^8\times [0.1\pi, 0.3\pi]$ sampled at a Cartesian grid with 
$K = 3^8 \times 9 = 59,049$ points to obtain the sampling set $\hcA$.
We draw $N_r = 200$ out-of-sample realizations $\balpha^{(r)}$ from $\cA$ with each
$\alpha_i^{(r)}$, $i=1,\ldots,9$, distributed uniformly in its corresponding interval.
The averaged relative $L^\infty(0,T;L^2(\Omega))$ finite element error of HOSVD-TROM is evaluated as 
$N_r^{-1}\sum\limits_{r=1}^{N_r}R_{\rm HOSVD}(\balpha^{(r)})$ with $R_{\rm HOSVD}(\balpha)$ 
defined in \eqref{eqn:Rx}. The average relative $L^2(0,T;H^1(\Omega))$ finite element error is computed 
in the same way after modifying  $R_{\rm HOSVD}(\balpha)$ accordingly. 
The TROM vs POD gain statistic was defined in \eqref{eqn:relgain}.

For such large values of $K$ and therefore large $\bPhi$, $\#(\bPhi)=$ 3.3534e+09, 
the algorithm we use for finding CP decomposition turns out to be the most memory-intensive and runs out of memory.
Thus, in what follows we only report the results for HOSVD- and TT-TROM approaches.

\begin{table}[h!]
\caption{Tensor compression  ranks, averaged relative error of the HOSVD-TROM finite element solutions,  statistics of the gain \eqref{eqn:relgain} for $\nu = 0.1$, $K = 3^8 \times 9 = 59,049$.}
\label{tab:relgainad1}
\small
\begin{center}
\begin{tabular}{|c|c|c|c|c|c|c|c|c|c|}
\hline\\[-2.2ex]
 & \multicolumn{3}{c|}{ $n = 10$, $\widetilde{\eps} = 10^{-5}$ } 
 & \multicolumn{3}{c|}{ $n = 12$, $\widetilde{\eps} = 10^{-6}$ } 
 & \multicolumn{3}{c|}{ $n = 13$, $\widetilde{\eps} = 10^{-7}$ } \\
\hline
HOSVD ranks
& \multicolumn{3}{c|}{\small [78,3,3,3,3,3, } 
& \multicolumn{3}{c|}{\small [116,3,3,3,3,3, } 
& \multicolumn{3}{c|}{\small [153,3,3,3,3,3, } \\
& \multicolumn{3}{c|}{\small 3,3,3,6,11] } 
& \multicolumn{3}{c|}{\small 3,3,3,7,12] } 
& \multicolumn{3}{c|}{\small 3,3,3,9,13] } \\
TT ranks
& \multicolumn{3}{c|}{\small [77,99,117,121,121, } 
& \multicolumn{3}{c|}{\small [116,162,204,219,218,} 
& \multicolumn{3}{c|}{\small [179,251,330,364,362, } \\
& \multicolumn{3}{c|}{\small 107,85,62,37,11] } 
& \multicolumn{3}{c|}{\small 186,139,91,50,12] } 
& \multicolumn{3}{c|}{\small 298,209,126,63,14] } \\
\hline
{\small TROM FE error}
&\multicolumn{3}{c|}{} &\multicolumn{3}{c|}{}& \multicolumn{3}{c|}{}\\
{\footnotesize $L^\infty(0,T;L^2(\Omega))$ }
& \multicolumn{3}{c|}{1.15e-03} 
& \multicolumn{3}{c|}{8.36e-04} 
& \multicolumn{3}{c|}{7.95e-04} \\
{\footnotesize $L^2(0,T;H^1(\Omega))$ }
& \multicolumn{3}{c|}{9.86e-04} 
& \multicolumn{3}{c|}{9.46e-04} 
& \multicolumn{3}{c|}{9.14e-04} \\\hline
$G_X$ & mean & std & min &  mean & std & min &  mean & std & min  \\
\hline
HOSVD & \bf ~9.17~ &~5.87~ & 3.11 &  \bf ~11.32~ & ~10.42~ & 1.16 & \bf ~10.69~ & ~9.06~ & 1.28  \\
TT & \bf 9.13 & 5.85 & 3.11 & \bf 11.40 & 10.61 & 1.17  & \bf 10.65 & 9.06 & 1.22 \\
\hline
\end{tabular}
\end{center}
\end{table} 

We present in Table~\ref{tab:relgainad1} the averaged relative error of the HOSVD-TROM finite element 
solutions (for TT-TROM the errors were very close and so are skipped) and the behavior of the gain statistics 
when tensor compression error $\widetilde{\eps}$ decreases,
while simultaneously increasing $n$ to be slightly less or equal to Tucker rank $\widetilde{N}$
for HOSVD or compression rank $\widetilde{r}_{D+1}$ for TT, respectively. We observe in 
Table~\ref{tab:relgainad1} a relatively weak dependence of the errors and  TROM vs POD gain mean 
on the choice of $\widetilde{\eps}$ and $n$, hence, we conclude that higher tensor compression error 
and smaller $n$ are more beneficial, since they correspond to higher compression factors and possible
faster run times for the offline stage of TROM algorithms.

For the second example we choose an advection-dominated flow with a smaller diffusion 
coefficient $\nu = 0.01$ and $M=4,797$. The parameter domain is 
$\cA = [-0.01, 0.01]^8\times[0.1\pi, 0.5\pi]$ sampled on a Cartesian grid with $K = 20 \times 2^8 = 5,120$ 
points to obtain $\hcA$. This gives the snapshot tensor with $\#({\bPhi}) = $ 1.2280e+09 entries.
We draw $N_r = 100$ out-of-sample realizations $\balpha^{(r)}$ from $\cA$ with each 
$\alpha_i^{(r)}$, $i=1,\ldots,9$, distributed uniformly in its corresponding interval.

\begin{table}[h!]
\caption{Tensor compression  ranks, number of elements passed to the online stage, 
averaged relative error of the HOSVD-TROM finite element solutions, and
statistics of relative gain \eqref{eqn:relgain} for $\nu = 0.01$, $K = 20\times2^8 = 5,120$.}
\label{tab:relgainad3}
\small
\begin{center}
$\widetilde{\eps} = 10^{-3}$\\
\begin{tabular}{|c|c|c|c|c|c|c|c|c|c|}
\hline
HOSVD ranks 
& \multicolumn{6}{l}{\small [76,2,2,2,2,2,2,2,2,11,12] } 
& \multicolumn{3}{l|}{\small $\#(\mbox{online}(\widetilde{\bPhi}))$=2568444} \\
TT ranks
& \multicolumn{6}{l}{\small [75,77,79,76,77,75,71,67,61,11] } 
& \multicolumn{3}{l|}{\small $\#(\mbox{online}(\widetilde{\bPhi}))$=100747} \\ \hline
$n$ & \multicolumn{3}{c|}{5} & \multicolumn{3}{c|}{8} & \multicolumn{3}{c|}{10} \\
\hline
{\small TROM FE error}
&\multicolumn{3}{c|}{} &\multicolumn{3}{c|}{}& \multicolumn{3}{c|}{}\\
{\footnotesize $L^\infty(0,T;L^2(\Omega))$ }
& \multicolumn{3}{c|}{3.45e-2} 
& \multicolumn{3}{c|}{7.07e-3} 
& \multicolumn{3}{c|}{5.58e-3} \\
{\footnotesize $L^2(0,T;H^1(\Omega))$ }
& \multicolumn{3}{c|}{5.59e-2} 
& \multicolumn{3}{c|}{1.10e-2} 
& \multicolumn{3}{c|}{5.86e-3} \\\hline
$G_X$ & mean & std & min &   mean & std & min &  mean & std & min  \\
\hline
HOSVD& \bf 6.95 & 1.10 & 5.41 & \bf 22.56 & 6.97 & 12.02  & \bf ~32.66 & ~24.70 & 10.00  \\
TT   & \bf 6.95 & 1.09 & 5.41 & \bf 22.54 & 6.94 & 11.97  & \bf ~31.83 & ~23.58 & 9.99  \\
\hline
\end{tabular}
\vskip0.05in
$\widetilde{\eps} = 10^{-5}$\\
\begin{tabular}{|c|c|c|c|c|c|c|c|c|c|}
\hline
HOSVD ranks
& \multicolumn{6}{l}{\small [184,2,2,2,2,2,2,2,2,15,18] } 
& \multicolumn{3}{l|}{\small $\#(\mbox{online}(\widetilde{\bPhi}))$=12718412} \\
TT ranks
& \multicolumn{6}{l}{\small [183,235,288,305,319,295,246,187,128,18] } 
& \multicolumn{3}{l|}{\small $\#(\mbox{online}(\widetilde{\bPhi}))$=1110964} \\
\hline
$n$ & \multicolumn{3}{c|}{5} & \multicolumn{3}{c|}{10} & \multicolumn{3}{c|}{15} \\
\hline
{\small TROM FE error}
&\multicolumn{3}{c|}{} &\multicolumn{3}{c|}{}& \multicolumn{3}{c|}{}\\
{\footnotesize $L^\infty(0,T;L^2(\Omega))$ }
& \multicolumn{3}{c|}{3.45e-2} 
& \multicolumn{3}{c|}{5.56e-3} 
& \multicolumn{3}{c|}{5.51e-3} \\
{\footnotesize $L^2(0,T;H^1(\Omega))$ }
& \multicolumn{3}{c|}{5.59e-2} 
& \multicolumn{3}{c|}{5.80e-3} 
& \multicolumn{3}{c|}{5.08e-3} \\
\hline
 $G_X$ & mean & std & min &    mean & std & min &   mean & std & min  \\
\hline
HOSVD & \bf 6.95 & 1.10 & 5.41  & \bf 33.34 & 25.90 & 10.01 & ~\bf 19.45 & ~16.33 & 5.60  \\
TT    & \bf 6.95 & 1.10 & 5.41  & \bf 33.34 & 25.90 & 10.01 & ~\bf 19.45 & ~16.33 & 5.60  \\
\hline
\end{tabular}
\vskip0.05in
$\widetilde{\eps} = 10^{-7}$\\
\begin{tabular}{|c|c|c|c|c|c|c|c|c|c|}
\hline
HOSVD ranks
& \multicolumn{6}{l}{\small [476,2,2,2,2,2,2,2,2,18,25] } 
& \multicolumn{3}{l|}{\small $\#(\mbox{online}(\widetilde{\bPhi}))$=54835592} \\
TT ranks
& \multicolumn{6}{l}{\small [528,618,753,858,866,725,520,341,211,25] } 
& \multicolumn{3}{l|}{\small $\#(\mbox{online}(\widetilde{\bPhi}))$=6975287} \\
\hline
$n$ & \multicolumn{3}{c|}{5} & \multicolumn{3}{c|}{10} & \multicolumn{3}{c|}{15} \\
\hline
{\small TROM FE error}
&\multicolumn{3}{c|}{} &\multicolumn{3}{c|}{}& \multicolumn{3}{c|}{}\\
{\footnotesize $L^\infty(0,T;L^2(\Omega))$ }
& \multicolumn{3}{c|}{3.44e-2} 
& \multicolumn{3}{c|}{5.56e-3} 
& \multicolumn{3}{c|}{5.51e-3} \\
{\footnotesize $L^2(0,T;H^1(\Omega))$ }
& \multicolumn{3}{c|}{5.59e-2} 
& \multicolumn{3}{c|}{5.80e-3} 
& \multicolumn{3}{c|}{5.08e-3} \\
\hline
 $G_X$ & mean & std & min &  mean & std & min &   mean & std & min  \\
\hline
HOSVD & \bf 6.95 & 1.10 & 5.41  & \bf 33.34 & 25.90 & 10.01  & ~\bf 19.45 & ~16.33 & 5.60  \\
TT    & \bf 6.95 & 1.10 & 5.41  & \bf 33.34 & 25.90 & 10.01  & ~\bf 19.45 & ~16.33 & 5.60  \\
\hline
\end{tabular}
\end{center}
\end{table} 

Table~\ref{tab:relgainad3} shows tensor compression ranks,  number of elements passed to the online stage, 
the averaged relative error of the HOSVD-TROM finite element solutions, the relative gain statistics for three different 
levels of  tensor compression error $\widetilde{\eps} = 10^{-3}, 10^{-5}, 10^{-7}$ with
three different reduced space dimensions $n$ for each case. We observe that it is possible to achieve
performance that is very close to the best one with $\widetilde{\eps}$ as large as $10^{-3}$,
provided $n$ is large enough. Similarly to the results for $\nu = 0.1$, we suggest that for the given problem, 
discretization and parameter sampling, the finite element error is dominated by the interpolation error of 
the TROM and using a tighter compression threshold \rev{or larger $n$} does not lead to more accurate ROM solutions. 
It seems beneficial to use low-accuracy tensor decompositions to save on both the computation and storage, 
while not losing much in terms of relative gain compared to more expensive options. As expected, the TT format 
becomes more cost-efficient for the higher parameter space dimension $D$.  

Overall, while the accuracy increase of HOSVD- and TT-TROM over POD-ROM is still substantial in
the advection-diffusion setting with $D=9$ parameters, it is smaller than the one for the heat equation
considered in Section~\ref{sec:linpde}. This is most probably caused by larger variability of the
snapshots with respect to parameter variation making the problem a good candidate for 
non-interpolatory TROM (not studied here).  
 
\section{Conclusions}

Summarizing the findings of the paper, the tensorial projection ROM for parametric dynamical systems 
builds on several new ideas:\\
(i) To approximately represent the set of observed snapshots, it uses low-rank tensor formats, 
rather than a truncated SVD of the snapshot matrix.   
The corresponding tensor decompositions  provide POD-type universal basis while preserving information 
about solution variation with respect to parameters.\\
(ii) This additional information is used to find a local (parameter-specific) ROM basis for any incoming parameter
that is not necessarily from the training/sampling set.\\
(iii) The local basis can be represented by its coordinates in the universal low-dimensional basis allowing an 
effective split of the ROM evaluation between the online and offline phases.

An interpolation procedure was suggested to extract the information about parameter dependence of the solutions, 
and thus of the ROM spaces, from the low-rank tensor decompositions. Online stage uses fast linear algebra with 
complexity depending only on the compression ranks. Non-interpolatory or hybrid approaches are also possible 
and in fact can produce even more accurate and robust TROMs. We will study these options elsewhere. 
For interpolatory TROMs, Theorem~\ref{Th1} proves an estimate on the representation power of the local ROM bases.
Numerical experiment with parameterized heat equation supported the estimate and illustrated the role of each of its terms.   

Three popular compressed tensor formats were considered to represent the low-rank tensor in the TROM. 
Of course, other low-rank tensor decompositions can be used within the general framework of TROM. 
Out of the three tested, we found HOSVD to be most user-friendly and cost-efficient provided either the 
dimension of the parameter space is not too large or no Cartesian structure is exploited in organizing the snapshots. 
Otherwise, TT-TROM provides necessary tools to handle higher-dimensional parameter spaces. 
We also observed that the accuracy of TROMs crucially depend on $n$, $\widetilde{\eps}$ and parameter 
domain sampling, but not as much on the particular low-rank tensor format employed.

Finally, for higher-dimensional parameter spaces a grid-based sampling of the parameter domain becomes 
prohibitively expensive in terms of offline computation costs. \rev{Significant offline costs also 
incur for problems with less smooth dependence of solution on parameters, which would require a denser sampling, and for  problems where each high-fidelity solve is expensive because of fine spatial or temporal resolution.} We see several ways to develop TROMs \rev{addressing these challenges}: 
(i)~use a sophisticated sampling, e.g., based on a greedy strategy, and organize the snapshots in 3D tensors, 
(ii)~to benefit from Cartesian structure and higher order tensor decompositions, apply a tensor completion 
method to find a low-rank representation of the snapshot tensor sampled at a few nodes of the parameter grid,
\rev{ and
(iii)~combine TROMs with compressed formats to represent high-fidelity snapshots.}
We leave these options for a future research.    

\section*{Acknowledgments}
 M.O. was supported in  part by the U.S. National Science Foundation under awards DMS-2011444 and DMS-1953535.
This material is based upon research supported in part by the U.S. Office of Naval Research 
under award number N00014-21-1-2370 to A.M. The authors thank Vladimir Druskin, Traian Iliescu, and Vladimir Kazeev for their comments on the first draft of this paper. 

\bibliographystyle{siam}
\bibliography{literatur}{}

\begin{thebibliography}{10}

\bibitem{amsallem2008interpolation}
{\sc D.~Amsallem and C.~Farhat}, {\em Interpolation method for adapting
  reduced-order models and application to aeroelasticity}, AIAA journal, 46
  (2008), pp.~1803--1813.

\bibitem{amsallem2012nonlinear}
{\sc D.~Amsallem, M.~J. Zahr, and C.~Farhat}, {\em Nonlinear model order
  reduction based on local reduced-order bases}, International Journal for
  Numerical Methods in Engineering, 92 (2012), pp.~891--916.

\bibitem{antoulas2000survey}
{\sc A.~C. Antoulas, D.~C. Sorensen, and S.~Gugercin}, {\em A survey of model
  reduction methods for large-scale systems}, tech. report, 2000.

\bibitem{ballani2017multilevel}
{\sc J.~Ballani, D.~Kressner, and M.~D. Peters}, {\em Multilevel tensor
  approximation of pdes with random data}, Stochastics and Partial Differential
  Equations: Analysis and Computations, 5 (2017), pp.~400--427.

\bibitem{baur2011interpolatory}
{\sc U.~Baur, C.~Beattie, P.~Benner, and S.~Gugercin}, {\em Interpolatory
  projection methods for parameterized model reduction}, SIAM Journal on
  Scientific Computing, 33 (2011), pp.~2489--2518.

\bibitem{bengua2017efficient}
{\sc J.~A. Bengua, H.~N. Phien, H.~D. Tuan, and M.~N. Do}, {\em Efficient
  tensor completion for color image and video recovery: Low-rank tensor train},
  IEEE Transactions on Image Processing, 26 (2017), pp.~2466--2479.

\bibitem{benner2016low}
{\sc P.~Benner, S.~Dolgov, A.~Onwunta, and M.~Stoll}, {\em Low-rank solvers for
  unsteady stokes--brinkman optimal control problem with random data}, Computer
  Methods in Applied Mechanics and Engineering, 304 (2016), pp.~26--54.

\bibitem{benner2017solving}
\leavevmode\vrule height 2pt depth -1.6pt width 23pt, {\em Solving optimal
  control problems governed by random navier-stokes equations using low-rank
  methods}, arXiv preprint arXiv:1703.06097,  (2017).

\bibitem{benner2014robust}
{\sc P.~Benner and L.~Feng}, {\em A robust algorithm for parametric model order
  reduction based on implicit moment matching}, in Reduced order methods for
  modeling and computational reduction, Springer, 2014, pp.~159--185.

\bibitem{benner2015survey}
{\sc P.~Benner, S.~Gugercin, and K.~Willcox}, {\em A survey of projection-based
  model reduction methods for parametric dynamical systems}, SIAM review, 57
  (2015), pp.~483--531.

\bibitem{benner2015low}
{\sc P.~Benner, A.~Onwunta, and M.~Stoll}, {\em Low-rank solution of unsteady
  diffusion equations with stochastic coefficients}, SIAM/ASA Journal on
  Uncertainty Quantification, 3 (2015), pp.~622--649.

\bibitem{BrennerScott}
{\sc S.~Brenner and L.~Scott}, {\em The Mathematical Theory of Finite Element
  Methods}, Springer, New York, second~ed., 2002.

\bibitem{brunton2016discovering}
{\sc S.~L. Brunton, J.~L. Proctor, and J.~N. Kutz}, {\em Discovering governing
  equations from data by sparse identification of nonlinear dynamical systems},
  Proceedings of the national academy of sciences, 113 (2016), pp.~3932--3937.

\bibitem{bui2008model}
{\sc T.~Bui-Thanh, K.~Willcox, and O.~Ghattas}, {\em Model reduction for
  large-scale systems with high-dimensional parametric input space}, SIAM
  Journal on Scientific Computing, 30 (2008), pp.~3270--3288.

\bibitem{carlberg2013gnat}
{\sc K.~Carlberg, C.~Farhat, J.~Cortial, and D.~Amsallem}, {\em The gnat method
  for nonlinear model reduction: effective implementation and application to
  computational fluid dynamics and turbulent flows}, Journal of Computational
  Physics, 242 (2013), pp.~623--647.

\bibitem{carroll1970analysis}
{\sc J.~D. Carroll and J.-J. Chang}, {\em Analysis of individual differences in
  multidimensional scaling via an n-way generalization of “eckart-young”
  decomposition}, Psychometrika, 35 (1970), pp.~283--319.

\bibitem{chaturantabut2010nonlinear}
{\sc S.~Chaturantabut and D.~C. Sorensen}, {\em Nonlinear model reduction via
  discrete empirical interpolation}, SIAM Journal on Scientific Computing, 32
  (2010), pp.~2737--2764.

\bibitem{chinesta2010recent}
{\sc F.~Chinesta, A.~Ammar, and E.~Cueto}, {\em Recent advances and new
  challenges in the use of the proper generalized decomposition for solving
  multidimensional models}, Archives of Computational methods in Engineering,
  17 (2010), pp.~327--350.

\bibitem{chinesta2013proper}
{\sc F.~Chinesta, R.~Keunings, and A.~Leygue}, {\em The proper generalized
  decomposition for advanced numerical simulations: a primer}, Springer Science
  \& Business Media, 2013.

\bibitem{cohen2011analytic}
{\sc A.~Cohen, R.~Devore, and C.~Schwab}, {\em Analytic regularity and
  polynomial approximation of parametric and stochastic elliptic pde's},
  Analysis and Applications, 9 (2011), pp.~11--47.

\bibitem{de2000multilinear}
{\sc L.~De~Lathauwer, B.~De~Moor, and J.~Vandewalle}, {\em A multilinear
  singular value decomposition}, SIAM journal on Matrix Analysis and
  Applications, 21 (2000), pp.~1253--1278.

\bibitem{de2008tensor}
{\sc V.~De~Silva and L.-H. Lim}, {\em Tensor rank and the ill-posedness of the
  best low-rank approximation problem}, SIAM Journal on Matrix Analysis and
  Applications, 30 (2008), pp.~1084--1127.

\bibitem{dolgov2018direct}
{\sc S.~V. Dolgov, V.~A. Kazeev, and B.~N. Khoromskij}, {\em Direct
  tensor-product solution of one-dimensional elliptic equations with
  parameter-dependent coefficients}, Mathematics and computers in simulation,
  145 (2018), pp.~136--155.

\bibitem{drohmann2012reduced}
{\sc M.~Drohmann, B.~Haasdonk, and M.~Ohlberger}, {\em Reduced basis
  approximation for nonlinear parametrized evolution equations based on
  empirical operator interpolation}, SIAM Journal on Scientific Computing, 34
  (2012), pp.~A937--A969.

\bibitem{eftang2011hp}
{\sc J.~L. Eftang, D.~J. Knezevic, and A.~T. Patera}, {\em An hp certified
  reduced basis method for parametrized parabolic partial differential
  equations}, Mathematical and Computer Modelling of Dynamical Systems, 17
  (2011), pp.~395--422.

\bibitem{eftang2010hp}
{\sc J.~L. Eftang, A.~T. Patera, and E.~M. R{\o}nquist}, {\em An" hp" certified
  reduced basis method for parametrized elliptic partial differential
  equations}, SIAM Journal on Scientific Computing, 32 (2010), pp.~3170--3200.

\bibitem{eigel2017adaptive}
{\sc M.~Eigel, M.~Pfeffer, and R.~Schneider}, {\em Adaptive stochastic galerkin
  fem with hierarchical tensor representations}, Numerische Mathematik, 136
  (2017), pp.~765--803.

\bibitem{gandy2011tensor}
{\sc S.~Gandy, B.~Recht, and I.~Yamada}, {\em Tensor completion and low-n-rank
  tensor recovery via convex optimization}, Inverse problems, 27 (2011),
  p.~025010.

\bibitem{grasedyck2013literature}
{\sc L.~Grasedyck, D.~Kressner, and C.~Tobler}, {\em A literature survey of
  low-rank tensor approximation techniques}, GAMM-Mitteilungen, 36 (2013),
  pp.~53--78.

\bibitem{grepl2007efficient}
{\sc M.~A. Grepl, Y.~Maday, N.~C. Nguyen, and A.~T. Patera}, {\em Efficient
  reduced-basis treatment of nonaffine and nonlinear partial differential
  equations}, ESAIM: Mathematical Modelling and Numerical Analysis, 41 (2007),
  pp.~575--605.

\bibitem{grepl2005posteriori}
{\sc M.~A. Grepl and A.~T. Patera}, {\em A posteriori error bounds for
  reduced-basis approximations of parametrized parabolic partial differential
  equations}, ESAIM: Mathematical Modelling and Numerical Analysis, 39 (2005),
  pp.~157--181.

\bibitem{griebel2021analysis}
{\sc M.~Griebel and H.~Harbrecht}, {\em Analysis of tensor approximation
  schemes for continuous functions}, Foundations of Computational Mathematics,
  (2021), pp.~1--22.

\bibitem{gugercin2004survey}
{\sc S.~Gugercin and A.~C. Antoulas}, {\em A survey of model reduction by
  balanced truncation and some new results}, International Journal of Control,
  77 (2004), pp.~748--766.

\bibitem{hackbusch2012tensor}
{\sc W.~Hackbusch}, {\em Tensor spaces and numerical tensor calculus}, vol.~42,
  Springer, 2012.

\bibitem{hackbusch2007tensor}
{\sc W.~Hackbusch and B.~N. Khoromskij}, {\em Tensor-product approximation to
  operators and functions in high dimensions}, Journal of Complexity, 23
  (2007), pp.~697--714.

\bibitem{harshman1970foundations}
{\sc R.~A. Harshman}, {\em Foundations of the parafac procedure: Models and
  conditions for an" explanatory" multimodal factor analysis},  (1970).

\bibitem{Hartman}
{\sc P.~Hartman}, {\em Ordinary Differential Equations}, vol.~590, John Wiley
  and Sons, 1964.

\bibitem{haastad1990tensor}
{\sc J.~H{\aa}stad}, {\em Tensor rank is np-complete}, Journal of Algorithms,
  11 (1990), pp.~644--654.

\bibitem{hesthaven2016certified}
{\sc J.~S. Hesthaven, G.~Rozza, and B.~Stamm}, {\em Certified reduced basis
  methods for parametrized partial differential equations}, vol.~590, Springer,
  2016.

\bibitem{hitchcock1927expression}
{\sc F.~L. Hitchcock}, {\em The expression of a tensor or a polyadic as a sum
  of products}, Journal of Mathematics and Physics, 6 (1927), pp.~164--189.

\bibitem{huang2014provable}
{\sc B.~Huang, C.~Mu, D.~Goldfarb, and J.~Wright}, {\em Provable low-rank
  tensor recovery}, Optimization-Online, 4252 (2014), pp.~455--500.

\bibitem{kastian2020two}
{\sc S.~Kastian, D.~Moser, L.~Grasedyck, and S.~Reese}, {\em A two-stage
  surrogate model for neo-hookean problems based on adaptive proper orthogonal
  decomposition and hierarchical tensor approximation}, Computer Methods in
  Applied Mechanics and Engineering, 372 (2020), p.~113368.

\bibitem{kerschen2005method}
{\sc G.~Kerschen, J.-c. Golinval, A.~F. Vakakis, and L.~A. Bergman}, {\em The
  method of proper orthogonal decomposition for dynamical characterization and
  order reduction of mechanical systems: an overview}, Nonlinear dynamics, 41
  (2005), pp.~147--169.

\bibitem{khoromskij2011tensor}
{\sc B.~N. Khoromskij and C.~Schwab}, {\em Tensor-structured galerkin
  approximation of parametric and stochastic elliptic pdes}, SIAM Journal on
  Scientific Computing, 33 (2011), pp.~364--385.

\bibitem{kiers2000towards}
{\sc H.~A. Kiers}, {\em Towards a standardized notation and terminology in
  multiway analysis}, Journal of Chemometrics: A Journal of the Chemometrics
  Society, 14 (2000), pp.~105--122.

\bibitem{ReviewTensor}
{\sc T.~G. Kolda and B.~W. Bader}, {\em Tensor decompositions and
  applications}, SIAM review, 51 (2009), pp.~455--500.

\bibitem{kramer2019nonlinear}
{\sc B.~Kramer and K.~E. Willcox}, {\em Nonlinear model order reduction via
  lifting transformations and proper orthogonal decomposition}, AIAA Journal,
  57 (2019), pp.~2297--2307.

\bibitem{kressner2011low}
{\sc D.~Kressner and C.~Tobler}, {\em Low-rank tensor krylov subspace methods
  for parametrized linear systems}, SIAM Journal on Matrix Analysis and
  Applications, 32 (2011), pp.~1288--1316.

\bibitem{lee2019low}
{\sc K.~Lee, H.~C. Elman, and B.~Sousedik}, {\em A low-rank solver for the
  navier--stokes equations with uncertain viscosity}, SIAM/ASA Journal on
  Uncertainty Quantification, 7 (2019), pp.~1275--1300.

\bibitem{liang2002properi}
{\sc Y.~Liang, H.~Lee, S.~Lim, W.~Lin, K.~Lee, and C.~Wu}, {\em Proper
  orthogonal decomposition and its applications—part i: Theory}, Journal of
  Sound and vibration, 252 (2002), pp.~527--544.

\bibitem{liang2002proper}
\leavevmode\vrule height 2pt depth -1.6pt width 23pt, {\em Proper orthogonal
  decomposition and its applications—part i: Theory}, Journal of Sound and
  vibration, 252 (2002), pp.~527--544.

\bibitem{liang2002properii}
{\sc Y.~Liang, W.~Lin, H.~Lee, S.~Lim, K.~Lee, and H.~Sun}, {\em Proper
  orthogonal decomposition and its applications--part ii: Model reduction for
  mems dynamical analysis}, Journal of Sound and Vibration, 256 (2002),
  pp.~515--532.

\bibitem{liu2012tensor}
{\sc J.~Liu, P.~Musialski, P.~Wonka, and J.~Ye}, {\em Tensor completion for
  estimating missing values in visual data}, IEEE transactions on pattern
  analysis and machine intelligence, 35 (2012), pp.~208--220.

\bibitem{lumley1967structure}
{\sc J.~L. Lumley}, {\em The structure of inhomogeneous turbulent flows},
  Atmospheric turbulence and radio wave propagation,  (1967).

\bibitem{nouy2015low}
{\sc A.~Nouy}, {\em Low-rank tensor methods for model order reduction}, arXiv
  preprint arXiv:1511.01555,  (2015).

\bibitem{nouy2017low}
\leavevmode\vrule height 2pt depth -1.6pt width 23pt, {\em Low-rank methods for
  high-dimensional approximation and model order reduction}, Model reduction
  and approximation, P. Benner, A. Cohen, M. Ohlberger, and K. Willcox, eds.,
  SIAM, Philadelphia, PA,  (2017), pp.~171--226.

\bibitem{TT2}
{\sc I.~Oseledets and E.~Tyrtyshnikov}, {\em Tt-cross approximation for
  multidimensional arrays}, Linear Algebra and its Applications, 432 (2010),
  pp.~70--88.

\bibitem{TT1}
{\sc I.~V. Oseledets}, {\em Tensor-train decomposition}, SIAM Journal on
  Scientific Computing, 33 (2011), pp.~2295--2317.

\bibitem{patera2007reduced}
{\sc A.~T. Patera, G.~Rozza, et~al.}, {\em Reduced basis approximation and a
  posteriori error estimation for parametrized partial differential equations},
  2007.

\bibitem{rathinam2003new}
{\sc M.~Rathinam and L.~R. Petzold}, {\em A new look at proper orthogonal
  decomposition}, SIAM Journal on Numerical Analysis, 41 (2003),
  pp.~1893--1925.

\bibitem{rowley2005model}
{\sc C.~W. Rowley}, {\em Model reduction for fluids, using balanced proper
  orthogonal decomposition}, International Journal of Bifurcation and Chaos, 15
  (2005), pp.~997--1013.

\bibitem{schneider2014approximation}
{\sc R.~Schneider and A.~Uschmajew}, {\em Approximation rates for the
  hierarchical tensor format in periodic sobolev spaces}, Journal of
  Complexity, 30 (2014), pp.~56--71.

\bibitem{sidiropoulos2017tensor}
{\sc N.~D. Sidiropoulos, L.~De~Lathauwer, X.~Fu, K.~Huang, E.~E. Papalexakis,
  and C.~Faloutsos}, {\em Tensor decomposition for signal processing and
  machine learning}, IEEE Transactions on Signal Processing, 65 (2017),
  pp.~3551--3582.

\bibitem{sirovich1987turbulence}
{\sc L.~Sirovich}, {\em Turbulence and the dynamics of coherent structures. i.
  coherent structures}, Quarterly of applied mathematics, 45 (1987),
  pp.~561--571.

\bibitem{son2013real}
{\sc N.~T. Son}, {\em A real time procedure for affinely dependent parametric
  model order reduction using interpolation on grassmann manifolds},
  International Journal for Numerical Methods in Engineering, 93 (2013),
  pp.~818--833.

\bibitem{temlyakov1988estimates}
{\sc V.~N. Temlyakov}, {\em Estimates for the best bilinear approximations of
  periodic functions}, Trudy Matematicheskogo Instituta imeni VA Steklova, 181
  (1988), pp.~250--267.

\bibitem{trefethen2017multivariate}
{\sc L.~Trefethen}, {\em Multivariate polynomial approximation in the
  hypercube}, Proceedings of the American Mathematical Society, 145 (2017),
  pp.~4837--4844.

\bibitem{tucker1966some}
{\sc L.~R. Tucker}, {\em Some mathematical notes on three-mode factor
  analysis}, Psychometrika, 31 (1966), pp.~279--311.

\bibitem{yuan2016tensor}
{\sc M.~Yuan and C.-H. Zhang}, {\em On tensor completion via nuclear norm
  minimization}, Foundations of Computational Mathematics, 16 (2016),
  pp.~1031--1068.

\end{thebibliography}

\end{document}